\documentclass[twoside,a4paper]{amsart}

\makeatletter
\renewcommand{\tocsection}[3]{%
  \indentlabel{\@ifnotempty{#2}{\bfseries\ignorespaces#1 #2\quad}}\bfseries#3}
\renewcommand{\tocsubsection}[3]{%
  \indentlabel{\@ifnotempty{#2}{\ignorespaces#1 #2\quad}}#3}

\newcommand\@dotsep{4.5}
\def\@tocline#1#2#3#4#5#6#7{\relax
  \ifnum #1>\c@tocdepth 
  \else
    \par \addpenalty\@secpenalty\addvspace{#2}%
    \begingroup \hyphenpenalty\@M
    \@ifempty{#4}{%
      \@tempdima\csname r@tocindent\number#1\endcsname\relax
    }{%
      \@tempdima#4\relax
    }%
    \parindent\z@ \leftskip#3\relax \advance\leftskip\@tempdima\relax
    \rightskip\@pnumwidth plus1em \parfillskip-\@pnumwidth
    #5\leavevmode\hskip-\@tempdima{#6}\nobreak
    \leaders\hbox{$\m@th\mkern \@dotsep mu\hbox{.}\mkern \@dotsep mu$}\hfill
    \nobreak
    \hbox to\@pnumwidth{\@tocpagenum{\ifnum#1=1\bfseries\fi#7}}\par
    \nobreak
    \endgroup
  \fi}
\AtBeginDocument{%
\expandafter\renewcommand\csname r@tocindent0\endcsname{0pt}
}
\def\l@subsection{\@tocline{2}{0pt}{2.5pc}{5pc}{}}
\makeatother

\usepackage{xcolor}
\definecolor{webgreen}{rgb}{0,.5,0}
\definecolor{burgundy}{rgb}{.6,0,0}
\definecolor{RoyalBlue}{cmyk}{1, 0.50, 0, 0}
\usepackage[colorlinks=true, breaklinks=true, urlcolor=burgundy, linkcolor=RoyalBlue, citecolor=webgreen,backref=page]{hyperref}
\usepackage{tikz,epsfig, graphicx, subfig}
\usepackage{verbatim, setspace}
\usepackage{amsmath, amssymb}

\usepackage{newtxtext,newtxmath}

\usepackage{mathabx}

\numberwithin{equation}{section}

\newcommand{\R}		{\mathbb{R}}
\newcommand{\C}		{\mathbb{C}}
\newcommand{\N}		{\mathbb{N}}
\newcommand{\Z}		{\mathbb{Z}}

\newcommand{\diag}{\mathrm{diag}}

\newcommand{\im}{\mathsf{Im}}

\renewcommand{\arg}{\mathrm{arg}}
\renewcommand{\det}{\mathrm{det}}

\newcommand{\qasq}{\quad \text{as} \quad}
\newcommand{\qandq}{\quad \text{and} \quad}

\newcommand{\dd}{\mathrm{d}}
\newcommand{\ic}{\mathrm{i}}
\newcommand{\RS}{\mathfrak S}
\newcommand{\bd}{\boldsymbol{\Delta}}
\newcommand{\ualpha}{\boldsymbol\alpha}
\newcommand{\ubeta}{\boldsymbol\beta}

\newcommand{\z}	{{\boldsymbol z}}
\newcommand{\x}	{{\boldsymbol x}}
\newcommand{\tr} {{\boldsymbol t}}
\newcommand{\s}	{{\boldsymbol s}}

\newcommand{\rhy}   {\textnormal{RHP}-${\boldsymbol Y}$}
\newcommand{\rhx}   {\textnormal{RHP}-${\boldsymbol X}$}
\newcommand{\rhn}   {\textnormal{RHP}-${\boldsymbol N}$}
\newcommand{\rhz}   {\textnormal{RHP}-${\boldsymbol Z}$}
\newcommand{\rhpi}   {\textnormal{RHP}-${\boldsymbol P_i}$}
\newcommand{\rhpo}   {\textnormal{RHP}-${\boldsymbol P_0}$}

\newcommand{\rhphi}   {\textnormal{RHP}-${\boldsymbol \Phi}_\alpha$}
\newcommand{\rhpsi}   {\textnormal{RHP}-${\boldsymbol \Psi}_\nu$}
\newcommand{\rhphis}   {\textnormal{RHP}-${\boldsymbol \Phi}_\alpha^\mathsf{sing}$}

\newtheorem{theorem}{Theorem}
\newtheorem{proposition}{Proposition}
\newtheorem*{definition}{Definition}
\newtheorem*{bvproblem}{Boundary Value Problem: BVP-\(\Psi\)}

\begin{document}

\title[]{Non-Hermitian Orthogonal Polynomials on a Trefoil}

\author[A. Barhoumi]{Ahmad B. Barhoumi}

\address{Department of Mathematics, University of Michigan, 530 Church Street, Ann Arbor, MI 48109, USA}
\email{\href{mailto:barhoumi@umich.edu}{barhoumi@umich.edu}}

\author[M. Yattselev]{Maxim L. Yattselev}
\address{Department of Mathematical Sciences, Indiana University-Purdue University Indianapolis, 402~North Blackford Street, Indianapolis, IN 46202, USA}
\email{\href{mailto:maxyatts@iu.edu}{maxyatts@iu.edu}}

\begin{abstract}
We investigate asymptotic behavior of polynomials \( Q_n(z) \) satisfying non-Hermitian orthogonality relations
\[
\int_\Delta s^kQ_n(s)\rho(s)\dd s =0, \quad k\in\{0,\ldots,n-1\},
\]
where \( \Delta \) is a Chebotar\"ev (minimal capacity) contour connecting three non-collinear points and \( \rho(s) \) is a Jacobi-type weight including a possible power-type singularity at the Chebotar\"ev center of \( \Delta \). 
\end{abstract}



\thanks{The first author was partially supported by the National Science Foundation under DMS-1812625. The research of the second author was supported in part by a grant from the Simons Foundation, CGM-706591.}

\maketitle

\tableofcontents

\section{Introduction}

Let \( f(z)=\sum_{k=1}^\infty f_kz^{-k} \) be a convergent power series. The \( n \)-th diagonal Pad\'e approximant of \( f(z) \) is a rational function \( [n/n]_f(z) = P_n(z)/Q_n(z) \), where \( \deg(P_n),\deg(Q_n)\leq n\) and \( Q_n(z) \) is not identically zero, such that the linearized error function \( R_n(z) \) satisfies
\begin{equation}
\label{linear}
R_n(z) := (Q_nf-P_n)(z) = \mathcal O\left(z^{-(n+1)}\right)
\end{equation}
as \( z\to\infty \). One can readily check that the above condition is nothing but a system of linear equations on the coefficients of the polynomials of \( P_n(z),Q_n(z) \) that always has a non-trivial solution. This system is not necessarily uniquely solvable, but it is known that the rational function \( [n/n]_f(z) \) is indeed unique. Hereafter, we shall understand that \( Q_n(z) \) in \eqref{linear} is monic and of minimal possible degree, which does make it unique.

Our goal is to understand convergence properties of \( \smash{[n/n]_f(z)} \). It was shown by Herbert Stahl \cite{St85,St85b,St97} that if \( f(z) \) can be meromorphically continued along any path \( \overline\C\setminus A \) for a polar set \( A \) and there exists a point in \(\C\setminus A \) with at least two distinct continuations, then there exists a compact set \( \Delta_f \) such that  \( f(z) \) has a single-valued meromorphic continuation into the complement of \( \Delta_f \) (which is a branch cut for \( f(z) \)) and the diagonal Pad\'e approximants converge to this continuation in logarithmic capacity. The set \( \Delta_f \) is uniquely characterized as the branch cut of smallest logarithmic capacity that is set-theoretically minimal (if \( K \) is a branch cut with the same logarithmic capacity, then \( \Delta_f\subseteq K \)).

The minimal capacity contour \( \Delta_f \) can also be characterized from the point of view of quadratic differentials. Assume for simplicity that \( A=\{a_1,\ldots,a_m\} \), \( m\geq2 \), is a finite set. Suppose that every element of \( A \) is a branch point of \( f(z) \). Then there exist auxiliary points \( b_i \), \( i\in\{1,\ldots,m-2\} \), sometimes called Chebotar\"ev centers\footnote{In a somewhat different language, Chebotar\"ev posed a problem of finding a connected set of minimal logarithmic capacity containing a given finite set of points; descriptions of this set were independently given by Gr\"otzsch \cite{Grot30} and Lavrentiev \cite{Lav30,Lav34}.}, which are not necessarily distinct nor disjoint from the elements of \( A \), such that the set \( \Delta_f \) consists of the critical  trajectories of a rational quadratic differential
\[
\frac{\prod_{i=1}^{m-2}(z-b_i)}{\prod_{i=1}^m(z-a_i)}\dd z^2.
\]

Suppose that all the points \( b_i \) are disjoint from \( A \) and that each \( b_i \) appears either once or an even number of times and in the latter case does not belong to \( \Delta_f \) (a generic situation). Assume further that \( f(z) \) has either logarithmic or power branching at each \( a_i \) (i.e., behaves like \( \log(z-a_i) \) or \( (z-a_i)^{\alpha_i} \), \( \alpha_i>-1 \)). Then the strong asymptotics of the corresponding diagonal Pad\'e approximants was investigated by Aptekarev and the second author in \cite{ApYa15}. Our overarching goal is to remove all the assumptions on contours \( \Delta_f \) corresponding to finite sets \( A \). As will become clear later, the main difficulty lies in the local analysis of the polynomials \( Q_n(z) \) around the auxiliary points \( b_i \). The first step in this direction was taken by the authors in \cite{BarYa20} where we considered the case \( m=4 \) and \( b_1=b_2\in\Delta_f\setminus A \). Here, we consider that case \( m=4 \) and \( b_1=b_2\in A \).

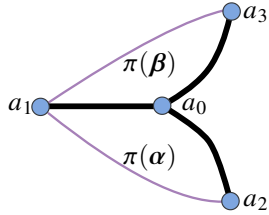
\begin{figure}[!ht]
\includegraphics[scale=.6]{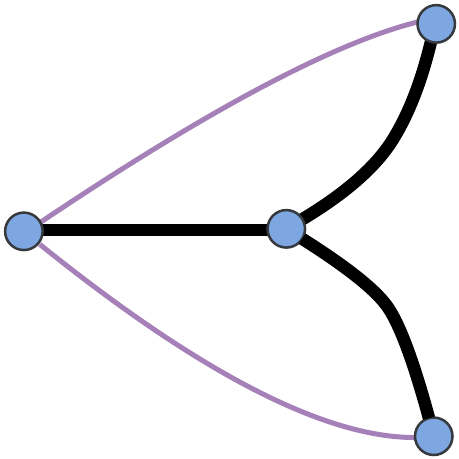}
\begin{picture}(0,0)
\put(-25,38){\(a_0\)}
\put(-90,38){\(a_1\)}
\put(-2,73){\(a_3\)}
\put(-2,2){\(a_2\)}
\put(-47,18){\(\pi(\ualpha)\)}
\put(-47,53){\(\pi(\ubeta)\)}
\end{picture}
\caption{Contour \( \Delta \).}
\label{fig:1}
\end{figure}

Let us now change the notation slightly. Given three distinct non-collinear points \( a_1,a_2,a_3 \), denote by \( a_0 \) their Chebotar\"ev center. That is, there are three disjoint, except for \( a_0 \), analytic arcs, say \( \Delta_1,\Delta_2,\Delta_3 \) (\( \Delta_i \) has endpoints \( a_i \) and \( a_0 \)) such that
\begin{equation}
\label{differential}
\frac{z(t)-a_0}{(z(t)-a_1)(z(t)-a_2)(z(t)-a_3)}\big(z^\prime(t)\big)^2<0
\end{equation}
for any smooth parametrization \( z(t) \) of any of the arcs \( \Delta_i \). Let \( \Delta:=\Delta_1\cup\Delta_2\cup\Delta_3 \) (in our preceding notation \( A=\{a_0,a_1,a_2,a_3\} \), \( b_1=b_2=a_0 \), and \( \Delta_f=\Delta \)), see Figure~\ref{fig:1}. Examples of functions \( f(z) \) that lead to such minimal capacity contours include
\begin{equation}
\label{log-power}
\prod(z-a_i)^{\alpha_i} \qandq \sum c_i\log(z-a_i),
\end{equation}
where \( \alpha_i \not\in \Z \) while \( \sum_{i=0}^4\alpha_i\in\Z \) and \( c_i\neq0 \) while \( \sum_{i=0}^4 c_i=0 \). Let us point out that if \( f(z) \) was given by either of the expressions above, but with \( \alpha_0=c_0 =0 \), then the asymptotics of the diagonal Pad\'e approximants to \( f(z) \) was obtained in \cite{ApYa15}, see also \cite{M-FRakhSuet12}. Moreover, asymptotics of the approximants for the case \( \alpha_0=\alpha_1=\alpha_2=\alpha_3=-1/2 \) is contained \cite{Ya15}, see also \cite{BYa13,Suet00}. {The reader might want to consult the Appendix~\ref{ap:functions} to see how these functions relate to the definition below.}

In this work, we shall consider the following class of functions. Orient each arc \( \Delta_i \) towards \( a_0 \). Assume that the points \( a_1,a_2,a_3 \) are labeled counter-clockwise around \( a_0 \). 

\begin{definition}
We are interested in functions \( f(z) \) of the form
\begin{equation}
\label{CI}
f(z) = \frac1{2\pi\ic}\int_\Delta\frac{\rho(s)\dd s}{s-z},
\end{equation}
where there exist exponents \( \alpha_i>-1 \), \( i\in\{0,1,2,3\} \), and branches of \( (z-a_i)^{\alpha_i} \) holomorphic across \( \Delta\setminus\{a_i\} \) for which the restriction \( \rho_i(s) \) of \( \rho(s) \) to \( \Delta_i^\circ:=\Delta_i\setminus\{a_0,a_i\} \) is such that \( \rho_i(s)(s-a_0)^{-\alpha_0}(s-a_i)^{-\alpha_i} \) extends to a holomorphic and non-vanishing function in some neighborhood of \( \Delta_i \). 
\end{definition}

Notice that \( \rho(s) \) is not defined at \( a_0 \) even if \( \alpha_0=0 \). However, in the latter case the values \( \rho_i(a_0) \) are well-defined. It was assumed in \cite{ApYa15} that \( \rho_1(z)+\rho_2(z)+\rho_3(z)\equiv0 \) in some neighborhood of \( a_0 \) (no branching assumption). No such supposition is made here, neither in some neighborhood of \( a_0 \) nor at \( a_0 \) itself. The main advantage of the class of functions introduced in \eqref{CI} is that the denominator polynomials \( Q_n(z) \) can be equivalently characterized as non-Hermitian orthogonal polynomials satisfying
\begin{equation}
\label{ortho}
\int_\Delta s^kQ_n(s)\rho(s)\dd s=0, \quad k\in\{0,\ldots,n-1\}.
\end{equation}

This paper is organized as follows. The next section is devoted to the description of the main term of the asymptotics of the polynomials \( Q_n(z) \). The functions constructed in that section are known in integrable systems literature as Baker-Akhiezer functions and in the literature on non-Hermitian orthogonal polynomials are sometimes called Nuttall-Szeg\H{o} functions. The propositions stated in the next section are proven in Section~\ref{s:NS}. Our main results on the asymptotics of the polynomials \( Q_n(z) \) are stated in Section~\ref{sec:main}. Their proofs are given in Sections~\ref{s:RH1} and~\ref{s:RH2}. As it happens, local behavior of the polynomials \( Q_n(z) \) around \( a_0 \) is described by certain special functions that come from a \( 2\times 2 \) matrix function solving the so-called Painvlev\'e XXXIV Riemann-Hilbert problem with Stokes parameters that depend on the weight \( \rho(s) \) in a transcendental way. This connection is described in Appendices~\ref{sec:painleve} and~\ref{ap:34}. Appendix~\ref{ap:example} illustrates by an example that on the rotationally symmetric contour \( \Delta \) asymptotics of \( Q_n(z) \) can be very different for certain classes of weights from the rest of the cases.

\section{Nuttall-Szeg\H{o} Functions}
\label{sec:NS-BA}

Similarly to orthogonal polynomials on an interval, the asymptotics of the polynomials \( Q_n(z) \) is described by the term that captures the geometric rate of their growth, see \eqref{Phi}, and a Szeg\H{o} function of the weight \( \rho (s)\), see \eqref{szego}. Both of these functions naturally live on a genus one Riemann surface associated with the contour \( \Delta \), see \eqref{surface}. However, the nature of the theory of functions on a Riemann surface leads to the introduction of the third term, in fact, a normal family of functions that depends on \( n \), which are essentially ratios of the Riemann theta functions, see \eqref{Ti} (the necessity of these terms was already demonstrated by Akhiezer \cite{Akh60}).

\begin{figure}[!ht]
\includegraphics[scale=.6]{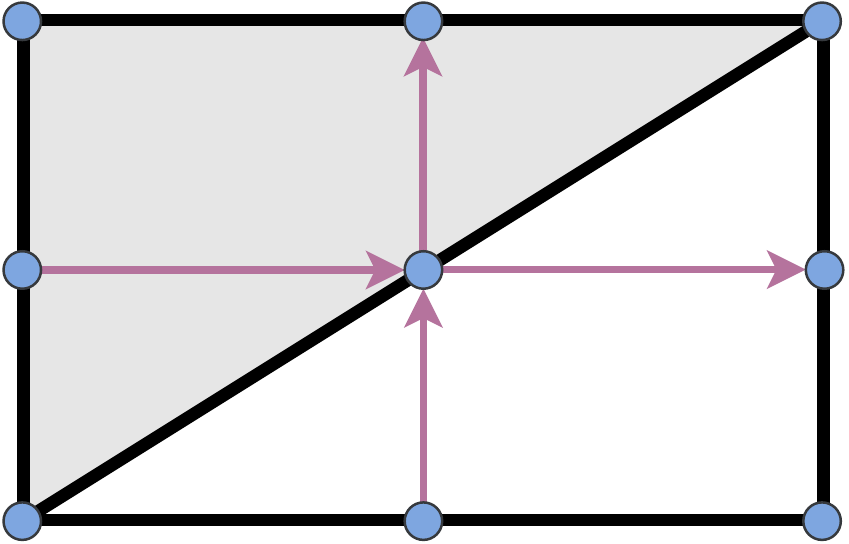}
\begin{picture}(0,0)
\put(-2,2){\(a_0\)}
\put(-2,89){\(a_0\)}
\put(-160,2){\(a_0\)}
\put(-160,89){\(a_0\)}
\put(-72,40){\(a_1\)}
\put(-2,45){\(a_2\)}
\put(-160,45){\(a_2\)}
\put(-79,-5){\(a_3\)}
\put(-79,96){\(a_3\)}
\put(-73,20){\(\ubeta\)}
\put(-85,68){\(\ubeta\)}
\put(-43,50){\(\ualpha\)}
\put(-113,50){\(\ualpha\)}
\put(-133,73){\(\RS^{(0)}\)}
\put(-33,13){\(\RS^{(1)}\)}
\put(-33,67){\(\bd_1\)}
\put(-133,22){\(\bd_1\)}
\put(-4,67){\(\bd_2\)}
\put(-160,22){\(\bd_2\)}
\put(-113,8){\(\bd_3\)}
\put(-43,82){\(\bd_3\)}
\end{picture}
\caption{Surface \( \RS \).}
\label{fig:torus}
\end{figure}

\subsection{Riemann surface}
\label{ssec:RS}

Let \( \RS \) be a genus one two-sheeted Riemann surface defined by
\begin{equation}
\label{surface}
\RS := \left\{\z=(z,w):~w^2=(z-a_0)(z-a_1)(z-a_2)(z-a_3)\right\}.
\end{equation}
We realize $\RS$ as a ramified cover of $\overline\C$ constructed in the following manner. Two copies of the extended complex plane are cut along each arc $\Delta_i^\circ$. These copies are glued together along the cuts in such a manner that the right (resp. left) side of the arc $\Delta_i^\circ$ belonging to one copy is joined with the left (resp. right) side of the same arc $\Delta_i^\circ$ only belonging to the other copy. Let
\[
\pi,w:\RS\to\overline\C, \quad \pi(\z) \mapsto z \qandq w(\z) \mapsto w.
\]
We call $\pi$ the canonical projection and label the sheets \( \RS^{(0)} \) and \( \RS^{(1)} \) (cut copies of the extended complex plane) so that \( w(\z) = w \) behaves like \( z^2 \) as \( \z\to\infty^{(0)} \). In our definition, the sheets are topologically closed and their intersection is equal to \( \bd \), where
\[
\bd:=\pi^{-1}(\Delta), \quad \bd_i:=\pi^{-1}(\Delta_i), \qandq \boldsymbol A :=\{ \boldsymbol a_0,\boldsymbol a_1,\boldsymbol a_2,\boldsymbol a_3\}, \quad \pi(\boldsymbol a_l)=a_l,
\]
\( i\in \{1,2,3\} \) and \( l\in\{0,1,2,3\} \). Thus, \( \boldsymbol A \) is the set of the ramification points of \( \RS \) and \( \bd \) is a union of three simple cycles that all intersect each other at \( \boldsymbol a_0\). We orient each \( \bd_i \) so that \( \RS^{(0)} \) remains on the left when $\bd_i$ is traversed in the positive direction. We write
\[
 z^{(k)} := \pi^{-1}(z)\cap\RS^{(k)}, \quad z\in\overline\C\setminus\Delta,
\]
where \( k\in\{0,1\} \), and use bold lower case letters such as $\z,\tr,\s$ to indicate generic points on $\RS$ with canonical projections $z,t,s$. We designate the symbol $\cdot^*$ to stand for the conformal involution that sends $z^{(k)}$ into $z^{(1-k)}$, $k\in\{0,1\}$ (we then extend this notion to \( \bd \) by continuity). Since \( \RS \) has genus \( 1 \), we need to select a homology basis. We choose cycles \( \ualpha,\ubeta \) to be involution-symmetric, i.e., \( \s\in\ualpha \) if and only if \( \s^*\in\ualpha\), and such that \( \pi(\ualpha),\pi(\ubeta) \) are smooth Jordan arcs joining \( a_1,a_2 \) and \( a_1,a_3 \) respectively, while \( \ualpha \) is oriented towards \( \boldsymbol a_1 \) and \( \ubeta \) away from \( \boldsymbol a_1 \) in \( \RS^{(0)} \), see Figures~\ref{fig:1} and~\ref{fig:torus}. In what follows it will be sometimes convenient to set
\[
\RS_{\ualpha} := \RS\setminus\ualpha \qandq \RS_{\ualpha,\ubeta} := \RS\setminus(\ualpha\cup\ubeta).
\]
Finally, given a function \( F(\z) \) defined on \( \RS\setminus\bd \), or on its subset, we shall set \( F(z):=F(z^{(0)}) \) and \( F^*(z) := F(z^{(1)}) \). That is, with a slight abuse of notation we use the same letter for a function on the top sheet as well as its pull-back to the complex plane while we use the involution symbol \( ^* \) to mark the pull-back to the complex plane from the bottom sheet. In particular,
\begin{equation}
\label{w}
w(z) = \sqrt{(z-a_0)(z-a_1)(z-a_2)(z-a_3)}, \quad z\in\C\setminus\Delta,
\end{equation}
is the branch such that \( w(z) = z^2 + \mathcal O(z) \) as \( z\to\infty \) and \( w^*(z)=-w(z) \). 

As it happens, the asymptotics of the polynomials \( Q_n(z) \) as well as of the linearized error functions \( R_n(z) \), see \eqref{linear},  is described by the solution of the following boundary value problem on \( \bd \).

\begin{bvproblem}
\label{bvp}
Given \( n\in \N \) and a function \( \rho(s) \) on \( \Delta \) as described after \eqref{CI}, find a function \( \Psi(\z) \) that is holomorphic in \( \RS\setminus\big(\bd\cup\big\{\infty^{(0)}\big\}\big) \), has a pole of order \( n \) at \( \infty^{(0)} \) and a zero of order at least \( n-1 \) at \( \infty^{(1)} \), and whose traces on \( \bd \) are continuous away from \( \boldsymbol A \) and satisfy
\begin{equation}
\label{BVP-Psi}
\left\{
\begin{array}{l}
\Psi_+(\s) =\Psi_-(\s)/(\rho w_+)(s), \quad \s\in\bd, \medskip \\
\big|\Psi(z)\big| = \mathcal O\left(|z-a_l|^{-(2\alpha_l+1)/4}\right) \qasq z\to a_l, \;\; l\in\{0,1,2,3\}.
\end{array}
\right.
\end{equation}
\end{bvproblem}

The solution of this boundary value problem is given in Theorem~\ref{thm:BANS} further below. It is rather explicit and constructed as a product of three terms that are introduced in the following three subsections.

\subsection{Geometric Term}
\label{ssec:map}

Let
\begin{equation}
\label{Phi}
\Phi(\z):= \exp\left\{\int_{\boldsymbol a_0}^\z G\right\}, \quad G(\z) := \frac{(z-a_0)\dd z}{w(\z)},
\end{equation}
for \( \z\in\RS_{\ualpha,\ubeta} \), where the path of integration lies entirely in \( \RS_{\ualpha,\ubeta} \) and \( G(\z) \) is a meromorphic differential on \( \RS \) having two simple poles at \( \infty^{(1)} \) and \( \infty^{(0)} \) with respective residues \( 1 \) and \( -1 \), whose period over any cycle on \( \RS \) is purely imaginary (the latter claim follows from \eqref{differential}). The function $\Phi(\z)$ is holomorphic and non-vanishing on $\RS_{\ualpha,\ubeta}$ except for a simple pole at $\infty^{(0)}$ and a simple zero at $\infty^{(1)}$.  Define (real) constants
\begin{equation}
\label{constants}
\omega:=-\frac1{2\pi\ic}\oint_{\ubeta}G \qandq \tau:=\frac1{2\pi\ic}\oint_{\ualpha}G.
\end{equation}
One can readily check that \( \Phi(\z) \) possesses continuous traces on both sides of each cycle of the canonical basis (in fact, it extends to a multiplicatively multi-valued function on the entire surface) that satisfy
\begin{equation}
\label{Phi-jump}
\Phi_+(\s) = \Phi_-(\s)\left\{
\begin{array}{ll}
e^{2\pi\ic \omega}, & \s \in \ualpha\setminus\{\boldsymbol a_1\}, \smallskip \\
e^{2\pi\ic \tau}, & \s \in \ubeta\setminus\{\boldsymbol a_1\}.
\end{array}
\right.
\end{equation}
Since \( w(z)=-w^*(z) \), it follows from \eqref{Phi} that \( \Phi(z)\Phi^*(z)\equiv 1 \). Moreover, since the constants in \eqref{constants} are real, it also holds that \( \log|\Phi(\z)| \) is harmonic on \( \RS\setminus\{\infty^{(0)},\infty^{(1)}\} \). Furthermore, since \( \Delta \) is the set of critical trajectories of the quadratic differential \( -(z-a_0)^2w^{-2}(z)\dd z^2 \) and the domain \( \overline\C\setminus\Delta \) contains a single critical point of this differential, namely, double pole at infinity, it follows from Basic Structure Theorem \cite[Theorem~3.5]{Jenkins} that \( \overline\C\setminus\Delta \) is a circle domain  for this differential, see \cite[Definition~3.9]{Jenkins}. In particular, it must be true that
\begin{equation}
\label{Phi-Green}
-\log|\Phi^*(z)| = \log|\Phi(z)|>0, \quad z\in \overline\C\setminus\Delta.
\end{equation}
In particular, this means that \( \log|\Phi(s)|\equiv0 \), \( s\in\Delta \), and that \( \log|\Phi(z)| \) is the Green function for \( \overline\C\setminus\Delta \) with pole at infinity.

\subsection{Szeg\H{o} Function}
\label{ssec:SF}

Since \( \RS \) has genus one, it has one (up to a constant factor) holomorphic differential. Since \( \RS \) is a two-sheeted surface, it is a standard fact of the theory of Riemann surfaces that
\begin{equation}
\label{B}
\Omega(\z) := \left(\oint_{\ualpha} \frac{\dd s}{w(\s)}\right)^{-1}\frac{\dd z}{w(\z)}, \quad \mathsf B := \oint_{\ubeta}\Omega,
\end{equation}
is the holomorphic differential on \( \RS \) normalized to have unit period on the \( \ualpha \) cycle of the homology basis (it was shown by Riemann that the constant \( \mathsf B \) has positive imaginary part). In particular, the function
\begin{equation}
\label{a}
a(\z) := \int_{\boldsymbol a_0}^\z \Omega, \quad \z\in\RS_{\ualpha,\ubeta},
\end{equation}
where the path of integration lies entirely in \( \RS_{\ualpha,\ubeta} \), extends to additively multi-valued analytic function on the entire surface \( \RS \). Denote further by \( \Omega_{\z,\z^*}(\s) \) the meromorphic differential with two simple poles at $\z$ and $\z^*$ with respective residues \( 1 \) and $-1$ normalized to have a zero period on the $\ualpha$-cycle. When \( \z \) does not lie on top of the point at infinity, it can be readily verified that
\begin{equation}
\label{cauchy-kernel}
\Omega_{\z,\z^*}(\s) = \dfrac{w(\boldsymbol z)}{s - z} \dfrac{\dd s}{w(\s)} - \left(\oint_{\ualpha} \dfrac{w(\z)}{t - z}\dfrac{\dd t}{w(\tr)} \right)\Omega(\s).
\end{equation}
The third kind differential \( \Omega_{\z,\z^*}(\s) \) can be thought of as a discontinuous Cauchy kernel on \( \RS \), see \cite{Zver71}. The following proposition elaborates this point.

\begin{proposition}
\label{prop:szego}
For a function \( \rho(s) \) as described after \eqref{CI} fix a branch of \( \log(\rho_iw_+)(s) \) that is continuous on \( \Delta_i^\circ \), \( i\in\{1,2,3\} \). Put
\begin{equation}
\label{szego}
S_\rho(\z) := \exp\left\{-\dfrac{1}{4 \pi \ic}\oint_{\bd}\log(\rho w_+)(s)\Omega_{\z,\z^*}(\s)\right\}.
\end{equation}
Then \( S_\rho(\z) \) is a holomorphic and non-vanishing function in \( \RS_{\ualpha}\setminus\bd \) with continuous traces on \( (\bd\cup\ualpha)\setminus\boldsymbol A \) that satisfy
\begin{equation}
\label{S-jump}
S_{\rho+}(\s) = S_{\rho-}(\s)\left\{
\begin{array}{rl}
\displaystyle \exp\big\{2\pi\ic c_\rho\big\}, &  \s\in\ualpha\setminus\boldsymbol A, \medskip \\
1/(\rho w_+)(s), & \s\in\bd\setminus\boldsymbol A,
\end{array}
\right.
\end{equation}
where the constant \( c_\rho \) is given by
\begin{equation}
\label{veccrho}
c_\rho := \frac1{2\pi\ic}\oint_{\bd}\log(\rho w_+)\Omega.
\end{equation}
It also holds that \( S_\rho(\z)S_\rho(\z^*)\equiv1 \), \( \z\in \RS_{\ualpha}\setminus\bd \), and
\begin{equation}
\label{Srho-ai}
\big|S_\rho(z)\big| \sim |z-a_l|^{-(2\alpha_l+1)/4}  \qasq z\to a_l, \quad l\in\{0,1,2,3\}.
\end{equation}
\end{proposition}

We call \( S_\rho(\z) \) the Szeg\H{o} function of the weight \( \rho(s) \). The Cauchy kernel \( \Omega_{\z,\z^*}(\s) \) is called discontinuous due to the presence of the jump of \( S_\rho(\z) \) across \( \ualpha \). We prove Proposition~\ref{prop:szego} in Section~\ref{ssec:szego}.

\subsection{Adjustment Terms}

Due to the presence of the jumps of both \( \Phi(\z) \) and \( S_\rho(\z) \) across the cycles of the homology basis, we need to introduce an adjustment term that will cancel out these jumps. To this end, let
\[
\theta(u) := \sum_{n\in\Z}\exp\big\{\pi\ic \mathsf Bn^2 + 2\pi\ic nu\big\}, \quad u\in\C,
\]
be the Riemann theta function associated with $\mathsf B$. Further, let \( \mathsf{Jac}(\RS) := \C/ \{\Z+\mathsf B\Z \} \) be the Jacobi variety of \( \RS \), where \( \mathsf B \) is given by \eqref{B}. We shall represent elements of \( \mathsf{Jac}(\RS) \) as equivalence classes \( [s\,] = \{s+ l + \mathsf Bm:l,m\in \Z\} \), where \( s\in\C \). It is known that Abel's map
\begin{equation}
\label{Amap}
\z\in\RS \mapsto \mathfrak a(\z) \in \mathsf{Jac}(\RS), \quad \mathfrak a(\z) :=\left[\int_{\boldsymbol a_0}^\z\Omega\right],
\end{equation}
is a bijection (clearly, it holds that \( [a(\z)]=\mathfrak a(\z) \) for \( \z\in\RS_{\ualpha,\ubeta} \), see \eqref{a}). Inverting Abel's map \eqref{Amap} is known as the Jacobi inversion problem, see \cite[Section~III.6]{FarkasKra}. Since \( \RS \) has genus one, every Jacobi inversion problem has a unique solution.

\begin{proposition}
\label{prop:Ti}
Let \( \omega,\tau \) and \( c_\rho \) be as in \eqref{constants} and \eqref{veccrho}, respectively. Set
\begin{equation}
\label{xn-yn}
x_n := n\omega - \lfloor n\omega \rfloor, \quad y_n := n\tau - \lfloor n\tau \rfloor, \qandq \mathsf A:= (1+\omega)/2 + \mathsf B(1+\tau)/2,
\end{equation}
where \( \lfloor \cdot \rfloor \) is equal to the largest integer smaller or equal to its argument. Define
\begin{equation}
\label{Ti}
T_n(\z) := \exp\big\{-2\pi\ic y_na(\z)\big\}\dfrac{\theta \big( a(\z) - c_\rho - x_n - \mathsf By_n + \mathsf A \big)}{\theta \big( a(\z) + \mathsf A \big)}, \quad \z\in\RS_{\ualpha,\ubeta}.
\end{equation}
Each function \( T_n(\z) \) is meromorphic in the domain of its definition with continuous traces on \( (\ualpha\cup\ubeta)\setminus\{\boldsymbol a_1\}\) that satisfy
\begin{equation}
\label{T-jump}
T_{n+}(\s) = T_{n-}(\s)\left\{
\begin{array}{rl}
\displaystyle \exp\big\{-2\pi \ic \big( c_\rho + n\omega\big) \big\}, &  \s\in\ualpha\setminus\{\boldsymbol a_1\}, \medskip \\
\displaystyle \exp\big\{-2\pi\ic n\tau \big\}, &  \s\in\ubeta\setminus\{\boldsymbol a_1\}.
\end{array}
\right.
\end{equation}
In fact, \( T_n(\z) \) can be holomorphically continued across each of these cycles. It has precisely one pole, namely \( \infty^{(1)} \), and exactly one zero, namely \( \z_n \), both simple (the function will become analytic and non-vanishing if \( \z_n=\infty^{(1)} \)), where \( \z_n=\z_n(\rho) \in \RS \) is the unique solution of the Jacobi inversion problem
\begin{equation}
\label{jip}
\mathfrak a(\z_n) = \big[c_\rho + (n-1/2)(\omega + \mathsf B\tau)\big].
\end{equation} 
\end{proposition}

We prove Proposition~\ref{prop:Ti} in Section~\ref{ssec:theta}. 

\subsection{Nuttall-Szeg\H{o} Functions}

Now we are ready to described the solutions of \hyperref[bvp]{BVP-\(\Psi\)}.

\begin{theorem}
\label{thm:BANS}
Fix \( n\in\N \). Let \( c_\rho \) be given by \eqref{veccrho} and \( \z_n \) be the solution of \eqref{jip}. If \( \z_n=\infty^{(0)} \), then the boundary value problem \hyperref[bvp]{BVP-\(\Psi\)} does not have a solution. Otherwise,  any solution of \hyperref[bvp]{BVP-\(\Psi\)} is of the form \( c\Psi_n(\z) \) for a non-zero constant \( c \), where
\begin{equation}
\label{Psin}
\Psi_n(\z) := \big(\Phi^nS_\rho T_n \big)(\z), \quad \z\in\RS\setminus\bd,
\end{equation}
and the functions \( \Phi(\z) \), \( S_\rho(\z) \), and \( T_n(\z) \) are given by \eqref{Phi}, \eqref{szego}, and \eqref{Ti}, respectively.
\end{theorem}

We prove Theorem~\ref{thm:BANS} in Section~\ref{ssec:BANS}. Notice that in addition to the zero of order \( n-1 \) at \( \infty^{(1)} \), \( \Psi_n(\z) \) also has a simple zero at \( \z_n \) (if \(\z_n\) belongs to \( \bd\setminus\boldsymbol A \) the simplicity of the zero is understood as the simplicity of the zero of an analytic continuation of \( \Psi_n(\z)\) to a neighborhood of \( \z_n \); of course, if \( \z_n=\infty^{(1)} \), then the zero at \( \infty^{(1)} \) has order \( n \)) and the asymptotics of the behavior of \( \Psi_n(\z) \) around the points in \( \boldsymbol A \) can be improved to 
\begin{equation}
\label{Psin-ai}
\big|\Psi_n(z)\big| \sim |z-a_l|^{m_n(\boldsymbol a_l)/2-\alpha_l/2-1/4} \qasq z\to a_l,
\end{equation}
\( l\in\{0,1,2,3\} \), where \( m_n(\boldsymbol a)=1 \) if \( \z_n=\boldsymbol a \) and \( m_n(\boldsymbol a)=0 \) otherwise. 

Of course, the function \( \Psi_n(\z) \) is defined even when \( \z_n=\infty^{(0)} \). In this case it simply has a pole of order \( n-1 \) at \( \infty^{(0)} \) and therefore does not solve \hyperref[bvp]{BVP-\(\Psi\)}.  It readily follows from \eqref{jip} that if \( \z_n=\z_m \) for some distinct indices \( n \) and  \( m \), then \( [0] = [(n-m)(\omega+\mathsf B\tau)] \), which means that \( \omega \) and \( \tau \) are rational numbers with the denominators that divide \( n-m \). Thus, either \( \omega \) and \( \tau \) are rational numbers, in which case \( \{ \z_n \} \) is a periodic sequence, or  \( \z_n=\infty^{(0)} \) for at most one index \( n \). Unfortunately, for our asymptotic analysis we need to exclude not only indices for which \( \z_n=\infty^{(0)} \), but also those for which \( \z_n \) is close \( \infty^{(0)} \).  To this end, given \( \varepsilon>0 \) such that \( \max_{s\in\Delta}|s|<1/\varepsilon \), we let
\begin{equation}
\label{Ne}
\N_\varepsilon^\mathsf{reg} := \left\{n\in\N:~\z_n\not\in \boldsymbol U_\varepsilon^{(0)}\right\},
\end{equation}
where \( \boldsymbol U_\varepsilon^{(k)} \) is a neighborhood of \( \infty^{(k)} \) with natural projection equal to \( \{|z|>1/\varepsilon\} \).

\begin{proposition}
\label{prop:jip}
If \( \omega \) and \( \tau \) are rational numbers with the smallest common denominator \( d \), then there exists \( \varepsilon_0>0 \) such that either \( \N_\varepsilon^\mathsf{reg} = \N \) or \( \N_\varepsilon^\mathsf{reg} = \N\setminus (n_0+d\N) \) for some \(n_0\in\{0,\ldots,d-1\} \) and all \( \varepsilon<\varepsilon_0 \). In the latter case it necessarily holds that \( c_\rho=(p+\mathsf Bq)/d \) for some integers \( p \) and \( q \). Moreover, it holds that \( d\geq3 \) and \( d = 3 \) only if \( \Delta \) is the image under a linear map of
\begin{equation}
\label{Delta_sym}
\Delta_\mathsf{sym} = \big\{s:s^3\in[-1,0]\big\}.
\end{equation}

On the other hand, if at least one of the numbers \( \omega \) and \( \tau \) is irrational, then for any \( N\in\N \) there exists \( \varepsilon_N>0 \) such that
\begin{equation}
\label{gap-size}
 n \notin \N_\varepsilon^\mathsf{reg} \quad \Rightarrow \quad n+1,\ldots,n+N \in \N_\varepsilon^\mathsf{reg}
\end{equation}
for any \( 0<\varepsilon\leq\varepsilon_N \).  Moreover, there exists \( \delta=\delta(\varepsilon) \) such that
\begin{equation}
\label{n-1}
n \in \N_\varepsilon^\mathsf{reg} \quad \Rightarrow \quad \z_{n-1} \not\in \boldsymbol U_\delta^{(1)}.
\end{equation}
Furthermore, it is always true that \( \z_n=\infty^{(0)} \) if and only if \( \z_{n-1} = \infty^{(1)} \).
\end{proposition}

We prove Proposition~\ref{prop:jip} in Section~\ref{ssec:jip}. 

The reason we need to go from excluding indices for which \( \z_n \) is equal to \( \infty^{(0)} \) to sequences \( \N_\varepsilon^\mathsf{reg} \) is that we need to control asymptotic behavior of \( \Psi_n(\z) \).  More precisely, for our analysis it is important to know that functions \(  T_n(\z) \) and their reciprocals form normal families in certain subdomains of \( \RS \).

\begin{proposition}
\label{prop:Ti-as}
Let \( R>\max_{s\in\Delta}|s| \). There exists a constant \( M_R>0 \) such that
\begin{equation}
\label{T-bounds1}
\left\{
\begin{array}{ll}
|T_n(\z)| \leq  M_R, & \z\not\in \boldsymbol U_{1/R}^{(1)}, \medskip \\
|T_n(\z)| \leq M_R\frac{|z|}{R}, & \z\in \boldsymbol U_{1/R}^{(1)}, 
\end{array}
\right.
\end{equation}
for all \( n\in\N \). Furthermore, let \( \N_\varepsilon^\mathsf{reg} \) and \( \delta(\varepsilon) \) be as in Proposition~\ref{prop:jip}. Then there exists a constant \( M_\varepsilon>0 \) such that 
\begin{equation}
\label{T-bounds2}
\left\{
\begin{array}{ll}
|T_n^{-1}(\z)| \leq M_\varepsilon, & \z\in \boldsymbol U_{\varepsilon/2}^{(0)},  \medskip \\
|T_{n-1}^{-1}(\z)| \leq M_\varepsilon|z|^{-1}, & \z\in \boldsymbol U_{\delta(\varepsilon)/2}^{(1)},
\end{array}
\right.
\end{equation}
for \( n\in\N_\varepsilon^\mathsf{reg} \).
\end{proposition}

We prove Proposition~\ref{prop:Ti-as} in Section~\ref{ssec:Ti}.

\subsection{Local Obstruction}

To prove our main results on asymptotic behavior of orthogonal polynomials we utilize matrix Riemann-Hilbert analysis. As typical for such an approach, we need to construct a local parametrix around \( a_0 \) that models the behavior of the polynomials and their functions of the second kind there. In our case, this is done with the help of the matrix that solves Riemann-Hilbert problems that characterizes solutions of Painlev\'e XXXIV equation with the Painlev\'e variable equal to zero. It is known that such a matrix does not exist when the corresponding Painlev\'e function has a pole at zero. In this case we need to modify our construction and this modification requires placing an additional requirement on the sequence of allowable indices. More precisely, let us write
\begin{equation}
\label{Tn-exp}
T_n(z) := t_0^{(n)} + t_1^{(n)}(z-a_0)^{1/2} + t_2^{(n)}(z-a_0) + \cdots
\end{equation}
for \( z \) close to \( a_0 \) lying in the sector between \( \Delta_1 \) and \( \Delta_3 \) for some fixed determination of \( (z-a_0)^{1/2} \) that is analytic in this sector. Let \( h_1 \) be the constant defined further below in \eqref{GandH} and \eqref{H-exp} (it is the second coefficient of the Puiseux series at \( a_0 \) of a certain function that depends on the geometry via \( w(z) \) and the weight \( \rho(s) \) via \( S_\rho(z) \) and the function \( r(z) \) defined in \eqref{ri}). Define
\begin{equation}
\label{Nes}
\boldsymbol T_n := \left( \begin{matrix} 2t_2^{(n)} - h_1t_1^{(n)} & 2t_2^{(n-1)} - h_1t_1^{(n-1)} \medskip \\ t_0^{(n)} & t_0^{(n-1)} \end{matrix} \right).
\end{equation}
We need to require that the determinants of these matrices are separated away from zero.

\begin{proposition}
\label{prop:lo}
For each \( \varepsilon>0 \) small enough there exists \( \boldsymbol V_\varepsilon \), a neighborhood of \( \boldsymbol a_0 \) that can possibly be empty, such that if \( \z_n\notin \boldsymbol V_\varepsilon \), then either \( n \) or \( n+1 \) belongs to \( \N_\varepsilon^\mathsf{sing} \), where
\[
\N_\varepsilon^\mathsf{sing} := \Big\{n\in\N_\varepsilon^\mathsf{reg}:~| \det(\boldsymbol T_n) |\geq\varepsilon \Big\}.
\]
\end{proposition}

We prove Proposition~\ref{prop:lo} in Section~\ref{ssec:lo}. Let us note that there are no particular reasons to couple the constants in the definition of \( \N_\varepsilon^\mathsf{sing} \) and \( \N_\varepsilon^\mathsf{reg} \) except esthetic ones. That is, we could have defined \( \N_{\varepsilon_1,\varepsilon_2}^\mathsf{sing}\subseteq \N_{\varepsilon_1}^\mathsf{reg} \). 

Proposition~\ref{prop:lo} has the following implications. When at least one of the constants \( \omega \) and \( \tau \) is irrational, given a natural number \( N \),  \( \varepsilon \) can be taken small enough  so that \( \N_\varepsilon^\mathsf{reg} \) contains a consecutive block of \( N \) integers after each missing index by Proposition~\ref{prop:jip}. The same arguments as in the proof of that proposition can be employed to argue that \(\boldsymbol z_n \in \boldsymbol V_\varepsilon \) for at most one index in this block. Thus, according to Proposition~\ref{prop:lo}, essentially at least half of the indices in this block will belong to \( \N_\varepsilon^\mathsf{sing} \). If the constants \( \omega \) and \( \tau \) are rational with the smallest common denominator \( d\geq 3 \), then the sequence \( \{ \z_n \} \) is \( d \)-periodic and at most one index in this period can be such that the corresponding solution of the Jacobi inversion problem \eqref{jip} is equal to \( \boldsymbol a_0 \). Hence, it is always true that
\[
\N\setminus \big(n_0+d\N) \cup (n_1+d\N) \cup (n_1+1+d\N) \big) \subseteq \N_\varepsilon^\mathsf{sing},
\]
where \( n_0,n_1 \in \{0,\ldots,d-1\} \) are such that \( \z_{n_0}=\infty^{(0)} \) and \( \z_{n_1}=\boldsymbol a_0 \). Thus, \( \N_\varepsilon^\mathsf{sing} \) is again necessarily infinite, unless of course \( d=3 \). As we prove in Proposition~\ref{prop:jip}, in the latter case \( \Delta \) is a linear image of \( \Delta_\mathsf{sym} \). Further below in Appendix~\ref{ap:example} we construct a class of examples of weights of orthogonality on \( \Delta_\mathsf{sym} \) that do require the restriction to \( \N_\varepsilon^\mathsf{sing} \) and for which it is indeed the case that \( \z_0=\infty^{(0)} \), \( \z_1=\boldsymbol a_0 \), and \( \z_2=\infty^{(1)} \), while \( \N_\varepsilon^\mathsf{sing} = \varnothing \). Therefore neither of the forthcoming asymptotic results applies. Of course, one might surmise that this is simply due to the technical limitations of our methods. However, we do analyze the orthogonal polynomials for this class of weights and show that the main term of their asymptotics is actually \emph{different} from \( \Psi_n(z) \) constructed in Theorem~\ref{thm:BANS}. That is, the polynomials considered in Appendix~\ref{ap:example} do have different asymptotic behavior as compared to the polynomials in the rest of the cases, which does indicate that the obstruction that occurs during the construction of local parametrices is not a mere technicality.

It is also curious to observe that the condition defining \( \N_\varepsilon^\mathsf{reg} \) can be equivalently restated, see \eqref{curious0}, \eqref{curious1}, and \eqref{curious2} further below, as
\[
\varepsilon\leq \left| \det\left( \begin{matrix} t_1^{(n)} & t_1^{(n-1)} \medskip \\ t_0^{(n)} & t_0^{(n-1)} \end{matrix} \right) \right| \leq \varepsilon^{-1}.
\]

\section{Main Results}
\label{sec:main}

To state our main results on the asymptotics of the diagonal Pad\'e approximants to functions of the form \eqref{CI}, we need to separate the weights \( \rho(s) \) in \eqref{CI} into two classes. To this end, let \( \rho_i(s) \) be the restriction to \( \Delta_i^\circ \) of \( \rho(s) \) and \( (z-a_0)^{\alpha_0} \) be a fixed branch whose branch cut splits the sector between \( \Delta_2 \) and \( \Delta_3 \). Define
\begin{equation}
\label{varrho}
\varrho_i(s) := \rho_i(s)/(s-a_0)^{\alpha_0}, \quad s\in\Delta_i^\circ, \quad i\in\{1,2,3\},
\end{equation}
which, according to our assumptions stated right after \eqref{CI}, extend to analytic non-vanishing functions in some neighborhood of \( a_0 \). Define constants
\begin{equation}
\label{stokes}
\begin{cases}
\displaystyle b_1(\rho) := - \frac{e^{\pi\ic\alpha_0}\varrho_2(a_0)+e^{-\pi\ic\alpha_0}\varrho_3(a_0)}{\varrho_1(a_0)}, \medskip \\ 
\displaystyle b_i(\rho) := - \frac{\varrho_{i-1}(a_0)+\varrho_{i+1}(a_0)}{e^{(-1)^i\pi\ic\alpha_0}\varrho_i(a_0)}, \quad i\in\{2,3\}, 
\end{cases}
\end{equation}
where subindices are understood cyclically in \( \{1,2,3\} \). The parameters \( \alpha_0 \) and \( b_i(\rho) \), \( i\in\{1,2,3\} \), uniquely determine a \( 2\times2 \) sectionally analytic matrix function \( \boldsymbol \Phi_{\alpha_0}(\zeta;x) \) as a solution of the Riemann-Hilbert problem \hyperref[rhphi]{\rhphi}, see Appendix~\ref{ap:34}, where \( x\in\C \) is a parameter that appears in the description of the asymptotic behavior at infinity. As we explain in Appendix~\ref{ap:34}, \( \boldsymbol \Phi_{\alpha_0}(\zeta;x) \) is meromorphic in \( x \). 

\begin{definition}
In what follows we shall say that a weight \( \rho(s) \) on \( \Delta \) belongs to the class \( W^\mathsf{reg} \) if the matrix function \( \boldsymbol \Phi_{\alpha_0}(\zeta;0) \) solving Riemann-Hilbert problem \hyperref[rhphi]{\rhphi} with \( \alpha=\alpha_0 \), \( b_i=b_i(\rho) \), \( i\in\{1,2,3\} \), and \( x=0 \) exists. That is, if \( \boldsymbol \Phi_{\alpha_0}(\zeta;x) \) does not have a pole at \( x=0 \). Otherwise, we shall say that the weight \( \rho(s) \) belongs to the class \( W^\mathsf{sing} \).
\end{definition}

It can be readily verified that parameters \eqref{stokes} satisfy the relation
\begin{equation}
\label{bi}
b_1+b_2+b_3-b_1b_2b_3 = 2\cos(\pi\alpha_0).
\end{equation}
We shall call any solution of \eqref{bi} a set of Stokes parameters. Any set of Stokes parameters satisfying \eqref{bi} leads to its own matrix \( \boldsymbol \Phi_{\alpha_0}(\zeta;x) \). Some of these matrices will have a pole at \( x=0 \) and some will not. Essentially any such set can be obtained via \eqref{stokes} (in particular, this means that the class \( W^\mathsf{sing} \) is non-empty). Indeed, given a set of Stokes parameters \eqref{bi} for which none is equal to \( e^{\alpha_0\pi\ic} \) (in this case it is also true that none is equal to \( e^{-\alpha_0\pi\ic} \) and a product of no two parameters is equal to \( 1 \), hence, the third one is uniquely determined by the other two), they are realized by a weight \( \rho(s) \) such that
\[
\frac{\varrho_1(a_0)}{\varrho_2(a_0)} = - \frac{1-b_2b_3}{1-e^{-\alpha_0\pi\ic}b_3} \qandq \frac{\varrho_3(a_0)}{\varrho_2(a_0)} =  \frac{1-e^{\alpha_0\pi\ic}b_2}{1-e^{-\alpha_0\pi\ic}b_3}. 
\]
On the other hand, if one of the Stokes parameters in \eqref{bi} is equal to \( e^{\alpha_0\pi\ic} \) or \( e^{-\alpha_0\pi\ic} \), then one (whichever) of the other two must be equal to \( e^{-\alpha_0\pi\ic} \) or \( e^{\alpha_0\pi\ic} \), respectively, and the third one could be absolutely arbitrary. However, this is not the case for the parameters \( b_i(\rho) \). If one of them happens to be \( e^{\alpha_0\pi\ic} \) or \( e^{-\alpha_0\pi\ic} \), then the other two are uniquely determined.

Unfortunately, currently it is unknown how to determine solely from the Stokes parameters \( b_i(\rho) \) whether \( \rho\in W^\mathsf{reg} \) or \( \rho\in W^\mathsf{sing} \). Such a determination amounts to distinguishing Stokes parameters that lead to finite initial conditions at \( x=0 \) for the corresponding solution of Painlev\'e XXXIV equation, see \eqref{Uu} and \eqref{ant-u}, further below, and those leading to solutions with a pole at \( x=0 \).

\begin{theorem}
\label{thm:asymp1}
Let \( f(z) \) be a function given by \eqref{CI} with \( \alpha_0\in(-1,1) \) and \( \rho\in W^\mathsf{reg} \). Let \( Q_n(z) \) be the minimal degree monic denominator of the diagonal Pad\'e approximants \( [n/n]_f(z) \) and \( R_n(z) \) be the corresponding linearized error function \eqref{linear}. Further, let \( \Psi_n(\z) \) be given by \eqref{Psin} and \( \N_\varepsilon^\mathsf{reg} \) be as in Proposition~\ref{prop:jip}. Then it holds for all \( n\in\N_\varepsilon^\mathsf{reg} \) large enough that
\begin{equation}
\label{Asymp1}
\begin{cases}
Q_n(z) & = \gamma_n \big(1 + \upsilon_{n1}(z) \big)\Psi_n(z) +\gamma_n \upsilon_{n2}(z) \Psi_{n-1}(z), \smallskip \\
(wR_n)(z) & =  \gamma_n \big(1 + \upsilon_{n1}(z) \big)\Psi_n^*(z) +\gamma_n \upsilon_{n2}(z) \Psi_{n-1}^*(z),
\end{cases}
\end{equation}
for $z\in\overline \C\setminus\Delta$ (in particular, \( \deg(Q_n)=n \) for such indices \( n \)), where \( \gamma_n:=\lim_{z\to\infty}z^n\Psi_n^{-1}(z) \) is the normalizing constant, \( w(z) \) is given by \eqref{w}, and the functions \( \upsilon_{ni}(s) \) vanish at infinity and satisfy\footnotemark 
\begin{equation}
\label{upsilons}
\upsilon_{ni}(z) = \mathcal O_\varepsilon\big( n^{-1/3} \big)
\end{equation}
locally uniformly on \( \overline\C\setminus\Delta \). Moreover, it holds for all \( n\in\N_\varepsilon^\mathsf{reg} \) large enough that
\begin{equation}
\label{Asymp2}
\begin{cases}
Q_n(s) & = \gamma_n\big(1 + \upsilon_{n1}(s) \big)\left(\Psi_{n+}(s) + \Psi_{n-}(s) \right) + \gamma_n \upsilon_{n2}(s)\left(\Psi_{n-1+}(s) + \Psi_{n-1-}(s) \right), \smallskip \\
(wR_n)_\pm(s) & = \gamma_n \big(1 + \upsilon_{n1}(s) \big)\Psi_{n\pm}^*(s) +\gamma_n \upsilon_{n2}(s) \Psi_{n-1\pm}^*(s),
\end{cases}
\end{equation}
for \( s\in\Delta^\circ \), where the functions \( \upsilon_{ni}(s) \) can be analytically continued into a neighborhood of \( \Delta^\circ \) and satisfy \eqref{upsilons} locally uniformly on \( \Delta^\circ \).
\end{theorem}
\footnotetext{The notation $f(n) = \mathcal{O}_\varepsilon (g(n))$ means that there exists a constant $C_\varepsilon$ depending on $\varepsilon$ such that $|f(n)| \leq C_\varepsilon |g(n)|$.}

Theorem~\ref{thm:asymp1} is proved in Section~\ref{s:RH1}. To prove it we use by now classical approach of Fokas, Its, and Kitaev \cite{FIK91,FIK92} connecting orthogonal polynomials to matrix Riemann-Hilbert problems. The RH problem is then analyzed via the nonlinear steepest descent method of Deift and Zhou \cite{DZ93}.

The error estimates in \eqref{upsilons} can be improved in some cases.
\begin{itemize}
\item  It is known, see \cite{ApYa15}, that if \( \alpha_0=0 \) and \( \rho_1(z)+\rho_2(z)+\rho_3(z) \equiv 0 \) in some neighborhood of \( a_0 \), then the error estimates improve to \( \mathcal O_\varepsilon (n^{-1}) \) (in this case \( b_i(\rho)=1 \) and \( \boldsymbol \Phi_0(\zeta;0) \) admits an explicit construction that uses Airy functions, so such weights do belong to \( W^\mathsf{reg}\)). 
\item If \( \rho(s) = \rho^*(s)/w_+(s) \) for some function \( \rho^*(s) \) that has an analytic extension to some neighborhood of \( \Delta \), then  the error estimates are geometric, see \cite{Ya15} (in this case \( b_i(\rho)=0 \) and \( \boldsymbol \Phi_{-1/2}(\zeta;x) \) is equal to the right-hand side of \hyperref[rhphi]{\rhphi}(d) with \( U_n=V_n\equiv 0 \)).
\end{itemize}

In the case of singular weights, the following theorem takes place.

\begin{theorem}
\label{thm:asymp2}
Let \( f(z) \) be a function given by \eqref{CI} with \( \alpha_0\in(-1,1) \) and \( \rho\in W^\mathsf{sing} \). {Assume further that
\begin{equation}
\label{bad-cond}
\varrho_1^\prime(a_0)/\varrho_1(a_0) = \varrho_2^\prime(a_0)/\varrho_2(a_0) = \varrho_3^\prime(a_0)/\varrho_3(a_0).
\end{equation}}
Let \( Q_n(z) \), \( R_n(z) \), and  \( \Psi_n(\z) \) be as in Theorem~\ref{thm:asymp1}. Let the matrices \( \boldsymbol T_n \) be as in \eqref{Nes} and \( \N_\varepsilon^\mathsf{sing} \) be as in Proposition~\ref{prop:lo}. Then it holds for all \( n\in\N_\varepsilon^\mathsf{sing} \) large enough that asymptotic formulae \eqref{Asymp1} and \eqref{Asymp2} remain valid but with \( \upsilon_{n1}(z) \) and \( \upsilon_{n2}(z) \) replaced by
\[
\upsilon_{n1}(z) -  \frac2{z-a_0}\frac{t_0^{(n-1)}t_0^{(n)}}{\det(\boldsymbol T_n)} \qandq \upsilon_{n2}(z) + \frac2{z-a_0}\frac{(t_0^{(n)})^2}{\det(\boldsymbol T_n)},
\]
respectively, where the error functions \( \upsilon_{ni}(z) \) possess the same properties as in Theorem~\ref{thm:asymp1}.
\end{theorem}

Let us point out that the new term added to the error functions \( \upsilon_{ni}(z) \) is not asymptotically small.  {The presence of the condition \eqref{bad-cond} is connected to the local analysis around \( a_0 \) and seems to be inescapable for our techniques. In particular, we needed this condition for our analysis of non-Hermitian orthogonal polynomials on a cross in \cite{BarYa20}, see the definitions of classes \( \mathcal W_\ell \) on page~2. Its nature is not entirely transparent to us at this stage. However, as we explain in Appendix~\ref{ap:functions}, functions of the form \eqref{log-power} do satisfy~\eqref{bad-cond}.} Theorem~\ref{thm:asymp2} is proved in Section~\ref{s:RH2}.

\section{Proof of Propositions~\ref{prop:szego}--\ref{prop:lo} and Theorem~\ref{thm:BANS}}
\label{s:NS}

\subsection{Proof of Proposition~\ref{prop:szego}}
\label{ssec:szego}

Orient \( \pi(\ualpha) \) towards \( a_1 \). Let
\[
H(z) := \frac{w(z)}{2\pi\ic}\int_{\pi(\ualpha)}\frac1{t-z}\frac{\dd t}{w(t)}, \quad z\in\C\setminus\big(\Delta\cup\pi(\ualpha)\big).
\]
Then \( H(z) \) is an analytic function in its domain of definition (the integral itself is analytic in \( \overline\C\setminus\pi(\ualpha) \)) that behaves like
\begin{equation}
\label{S1}
H(z) = -\frac z{2\pi\ic}\int_{\pi(\ualpha)}\frac{\dd s}{w(s)} + \mathcal O(1) = -\frac z{4\pi\ic}\oint_{\ualpha}\frac{\dd s}{w(\s)} + \mathcal O(1)
\end{equation}
as \( z\to\infty \). We get from Plemelj-Sokhotski formulae \cite[Section~I.4.2]{Gakhov} that
\begin{equation}
\label{S2}
H_+(s) - H_-(s) = 1, \quad s\in\pi(\ualpha)\setminus\{a_1,a_3\}.
\end{equation}
It further follows from \cite[Section~8.3]{Gakhov} that \( H(z) \) is bounded around \( a_1 \) and \( a_3 \). Set
\begin{equation}
\label{Lambda}
\Lambda(z) := \frac{w(z)}{2\pi\ic}\int_\Delta\frac{\log(\rho w_+)(s)}{s-z}\frac{\dd s}{w_+(s)}, \quad z\in\C\setminus\Delta. 
\end{equation}
This is a holomorphic function in its domain of definition that satisfies
\begin{equation}
\label{S3}
\Lambda(z) = -\frac z{2\pi\ic}\int_\Delta\frac{\log(\rho w_+)(s)}{w_+(s)}\dd s + \mathcal O(1) = -\frac{zc_\rho}{2}\oint_{\ualpha}\frac{\dd s}{w(\s)} + \mathcal O(1)
\end{equation}
as \( z\to\infty \), see \eqref{veccrho} and \eqref{B}. It follows from Plemelj-Sokhotski formulae that
\begin{equation}
\label{S4}
\Lambda_+(s) + \Lambda_-(s) = \log(\rho_i w_+)(s), \quad s\in\Delta_i^\circ, ~~i\in\{1,2,3\}.
\end{equation}
Write \( \Lambda(z) = \Lambda_1(z) + \Lambda_2(z) + \Lambda_3(z) \), where \( \Lambda_i(z) \) is given by \eqref{Lambda} with \( \Delta \) replaced by \( \Delta_i \). Fix \( i\in\{1,2,3\} \). Let \( \sqrt{z-a_i} \) and \( (z-a_i)^{\alpha_i} \) be branches whose branch cuts include \( \Delta_i \) and such that \( w(z) = v_i(z)\sqrt{z-a_i} \) and \( \rho_i(s) = \tilde\rho_i(s)(s-a_i)^{\alpha_i}_+ \) for some functions \( v_i(z) \) and \( \tilde\rho_i(z) \) analytic in a vicinity of \( a_i \). Then we get from \cite[Sections~I.8.3-6]{Gakhov} that
\begin{eqnarray}
\Lambda_i(z) &=& \frac{v_i(z)\sqrt{z-a_i}}{2\pi\ic}\int_{\Delta_i}\frac{(\alpha_i+1/2)\log_+(s-a_i) + \log(\tilde\rho_iv_i)(s)}{s-z}\frac{\dd s}{v_i(s)\sqrt{s-a_i}_+} \nonumber \\
& = & \frac{v_i(z)\sqrt{z-a_i}}{2\pi\ic}\int_{\Delta_i}\frac{(\alpha_i+1/2)\log_+(s-a_i)}{s-z}\frac{\dd s}{v_i(s)\sqrt{s-a_i}_+} + \mathcal O(1)\nonumber \\
\label{S5}
&=& (\alpha_i/2+1/4)\log(z-a_i) + \mathcal O(1)
\end{eqnarray}
as \( z\to a_i \), where \( \log(z-a_i) \) has a branch cut that includes \( \Delta_i \). 

Similarly, it follows from the conditions placed on \( \rho_i(s) \) that \( \rho_i(s)/(s-a_0)_+^{\alpha_0} \) extends to an analytic function in a vicinity of \( a_0 \), where \( (z-a_0)^{\alpha_0} \) is a branch with the cut along \( \Delta_i \). Moreover, to get a function analytic in a vicinity of \( a_0 \) out of \( w(z) \), we need to divide it by \( \sqrt{z-a_0} \) with the jump along \( \Delta_i \) and then multiply by \( -1 \) in the sector opposite to \( \Delta_i \). Hence, it holds that
\begin{equation}
\label{S6}
\Lambda_i(z) = \nu_i(\alpha_0/2+1/4)L_i(z) +\mathcal O(1)
\end{equation}
as \( z\to a_0 \), where \( L_i(z) := \log(z-a_0) \) has a branch cut that includes \( \Delta_i \) and \( \nu_i \) is equal to \( 1 \) in the sectors adjacent to \( \Delta_i \) and is equal to \( -1 \) in the sector opposite to \( \Delta_i \). Since
\[
S_\rho\big( z^{(k)}\big) = \exp\left\{(-1)^{k+1}\left(\Lambda(z)-2\pi\ic c_\rho H(z)\right)\right\},
\]
see \eqref{cauchy-kernel} and \eqref{szego}, the claims of the proposition follow (\eqref{S1} and \eqref{S3} give analyticity at infinity; \eqref{S2} and \eqref{S4} yield \eqref{S-jump}; \eqref{S5} and \eqref{S6} give \eqref{Srho-ai}).

\subsection{Proof of Proposition~\ref{prop:Ti}}
\label{ssec:theta}

Before we start proving this proposition in earnest, let us make an observation. Recall the differentials \( \Omega_{\z,\z^*}(\s) \) introduced in \eqref{szego} and the differential \( G(\s) \) defined in \eqref{Phi}. One can easily check that
\[
G(\s) = \Omega_{\infty^{(1)},\infty^{(0)}}(\s) + 2\pi\ic \tau\Omega(\s),
\]
where \( \Omega(\s) \) is given by \eqref{B}. It further follows from Riemann's relations \cite[Section~III.3.6]{FarkasKra} that
\[
\oint_{\ubeta}\Omega_{\infty^{(1)},\infty^{(0)}} = 2\pi\ic\int_{\infty^{(0)}}^{\infty^{(1)}}\Omega,
\]
where the path of integration lies entirely in \( \RS_{\ualpha,\ubeta} \). Therefore, we can immediately deduce from the last two identities, \eqref{B}, and \eqref{constants} that
\begin{equation}
\label{0infty}
\int_{\infty^{(1)}}^{\infty^{(0)}}\Omega = \omega + \mathsf B\tau \qandq \int_{\boldsymbol a_0}^{\infty^{(k)}} \Omega = (-1)^k(\omega+\mathsf B\tau)/2, \quad k\in\{0,1\},
\end{equation}
where we used symmetry \( \Omega(\z^*)=-\Omega(\z) \) for the last conclusion.

It is well known and can be checked directly that $\theta(u)$ enjoys the following periodicity properties:
\begin{equation}
\label{RTF-per}
\theta(u + l + \mathsf Bm) = \exp\big\{-\pi \ic\mathsf B m^2 - 2\pi \ic mu\big\}\theta(u), \quad l,m\in\Z.
\end{equation}
It is further known that \( \theta(u) = 0 \) if and only if \( [u] = [(\mathsf B+1)/2] \). As already mentioned, the function \( a(\z) \) from \eqref{a} extends to a multi-valued analytic function on \( \RS \). Therefore, it has continuous traces on \( (\ualpha\cup\ubeta)\setminus\{\boldsymbol a_1\} \) and these traces satisfy
\begin{equation}
\label{a-jump}
a_+(\s) - a_-(\s) = \left\{
\begin{array}{rl}
-\mathsf B, & \s\in\ualpha\setminus\{\boldsymbol a_1\}, \medskip \\
	1, & \s\in\ubeta\setminus\{\boldsymbol a_1\},
\end{array}
\right.
\end{equation}
due to normalization of \( \Omega(\z) \) and the definition of \( \mathsf B \), see \eqref{B}. Relations \eqref{RTF-per} and \eqref{a-jump} imply that \( T_n(\z) \) does indeed extend to a multiplicatively multi-valued function on \( \RS \). Furthermore, they also yield that
\begin{eqnarray*}
T_{n+}(\s) &=& \exp\big\{ -2\pi\ic y_n (a_-(\s) - \mathsf B)\big\}\frac{\theta(a_-(\s)-c_\rho-x_n-\mathsf By_n+\mathsf A-\mathsf B)}{\theta(a_-(\s)+\mathsf A-\mathsf B)} \\
& = & T_{n-}(\s)\exp\{2\pi\ic(-c_\rho-x_n)\}, \quad \s\in\ualpha\setminus\{\boldsymbol a_1\},
\end{eqnarray*}
from which the first relation in \eqref{T-jump} follows (recall that \( x_n \) and \( n\omega \) differ by an integer). The proof of the second one is analogous. Moreover, since
\[
[(1+\mathsf B)/2] = [a(\z)-c_\rho-x_n-\mathsf By_n+\mathsf A]  = [a(\z)-a(\z_n) + (1+\mathsf B)/2] 
\]
holds only when \( \z=\z_n \) by the unique solvability of the Jacobi inversion problem (we used \eqref{jip} and \eqref{xn-yn} for the second equality), \( T_n(\z) \) does indeed vanish solely at \( \z_n \). The location of the pole of \( T_n(\z) \) can be identified similarly using \eqref{0infty}.

\subsection{Proof of Theorem~\ref{thm:BANS}}
\label{ssec:BANS}

Fix \( n\in\N \) and let \( \Psi_n(\z) \) be given by \eqref{Psin}. It follows from the properties of \( \Phi(\z) \), \(S_\rho(\z) \), and \(T_n(\z) \), including jump relations \eqref{Phi-jump}, \eqref{S-jump}, and \eqref{T-jump}, and local behavior \eqref{Srho-ai} that \( \Psi_n(\z) \) solves \hyperref[bvp]{BVP-\(\Psi\)} when \( \z_n\neq\infty^{(0)} \) (when \( \z_n=\infty^{(0)} \), the pole at \( \infty^{(0)} \) is of order \( n-1 \)). Let \( \Psi(\z) \) be any solution of \hyperref[bvp]{BVP-\(\Psi\)}, if it exists. Consider the function \( F(\z) := \Psi(\z)/\Psi_n(\z) \). Then \( F(\z) \) is meromorphic in \( \RS\setminus\bd \) with at most one pole, namely \( \z_n \), which is simple if present. It also holds that \( F_+(\s) = F_-(\s) \) for \( \s\in\bd\setminus\boldsymbol A \). Since these traces are continuous on \( \bd\setminus\boldsymbol A \) (apart from a possible simple pole at \( \z_n \) when the latter lies on \( \bd\setminus\boldsymbol A \)), the analytic continuation principle yields that \( F(\z) \) is analytic in \( \RS\setminus(\{\z_n\}\cup\boldsymbol A) \) with at most a simple pole at \( \z_n \) when \( \z_n\not\in\boldsymbol A \). Finally, it follows from \eqref{BVP-Psi} and \eqref{Psin-ai} that every \( \boldsymbol a_i \) is a removable singularity of \( F(\z) \) unless it coincides with \( \z_n \), in which case it can be a simple pole. Altogether, \( F(\z) \) is a rational function on \( \RS \) with at most one pole, which is simple. Hence, \( F(\z) \) is a constant. That is, \( \Psi(\z)=c\Psi_n(\z) \), which finishes the proof of the theorem.

\subsection{Proof of Proposition~\ref{prop:jip}}
\label{ssec:jip}

Below we always assume that all \( \varepsilon_N \) satisfy \( \displaystyle\max_{s\in\Delta}|s|<\varepsilon_N^{-1} \). 

Because Jacobi inversion problems are uniquely solvable in the considered case, it follows from \eqref{Amap} and \eqref{0infty} that if \( \z_n = \infty^{(k)} \), then 
\[
\big[c_\rho + (n-1/2)(\omega + \mathsf B\tau)\big] = [(-1)^k(\omega + \mathsf B\tau)/2]
\]
and therefore \( \z_{n-(-1)^k}=\infty^{(1-k)} \), again by \eqref{Amap} and \eqref{0infty}, since
\[
\big[c_\rho + (n-1/2-(-1)^k)(\omega + \mathsf B\tau)\big] = [-(-1)^k(\omega + \mathsf B\tau)/2].
\]
This proves the very last claim of the proposition.

If \( \omega \) and \( \tau \) are rational numbers with the smallest common denominator \( d \), then
\[
 \big[c_\rho + (n+d-1/2)(\omega + \mathsf B\tau)\big] =  \big[c_\rho + (n-1/2)(\omega + \mathsf B\tau)\big]
\]
and therefore the sequence \( \{\z_n\} \) is \( d \)-periodic. Thus, if none of \( \z_0,\ldots,\z_{d-1} \) coincides with  \( \infty^{(0)} \), choose \( \varepsilon_0 \) small enough so that \( \z_0,\ldots,\z_{d-1}\not\in \boldsymbol U_{\varepsilon_0}^{(0)} \), in which case \( \N_\varepsilon^\mathsf{reg}=\N \) for all \( \varepsilon<\varepsilon_0 \). If \( \z_{n_0} = \infty^{(0)} \) for some \( n_0\in\{0,\ldots,d-1 \} \), then of course \( \z_{n_0+kd} = \infty^{(0)} \) for all \( k\geq 0 \) and \( \z_{l + kd} \neq \infty^{(0)} \) for \( l\in\{0,\ldots,d-1\}\setminus\{n_0\} \), since otherwise \( l-n_0 \) would be divisible by \( d \), which is impossible. Hence, by periodicity, \( \N_\varepsilon^\mathsf{reg} = \N\setminus (n_0+d\N) \) for all \( \varepsilon<\varepsilon_0 \), where \( \varepsilon_0 \) is such that \( \z_l\not\in \boldsymbol U_{\varepsilon_0}^{(0)} \) for \( l\in\{0,\ldots,d-1\}\setminus\{n_0\} \). Moreover, in this case we get that there exist integers \( l \) and \( m \) such that
\[
c_\rho = \big( (1-n_0)\omega + l\big) + \mathsf{B}\big( (1-n_0)\tau + m\big) = (p/d)+\mathsf B(q/d)
\]
by \eqref{Amap} and \eqref{0infty}, where \( p \) and \( q \) are again integers.

Now, recall the definition of the constants \( \omega \) and \( \tau \) in \eqref{constants}. We can equivalently write
\begin{equation}
\label{prop3a}
\omega = I_1+I_3 \qandq \tau = -(I_1+I_2), \quad I_i := -\frac1{\pi\ic} \int_{\Delta_i}\frac{(s-a_0)\dd s}{w_+(s)}
\end{equation}
(recall that each \( \Delta_i \) is oriented towards \( a_0 \)). It follows from the Cauchy integral formula that
\begin{equation}
\label{prop3b}
1 = \frac1{2\pi\ic}\oint_L\frac{(s-a_0)\dd s}{w(s)} = -\frac1{\pi\ic}\int_\Delta\frac{(s-a_0)\dd s}{w_+(s)} = I_1+I_2+I_3,
\end{equation}
where \( L \) is any smooth Jordan curve that encircles \( \Delta \) in the counter-clockwise direction. We also get from \eqref{differential} that the integrand in \eqref{prop3b}, i.e.,
\begin{equation}
\label{prop3c}
\frac{(s-a_0)}{w_+(s)}\frac{\dd s}{\pi\ic}
\end{equation}
is real and non-vanishing on each \( \Delta_i^\circ \). In particular, \( I_i\neq 0 \), \( i\in\{1,2,3\} \). Moreover, it is known, see \cite[Theorem~8.1]{Pommerenke2}, that there exists a conformal map \( \phi(z) \) in a neighborhood of \( a_0 \) such that the differential \eqref{differential} can be locally written as \( \phi\dd\phi^2 \). If we normalize  \( \phi(z) \) so that it maps \( \Delta_1 \) into non-positive reals (locally around \(a_0\)), then \( \sqrt\phi\dd\phi \) is equal to \eqref{prop3c} in the sector delimited by \( \Delta_2 \) and \( \Delta_3 \) that contains \( \Delta_1 \) and minus the differential in \eqref{prop3c} in the other sector. This yields that \( I_i>0 \), \( i\in\{1,2,3\} \). In particular, the latter claim implies that \( \omega\in(0,1)\), \( \tau \in (-1,0) \), and for no \( \Delta \) it can hold that \( \omega=-\tau=1/2 \) as \( \omega-\tau = 1 + I_1>1 \), that is, \( d\geq 3 \). 

Assume that \( d=3 \). Then \( I_i=1/3 \), \( i\in\{1,2,3\} \), as these are positive rational numbers smaller than \( 1 \) with denominator \( 3 \) by \eqref{prop3a} and \eqref{prop3b}. It readily follows from its symmetries that this is indeed the case for \( \Delta_\mathsf{sym} \). Let as before \( \Phi(z) \) stand for the pull-back of \( \Phi(\z) \) from \( \RS^{(0)} \) onto \( \C\setminus\Delta \). Denote by \( \Phi_\Delta(z) \) the analytic continuation of  \( \Phi(z) \) from infinity across \( \pi(\ualpha) \) and \( \pi(\ubeta) \) to a holomorphic function in \( \C\setminus\Delta \). In fact, \( \Phi_\Delta(z) \) is again given by \eqref{Phi}, where \( \z\in\RS^{(0)}\setminus\bd  \) and the path of integration also belongs to \( \RS^{(0)}\setminus\bd  \), except for \( \boldsymbol a_0 \), and proceeds from \( \boldsymbol a_0 \) into the sector delimited by \( \bd_2 \) and \( \bd_3 \), recall Figure~\ref{fig:1}. As was explained after \eqref{Phi-jump}, \( \Phi_\Delta(z) \) is a conformal map of \( \C\setminus\Delta \) onto \( \{|z|>1\} \) that maps infinity into infinity. This map also sends \( \partial(\C\setminus\Delta) \), the boundary of \( \C\setminus\Delta \) understood as the collection of limit points of sequences in \( \C\setminus\Delta \), bijectively onto the unit circle. Clearly, \( \partial(\C\setminus\Delta) \) consists of two copies of each \( s\in\Delta_i^\circ \), say \( s_+ \) and \( s_- \), as accessed from the \( + \) and \( - \) sides of \( \Delta_i \), one copy of each \( a_1,a_2,a_3 \), and three copies of \( a_0 \), say \( a_{0,12} \), \( a_{0,23} \), and \( a_{0,31} \), where \( a_{0,ij} \) is accessed from the sector delimited by \( \Delta_i \) and \( \Delta_j \). Then
\[
\begin{cases}
\Phi_\Delta(a_{0,23}) =1, \smallskip \\
\Phi_\Delta(a_{0,12}) = e^{\pi\ic (-2I_2)} = e^{-2\pi\ic/3}, \smallskip \\
\Phi_\Delta(a_{0,31}) = e^{\pi\ic (2I_3)} = e^{2\pi\ic/3},    
\end{cases}
\qandq
\begin{cases}
\Phi_\Delta(a_1) = e^{\pi\ic(2I_3+I_1)} = -1, \smallskip \\
\Phi_\Delta(a_2) = e^{\pi\ic (-I_2)} = e^{-\pi\ic/3}, \smallskip \\
\Phi_\Delta(a_3) = e^{\pi\ic I_3} = e^{\pi\ic/3}.
\end{cases}
\]
Moreover, as \( w_+(s) = -w_-(s) \), it also holds that \( \Phi_\Delta(s_+)/\Phi_\Delta(a_i) = \overline{\Phi_\Delta(s_-)/\Phi_\Delta(a_i)} \) for any \( s\in\Delta_i^\circ \). So, if  \( \Delta^1 \) and \( \Delta^2 \) are two such contours, then \( F(z) := \big(\Phi_{\Delta^2}^{-1}\circ \Phi_{\Delta^1}\big)(z) \) is a conformal map of \( \C\setminus\Delta^1 \) onto \( \C\setminus\Delta^2 \) that maps infinity into infinity and \( \partial(\C\setminus\Delta^1) \) continuously and bijectively onto \( \partial(\C\setminus\Delta^2) \). The just described properties also yield that
\[
F(s_+) = F(s_-), \quad s\in (\Delta_i^1)^{\circ}, \qandq F(a_l(\Delta^1)) = a_l(\Delta^2)
\]
for \( i\in\{1,2,3\} \) and \( l\in\{0,1,2,3\} \). Thus, \( F(z) \) is continuous and therefore holomorphic in the entire complex plane. Since \( F(z) \) is conformal around infinity, it has a simple pole there, which means that it is a linear function that maps \( \Delta^1 \) onto \( \Delta^2 \), as claimed.

Assume now that at least one the numbers \( \omega \) and \( \tau \) is irrational. As explained before the statement of the proposition, \( \z_n \) can be equal to \( \infty^{(0)} \) for at most one index \( n \). Fix \( N\geq 1 \). Let \( \x_{-1},\x_0,\x_1,\ldots,\x_N \) be the unique points on \( \RS \) such that
\[
\mathfrak a(\x_l) = \left[\frac{\omega+\mathsf B\tau}2 + \big(\omega+\mathsf B\tau \big)l \right], \quad l\in\{-1,0,\ldots,N\}.
\]
It follows from \eqref{0infty} that \( \x_{-1} = \infty^{(1)} \) and \( \x_0 = \infty^{(0)} \). It also holds that none of \( \x_1,\ldots,\x_N \) is equal to \( \infty^{(0)} \) and they are all distinct (if \( \x_i = \x_j \) for \( i\neq j \), then \( \omega \) and \( \tau \) are necessarily rational). Fix \( \epsilon_*>0 \) small enough so that the sets
\[
\mathsf O_l := \Big\{ [s]~\big|~\exists x:~[x]=\mathfrak a(\x_l),~~|s-x|<\epsilon_* \Big\} \subset \mathsf{Jac}(\RS)
\]
are disjoint. Recall that the Abel's map \( \mathfrak a(\z) \) is a holomorphic bijection between \( \RS \) and \( \mathsf{Jac}(\RS) \). Set
\[
\boldsymbol U_l := \mathfrak a^{-1} ( \mathsf O_l ), \quad l\in\{-1,0,\ldots,N\}.
\] 
Notice that if \( [s]\in \mathsf O_l \), then \( [s+\omega+\mathsf B\tau]\in \mathsf O_{l+1} \). Hence, the domains \( \boldsymbol U_l \) are disjoint and
\[
\z_{n+k} \in \boldsymbol U_k, \quad k\in\{-1,\ldots,N-1\}  \quad \Rightarrow \quad \z_{n+l} \in \boldsymbol U_l, \quad l\in\{k+1,\ldots,N\}.
\]
Thus, to prove \eqref{gap-size} it remains to choose \( \varepsilon_N \) small enough so that \( \boldsymbol U_{\varepsilon_N}^{(0)}\subseteq  \boldsymbol U_0 \), and to prove \eqref{n-1}, given \( 0<\varepsilon<\varepsilon_N \), to choose \( \epsilon_* \) and \( \delta(\varepsilon) \) so that \( \boldsymbol U_0 \subseteq \boldsymbol U_\varepsilon^{(0)} \) and \( \boldsymbol U_{\delta(\varepsilon)}^{(1)} \subseteq \boldsymbol U_{-1} \).

\subsection{Proof of Proposition~\ref{prop:Ti-as}}
\label{ssec:Ti}

Since \( a(\z) \) is a branch of an additively multi-valued function on \( \RS \), see \eqref{a-jump}, it holds that \( K_R \), the closure of \( a(\RS_{\ualpha,\ubeta}\setminus \boldsymbol U_{1/R}^{(1)}) \) is a compact set. Moreover, \( K_R \) is necessarily disjoint from any point on the lattice \( [-(\omega+\mathsf B\tau)/2] \) by \eqref{0infty}. Hence, there exists a constant \( M_{1,R} \) such that
\[
\Big|e^{-2\pi\ic y_nu}\theta^{-1}(u+\mathsf A)\Big| \leq M_{1,R}, \quad u\in K_R
\]
(recall also that \( y_n\in[0,1) \)). There also exists a constant \( M_{2,R} \) such that
\[
|\theta(u+\mathsf A)| \leq M_{2,R}, \quad u\in \Big\{v-x-\mathsf By - c_\rho:~v\in K_R,~x,y\in[0,1] \Big\},
\]
by compactness of the latter set. The first inequality in \eqref{T-bounds1} now holds with \( M_R=M_{1,R}M_{2,R} \). The second inequality follows from the maximum modulus principle applied in \( \boldsymbol U_{1/R}^{(1)} \). The proof of \eqref{T-bounds2} is identical.

\subsection{Proof of Proposition~\ref{prop:lo}}
\label{ssec:lo}

Given \( \x\in\RS_{\ualpha,\ubeta} \), let \( x(\x) \), \( y(\x) \) be real numbers such that
\[
\mathfrak a(\x) = \big[c_\rho + x(\x) + \mathsf By(\x) -(\omega+\mathsf B)/2\big], \quad x(\z_0) = y(\z_0) = 0, 
\]
recall \eqref{jip}, and that continuously depend on \( \x \). Observe that these functions continuously extend to the boundary of \( \RS_{\ualpha,\ubeta} \), \( x(\x) \) is in fact continuous across \( \ualpha \) and \( x_+(\x) - x_-(\x) = 1 \) for \( \x\in\ubeta \) while \( y(\x) \) is continuous across \( \ubeta \) and satisfies \( y_+(\x) - y_-(\x) = -1 \) for \( \x\in\ualpha \), see \eqref{a-jump}. In particular, it follows from the chosen normalization, \eqref{jip}, and \eqref{xn-yn} that either \( x(\z_n) = x_n \) or \( x(\z_n) = x_n-1 \) and the same is true about \( y(\z_n) \) and \( y_n \). Define
\[
T(\z;\x) := \exp\big\{-2\pi\ic y(\x) a(\z)\big\}\dfrac{\theta \big( a(\z) - c_\rho - x(\x) - \mathsf By(\x) + \mathsf A \big)}{\theta \big( a(\z) + \mathsf A \big)}, \quad \z\in\RS_{\ualpha,\ubeta}.
\]
In particular, it follows from \eqref{RTF-per} that
\begin{equation}
\label{prop5-1}
\begin{cases}
T(\z;\z_n)=T_n(\z), & y(\z_n)=y_n, \medskip \\
T(\z;\z_n)=e^{2\pi\ic(c_\rho+x_n+\mathsf By_n-\mathsf A-\mathsf B/2)}T_n(\z), & y(\z_n)=y_n-1.
\end{cases}
\end{equation}
Furthermore, similarly to \eqref{T-jump}, it holds that
\begin{equation}
\label{prop5-2}
T_+(\s;\x) = T_-(\s;\x)\left\{
\begin{array}{rl}
\displaystyle \exp\big\{-2\pi \ic \big( c_\rho + x(\x)\big) \big\}, &  \s\in\ualpha\setminus\{\boldsymbol a_1\}, \medskip \\
\displaystyle \exp\big\{-2\pi\ic y(\x) \big\}, &  \s\in\ubeta\setminus\{\boldsymbol a_1\}.
\end{array}
\right.
\end{equation}
It is also true that \( T(\s;\x) \) is a holomorphic non-vanishing function except for a simple zero at \( \x \) and a simple pole at \( \infty^{(1)} \). Finally, with the same determination of \( (z-a_0)^{1/2} \) as in \eqref{Tn-exp} we can write
\[
T(z;\x) := t_0(\x) + t_1(\x)(z-a_0)^{1/2} + t_2(\x)(z-a_0) + \cdots
\]
for \( z \) close to \( a_0 \) and lying in the sector between \( \Delta_1 \) and \( \Delta_3 \). Notice that the coefficients \( t_l(\x) \) are continuous functions of \( \x\in\RS_{\ualpha} \) (the change of \( x(\x) \) by an integer does not affect \( T(\z;\x) \)) that continuously extend to both sides of \( \ualpha \).

Recall \eqref{a}. Given \( \x\in\RS \), let \( \x_{\pm1}\in\RS \) be the unique points such that
\begin{equation}
\label{prop5-3}
\mathfrak a( \x_{\pm1}) = \left[a(\x) \pm \big( \omega + \mathsf B\tau\big)\right].
\end{equation}
Clearly, \( \x_1 \) and \( \x_{-1} \) continuously depend on \( \x \). Notice that \( (\x_1)_{-1} = \x \) and that \( (\z_n)_{\pm1} = \z_{n\pm1} \). Define
\[
\boldsymbol T(\x) := \begin{pmatrix} 2t_2(\x) - h_1t_1(\x) & 2t_2(\x_{-1}) - h_1t_1(\x_{-1}) \medskip \\ t_0(\x) & t_0(\x_{-1}) \end{pmatrix},
\]
where the constant \( h_1 \) is the same as in \eqref{Nes}, whatever it might be. Consider the function
\[
d(\x) := \left| \det(\boldsymbol T(\x)) \right| + \left| \det(\boldsymbol T(\x_1)) \right|.
\]
As we have explained, this is a continuous function of \( \x\in\RS_{\ualpha} \) that continuously extends to both sides of \( \ualpha \). Our next goal is to prove that \( \x=\boldsymbol a_0 \) is the only possible zero of \( d(\x) \) (it might happen that this function does not vanish at all).

Assume that \( d(\x)=0 \), in which case \( \det(\boldsymbol T(\x)) = \det(\boldsymbol T(\x_1)) =0 \). This means that rows of \( \boldsymbol T(\x) \) and \( \boldsymbol T(\x_1) \) are linearly dependent and since these matrices share a column, there exists a constant \( c \) such that
\[
2t_2(\boldsymbol y) - h_1t_1(\boldsymbol y) + c t_0(\boldsymbol y) = 0, \quad \boldsymbol y\in\{\x_{-1},\x,\x_1\}.
\]
Hence, the matrix \( [t_l(\boldsymbol y)] \) has zero determinant, where \( l \) is the column index and \(  \boldsymbol y\in\{\x_{-1},\x,\x_1\} \) is the row one. By taking linear combinations of rows rather than columns, we get that there exist constants \( c_{-1},c_0,c_1 \), not all equal to zero, such that
\begin{equation}
\label{prop5-4}
c_{-1}t_l(\x_{-1}) + c_0t_l(\x) + c_1t_l(\x_1) = 0, \quad l\in\{0,1,2\}.
\end{equation}
In another connection, since the integrand in \eqref{Phi} behaves like \( \sqrt{z-a_0} \) around \( \boldsymbol a_0 \), it holds that
\begin{equation}
\label{prop5-5}
\Phi(z) = 1 + \mathcal O\left( (z-a_0)^{3/2} \right)
\end{equation}
as \( z\to a_0 \). Consider  the function 
\begin{eqnarray*}
F(\z)  &:=&  c_{-1}\frac{T(\z;\x_{-1})}{\Phi(\z)T(\z;\x)} + c_0 + c_1 \frac{\Phi(\z)T(\z;\x_1)}{T(\z;\x)} \\
&=& \frac{c_{-1}T(\z;\x_{-1}) + c_0 \Phi(\z)T(\z;\x) + c_1\Phi^2(\z)T(\z;\x_1)}{\Phi(\z)T(\z;\x)}.
\end{eqnarray*}
It follows from \eqref{Phi-jump}, \eqref{prop5-2}, and \eqref{prop5-3} that this is a rational function on the entire surface \( \RS \). Moreover, the first representation readily yields that \( F(\z) \) has exactly three poles, all simple, at \( \x \), \( \infty^{(0)} \), and \( \infty^{(1)} \), while the second formula together with  \eqref{prop5-4} and \eqref{prop5-5} clearly shows that it has at least a triple zero at \( \boldsymbol a_0 \) (unless \( \x=\boldsymbol a_0 \), in which case it is at least a double zero). Since \( F(\z) \) must have equal number of poles and zeros, its zero/pole divisor is equal to
\[
3\boldsymbol a_0 - \x - \infty^{(0)} - \infty^{(1)}.
\]
However, this means that the function \( F(\z)/(z-a_0) \) has a zero/pole divisor \( \boldsymbol a_0 - \x \). Since there are no rational functions with at most one pole besides constants, \( F(\z)=C(z-a_0) \) for some constant \( C \) and it must hold that \( \boldsymbol a_0 = \x \) as claimed.

Now, it follows from \eqref{Nes} and \eqref{prop5-1} that
\[
\big| \det(\boldsymbol T_n) \big| = Y_{n-1}Y_n\big| \det(\boldsymbol T(\z_n)) \big|, \quad Y_n = \begin{cases} 1, & y(\z_n)=y_n, \smallskip \\ e^{\pi \im( (2-2y_n +\tau) \mathsf B - 2c_\rho)}, & y(\z_n)=y_n-1. \end{cases}
\]
Since $\im(\mathsf{B}) >0$, \( \tau\in(-1,0) \), and  $ y_n\in[0,1)$, we have that $Y_n \geq Y:= e^{-\pi \im(2c_\rho+\mathsf{B})}$ for all $n \in \N$. Let $ \boldsymbol V_\varepsilon := d^{-1}([0,2  Y^{-2} \varepsilon))$. Then for all \( \varepsilon \) small enough this set is either empty (when \( d(\x) \) does not vanish) or it is a neighborhood of \( \boldsymbol a_0 \) by continuity of \( d(\x) \). Then if \( \z_n\notin \boldsymbol V_\varepsilon \), it holds that
\[
| \det (\boldsymbol T_n)| + | \det (\boldsymbol T_{n + 1})|  \geq Y^2 d(\z_n) \geq 2\varepsilon,
\]
which yields that $ | \det(\boldsymbol T_m) | \geq | \det(\boldsymbol T(\z_m)) | \geq \varepsilon $ for either $ m=n $ or $m=n+1 $.

\section{Proof of Theorem~\ref{thm:asymp1}}
\label{s:RH1}

\subsection{Initial Riemann-Hilbert Problem}
\label{ssec:IRH}

Recall the definition of \( R_n(z) \) in \eqref{linear} as well as the definition of \( f(z) \) in \eqref{CI}. It follows from the known behavior of Cauchy integrals around the endpoints of arcs of integration \cite[Sections 8.2-4]{Gakhov} (we can consider \( f(z) \) as a sum of three Cauchy integrals over the arcs \( \Delta_1,\Delta_2,\Delta_3 \)) that
\[
|f(z)|,|R_n(z)| = \mathcal O(\psi_{\alpha_l}(z-a_l)), \quad \psi_\alpha(z):=\left\{
\begin{array}{rl}
1, & \alpha>0, \smallskip \\
\log|z|, & \alpha=0, \smallskip \\
|z|^\alpha, & \alpha<0,
\end{array}
\right.
\]
for \( l\in\{0,1,2,3\} \) (it was pointed out in \cite{ApYa15} that these functions are bounded around \( a_0 \) when \( \alpha_0=0 \) and \( \rho_1(0)+\rho_2(0)+\rho_3(0)=0\), where, as before, \( \rho_i(s) \) is the restriction of \( \rho(s) \) to \( \Delta_i^\circ = \Delta_i\setminus\{a_0,a_i\} \)). Assume now that
\begin{equation}
\label{assump}
\deg(Q_n) = n \qandq R_{n-1}(z) = k_n^{-1}z^{-n} + \mathcal O\left(z^{-n-1}\right),
\end{equation}
where the second relation must hold as \( z\to\infty \) for some finite non-zero constant \( k_n \). Let
\begin{equation}
\label{Y}
\boldsymbol Y(z) := \begin{pmatrix}
Q_n(z) & R_n(z)\\
k_{n-1} Q_{n-1}(z) & k_{n-1}R_{n-1}(z)
\end{pmatrix}.
\end{equation}
Then $\boldsymbol Y(z)$ solves the following Riemann-Hilbert problem (\rhy):

\begin{itemize}
\label{rhy}
	\item[(a)] $\boldsymbol Y(z)$ is analytic in $\C \setminus \Delta$ and $\displaystyle\lim_{z \to \infty} \boldsymbol Y(z)z^{-n \sigma_3} = \boldsymbol I$, where $\sigma_3 := \left(\begin{matrix} 1 & 0 \\ 0 & -1 \end{matrix}\right)$ and \( \boldsymbol I = \left(\begin{matrix} 1 & 0 \\ 0 & 1 \end{matrix}\right)\);
	\item[(b)] $\boldsymbol Y(z)$ has continuous traces on each $\Delta_i^\circ$, \( i\in\{1,2,3\} \), that satisfy 
	\[
	\boldsymbol Y_{+}(s) = \boldsymbol Y_- (s)\left( \begin{matrix}
	1 & \rho_i(s) \\ 0 & 1
	\end{matrix} \right), \quad s\in\Delta_i^\circ;
	\]
	\item[(c)] \(\displaystyle \boldsymbol Y(z) =	\mathcal O\left(\begin{matrix}  1 & \psi_{\alpha_l}(z-a_l) \\ 1 & \psi_{\alpha_l}(z-a_l)	\end{matrix}\right) \) as \( \Delta\not\ni z\to a_l \), \( l\in\{0,1,2,3 \} \).
\end{itemize}

It is by now standard to show, see for example \cite{ApYa15,Ya15}, that if a solution of \hyperref[rhy]{\rhy} exists then it has the form \eqref{Y} and \eqref{assump} holds. Notice that \( \det(\boldsymbol Y(z)) \) is holomorphic in \( \C\setminus\Delta \), equal to \( 1 \) at infinity, is continuous across \( \Delta^\circ \) and can have at most removable singularity at each \( a_l \). Hence, \( \det(\boldsymbol Y(z))\equiv1 \).

\subsection{Opening of Lenses}
\label{ss:ol}

\begin{figure}[!ht]
\includegraphics[scale=.9]{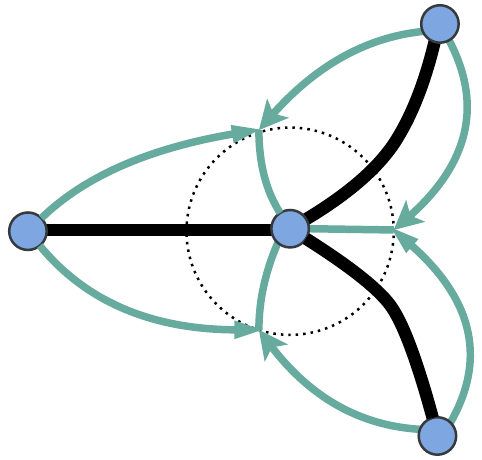}
\begin{picture}(0,0)
\put(-136,58){\(a_1\)}
\put(-40,58){\(\Delta_1^\star\)}
\put(-93,68){\(\Omega_1^+\)}
\put(-93,48){\(\Omega_1^-\)}
\put(-90,87){\(\Gamma_1^+\)}
\put(-90,26){\(\Gamma_1^-\)}
\put(-19,-4){\(a_2\)}
\put(-57,75){\(\Delta_2^\star\)}
\put(-40,25){\(\Omega_2^+\)}
\put(-25,40){\(\Omega_2^-\)}
\put(-58,10){\(\Gamma_2^+\)}
\put(-4,30){\(\Gamma_2^-\)}
\put(-19,122){\(a_3\)}
\put(-57,42){\(\Delta_3^\star\)}
\put(-40,92){\(\Omega_3^-\)}
\put(-24,77){\(\Omega_3^+\)}
\put(-58,104){\(\Gamma_3^-\)}
\put(-4,85){\(\Gamma_3^+\)}
\end{picture}
\caption{Contour \( \Gamma \). The dotted circle represents \( \partial U_0 \).}
\label{fig:2}
\end{figure}
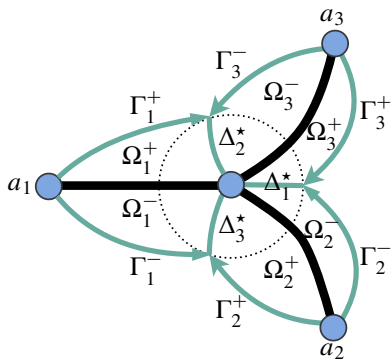

Now, we construct a system of arcs, the lens, around \( \Delta \) as on Figure~\ref{fig:2}. 
Fix \( \delta>0\) such that \( \Delta\cap U_0 \) is a smaller trefoil, where \( U_0 := \{z:|z-a_0|<\delta \} \), see the dotted circle on Figure~\ref{fig:2}. We are assuming that \( \delta \) is small enough so that \( \overline U_0 \cap (\pi(\ualpha)\cup\pi(\ubeta)) = \varnothing \). Denote by \( S_i \) the sector of \( U_0\setminus\Delta \) for which \( \partial S_i\cap\Delta_i = \{a_0\} \) (i.e., the sector opposite to \( \Delta_i \)). Further, let \( U_0^\pm \) be the connected components of \( U_0\setminus(\Delta_1\cup\Delta_1^\star) \) that border \( \Delta_1^\pm \) (\( U_0^+ \) intersects \( \Delta_3 \) and \( U_0^- \) intersects \( \Delta_2 \)), where \( \Delta_i^\star \) is an open analytic arc that splits the sector \( S_i \) into two halves (its endpoints are \( a_0 \) and some point on \( \partial S_i\setminus\Delta \)) and has the same tangent vector at \( a_0 \) as \( \Delta_i \) (we shall completely fix \( \Delta_i^\star \) later in Section~\ref{ssec:local}, see \eqref{conf-map-rays}), see Figure~\ref{fig:2}. We orient each \( \Delta_i^\star \) towards \( a_0 \). Further, let \( \Gamma_i^+ \) and \( \Gamma_i^- \) be open smooth arcs that connect \( a_i \) to \( \Delta_{i-1}^\star \) and \( \Delta_{i+1}^\star \), respectively, where the subindices are understood cyclically within \( \{1,2,3\} \), a convention we abide by throughout the rest of the paper. We assume that all these arcs do not intersect except for the common endpoints and let
\[
\Gamma := \Delta \bigcup\big(\cup_i\Delta_i^\star\big) \bigcup\big(\cup_i (\Gamma_i^+\cup\Gamma_i^-)\big).
\]
We let \( \Omega_i^\pm \) to denote bounded components of the complement of \( \Gamma \) that contain \( \Gamma_i^\pm \) as a part of their boundary and \( \Gamma^\circ \) be the subset of those points in \( \Gamma \) that possess a well defined tangent.
 
 Recall our assumption made after \eqref{CI} that \( \rho_i(s)(s-a_0)^{-\alpha_0}(s-a_i)^{-\alpha_i} \) extends to an analytic function in some neighborhood of \( \Delta_i \). We shall suppose that this neighborhood contains the closure of \( \Omega_i^+\cup\Omega_i^- \) and that the branch of \( (z-a_i)^{\alpha_i} \) is analytic throughout this closure, except, of course, at \( a_i \). Recall also that we fixed a branch of \( (z-a_0)^{\alpha_0} \) whose branch cut lies, in part, in \( S_1 \) and we shall further assume that this cut contains \( \Delta_1^\star \). Finally, recall \eqref{varrho} and suppose that \( U_0 \) is small enough so that every \( \varrho_i(z) \) is analytic in its closure.
 
Define
\begin{equation}
\label{X}
\boldsymbol X(z) := \boldsymbol Y(z)\left\{
\begin{array}{rl}
\left(\begin{matrix} 1 & 0 \smallskip \\ \mp1/\rho_i(z) & 1\end{matrix}\right), & z\in\Omega_i^\pm, \medskip \\
\boldsymbol I, & \text{otherwise}.
\end{array}
\right.
\end{equation}
Clearly, if \( \boldsymbol Y(z) \) is a solution of \hyperref[rhy]{\rhy}, then \( \boldsymbol X(z) \) is a solution of the following Riemann-Hilbert problem (\rhx):

\begin{itemize}
\label{rhx}
\item[(a)] $\boldsymbol X(z)$ is analytic in $\C \setminus \Delta$ and $\lim_{z \to \infty} \boldsymbol X(z)z^{-n \sigma_3} = \boldsymbol I$;
\item[(b)] $\boldsymbol X(z)$ has continuous traces on $\Gamma^\circ$ that satisfy 
\[
\boldsymbol X_{+}(s) = \boldsymbol X_- (s)
\left\{
\begin{array}{rl}	
\left( \begin{matrix} 0 & \rho_i(s) \\ -\rho_i^{-1}(s) & 0 \end{matrix} \right), & s\in\Delta_i^\circ, \medskip \\
\left( \begin{matrix} 1 & 0 \\ \rho_i^{-1}(s) & 1 \end{matrix} \right), & s\in\Gamma_i^\pm, \medskip \\
\left( \begin{matrix} 1 & 0 \\ \rho_{i+1+}^{-1}(s)+\rho_{i-1-}^{-1}(s) & 1 \end{matrix} \right), & s\in\Delta_i^\star,
\end{array}
\right.
\]
where taking traces is important only on \( \Delta_1^\star \);
\item[(c)] for each \( l\in\{0,1,2,3 \} \) it holds that
\[
\boldsymbol X(z) =	\mathcal O\left(\begin{matrix}  1 & |z-a_l|^{\alpha_l} \\ 1 & |z-a_l|^{\alpha_l} \end{matrix}\right) \qandq \boldsymbol X(z) =	\mathcal O\left(\begin{matrix}  \log|z-a_l| & \log|z-a_l| \\ \log|z-a_l| & \log|z-a_l| \end{matrix}\right)
\]
 as \( \Gamma\not\ni z\to a_l \) when \( \alpha_l<0 \) and \( \alpha_l=0 \), respectively, and
\[
\boldsymbol X(z) =	\mathcal O\left(\begin{matrix}  1 & 1 \\ 1 & 1 \end{matrix}\right)
\qandq
\boldsymbol X(z) =	\mathcal O\left(\begin{matrix}  |z-a_l|^{-\alpha_l} & 1 \\ |z-a_l|^{-\alpha_l} & 1 \end{matrix}\right)
\]
 as \( \Gamma\not\ni z\to a_l \), from within the unbounded and bounded components of the complement of \( \Gamma \), respectively, when \( \alpha_l>0 \).
\end{itemize}

Conversely, assume that \( \boldsymbol X(z) \) is a solution of \hyperref[rhx]{\rhx} and \( \boldsymbol Y(z) \) is defined by inverting \eqref{X}. It is quite easy to see that thus defined \( \boldsymbol Y(z) \) satisfies \hyperref[rhy]{\rhy}(a,b). Furthermore, since the second columns of \( \boldsymbol X(z) \) and \( \boldsymbol Y(z) \) coincide, the second column of \( \boldsymbol Y(z) \) satisfies \hyperref[rhy]{\rhy}(c) as well. In fact, the same holds for the first column of \( \boldsymbol Y(z) \) if \( \alpha_l<0 \). In other cases, we see from \hyperref[rhy]{\rhy}(b) that the first column of \( \boldsymbol Y(z) \) can have at most an isolated singularity at \( a_l \). However, if \( \alpha_l=0 \), then blowing up at \( a_l \) is at most logarithmic and therefore singularity must be removable. Similarly, when \( \alpha_l>0 \) and \( l\in\{1,2,3\} \), the first column of \( \boldsymbol Y(z) \) remains bounded as \( z\to a_l \) from outside the lens and therefore the singularity at \( a_l \) cannot be polar. Furthermore, in this case  \( |z-a_l|^{\alpha_l}\boldsymbol Y(z) \) is bounded around \( a_l \) and therefore the singularity cannot be essential. Thus, it is again removable. Finally, since \( a_0 \) cannot be approached from outside the lens, we rule out polar (and essential) singularity at \( a_0 \) due to the requirement \( \alpha_0<1 \).

\subsection{Parametrices}
\label{ssec:param}

Below we assume that the index \( n \) is such that \( \z_n\neq\infty^{(0)} \), see \eqref{jip}. 

The global parametrix is obtained from \hyperref[rhx]{\rhx} by ignoring local behavior around \( a_i \) and discarding the jumps on \( \Delta_i^\star \) and \( \Gamma_i^\pm \). Namely, let \( \Psi_n(\z) \) be given by \eqref{Psin}. Recall our convention that \( \Psi_n(z) \) and \( \Psi_n^*(z) \) stand for the pull-backs of \( \Psi_n(\z) \) from \( \RS^{(0)} \) and \( \RS^{(1)} \), respectively. Define
\[
\boldsymbol N(z) :=
\begin{pmatrix}
\gamma_n\Psi_n(z) & \gamma_n\Psi_n^*(z)/w(z) \smallskip \\
\gamma_{n-1}^*\Psi_{n-1}(z) & \gamma_{n-1}^*\Psi_{n-1}^*(z)/w(z) 
\end{pmatrix}
\]
for \( z\not\in\Delta \), where the constants \( \gamma_n,\gamma_{n-1}^* \) are chosen so that  $\lim_{z \to \infty} \boldsymbol N(z)z^{-n \sigma_3} = \boldsymbol I$. By the restriction placed on the indices  \( n\), these constants are finite and non-zero, see the last claim of Proposition~\ref{prop:jip}. Then \( \boldsymbol N(z) \) is a solution of the following Riemann-Hilbert problem (\rhn):

\begin{itemize}
\label{rhn}
	\item[(a)] $\boldsymbol N(z)$ is analytic in $\C \setminus \Delta$ and $\lim_{z \to \infty} \boldsymbol N(z)z^{-n \sigma_3} = \boldsymbol I$; \smallskip
	\item[(b)] $\boldsymbol N(z)$ has continuous traces on each $\Delta_i^\circ$, \( i\in\{1,2,3\} \), that satisfy 
	\[
	\boldsymbol N_+(s) = \boldsymbol N_-(s)\left( \begin{matrix}
	0 & \rho_i(s) \\ -\rho_i^{-1}(s) & 0
	\end{matrix} \right), \quad s\in\Delta_i^\circ;
	\]
\item[(c)] \(\displaystyle \boldsymbol N(z) = \mathcal O\left(\begin{matrix}  |z-a_l|^{-(2\alpha_l+1)/4} & |z-a_l|^{(2\alpha_l-1)/4} \smallskip \\ |z-a_l|^{-(2\alpha_l+1)/4} & |z-a_l|^{(2\alpha_l-1)/4}	\end{matrix}\right) \) as \( z\to a_l \), \( l\in\{0,1,2,3 \} \).
\end{itemize}
The properties of \( \boldsymbol N(z) \) described above follow easily from the requirements of \hyperref[bvp]{BVP-\(\Psi\)}.

In fact, it will be more convenient for us later to work with separate factors of \( \boldsymbol N(z) \) rather than \( \boldsymbol N(z) \) itself. To this end we write \( \boldsymbol N(z) := \boldsymbol C \boldsymbol M^\mathsf{reg}(z)\Phi^{n\sigma_3}(z)\), where
\begin{equation}
\label{Mr}
\boldsymbol C := \begin{pmatrix} \gamma_n & 0 \\ 0 & \gamma_{n-1}^* \end{pmatrix} \qandq
\displaystyle \boldsymbol M^\mathsf{reg}(z) := \begin{pmatrix} T_n(z) & T_n^*(z)/w(z) \smallskip \\ (T_{n-1}\Phi^*)(z) & (T_{n-1}^*\Phi)(z)/w(z) \end{pmatrix} S_\rho^{\sigma_3}(z),
\end{equation}
where we used the facts that \( \Phi(z)\Phi^*(z)\equiv1 \) and \( S_\rho(z)S_\rho^*(z) \equiv 1 \), see Proposition~\ref{prop:szego}.  As in the case of \( \boldsymbol Y(z) \), it is straightforward to argue that \( \det(\boldsymbol N(z))\equiv1 \) in~\( \C \). Hence, it holds that
\[
\det(\boldsymbol M^\mathsf{reg}(z))\equiv 1/\det(\boldsymbol C) = 1/(\gamma_n\gamma_{n-1}^*).
\]

Let \( \N_\varepsilon^\mathsf{reg} \) be as in \eqref{Ne}. Then it follows from the above identity, \eqref{Mr}, and \eqref{Psin} that
\begin{equation}
\label{curious0}
(\gamma_n\gamma_{n-1}^*)^{-1} = T_n(\infty)\lim_{z\to\infty}\frac{T_{n-1}^*(z)\Phi(z)}{w(z)}.
\end{equation}
It is well known (and is not very important for us here) that \( \lim_{z\to\infty}|\Phi(z)/z| = 1/\mathrm{cap}(\Delta) \) is the reciprocal of the logarithmic capacity of \( \Delta \) (in any case, it is a non-zero and finite constant). Hence, it follows from \eqref{T-bounds1}--\eqref{T-bounds2} that
\begin{equation}
\label{curious1}
\mathrm{cap}(\Delta) RM_R^{-2} \leq  |\gamma_n\gamma_{n-1}^*| \leq \mathrm{cap}(\Delta)M_\varepsilon^2, \quad n\in\N_\varepsilon^\mathsf{reg},
\end{equation}
for any given \( R>\max_{s\in\Delta}|s| \). We further get from \eqref{T-bounds1}--\eqref{T-bounds2} and the above inequalities that 
\begin{equation}
\label{Mr-bound}
\begin{cases}
\boldsymbol M^\mathsf{reg}(z) &= \boldsymbol{\mathcal O}_{K,\varepsilon}(1)\diag(1,1/w(z))S_\rho^{\sigma_3}(z) \smallskip \\ \boldsymbol M^\mathsf{reg}(z)^{-1} &= S_\rho^{-\sigma_3}(z)\diag(1/w(z),1)\boldsymbol{\mathcal O}_{K,\varepsilon}(1)
\end{cases}
\end{equation}
on any compact set \( K \) as \( \N_\varepsilon^\mathsf{reg}\ni n\to\infty \).

Assume now that \( \delta>0 \) is small enough so that \( U_i :=\{z:|z-a_i|<\delta\} \), \( i\in\{1,2,3\} \), intersects only \( \Delta_i \), i.e., \( U_i\cap\Delta_j = \varnothing \) for \( j\neq i \), and that \( \rho_i(s)(s-a_i)^{-\alpha_i} \) extends analytically to its closure, which is disjoint from the closure of \( U_0 \). Local parametrix around \( a_i \) is obtained by solving \hyperref[rhx]{\rhx} in \( U_i \). That is, we need to solve the following Riemann-Hilbert problem (\rhpi):

\begin{itemize}
\label{rhpi}
\item[(a)] $\boldsymbol P_i(z)$ is analytic in $U_i \setminus \Gamma$; \smallskip
\item[(b,c)] $\boldsymbol P_i(z)$ satisfies \hyperref[rhx]{\rhx}(b,c) within \( U_i \); \smallskip
\item[(d)] it holds that \( \boldsymbol P_i(s) = (\boldsymbol I+\boldsymbol{\mathcal O}(1/n)) \boldsymbol M^\mathsf{reg}(s) \Phi^{n\sigma_3}(s) \) uniformly on \( \partial U_i \).
\end{itemize}

The solution of \hyperref[rhpi]{\rhpi} is well known and was first derived in \cite{KMcLVAV04} in the context of orthogonal polynomials on a segment with the help of the modified Bessel and Hankel functions. Since then it appeared in countless papers and was adopted to the context of non-Hermitian orthogonal polynomials on symmetric contours in \cite{ApYa15}. Therefore, we shall skip the explicit construction.

Around \( a_0 \), we look only for an approximate parametrix. More precisely, we would like to solve the following Riemann-Hilbert problem (\rhpo):

\begin{itemize}
\label{rhpo}
\item[(a)] $\boldsymbol P_0(z)$ is analytic in $U_0 \setminus \Gamma$; \smallskip
\item[(b)] $\boldsymbol P_0(z)$ satisfies \hyperref[rhx]{\rhx}(b) within \( U_0 \) except the \(21\)-entries of the jump matrices on \( \Delta_i^\star \), \( i\in\{1,2,3\} \), are replaced by
\[
\kappa_{i,i-1}(s)\rho_{i+1+}^{-1}(s)+\kappa_{i,i+1}(s)\rho_{i-1-}^{-1}(s),
\]
where \( \kappa_{i,j}(z) := (\varrho_i(z)\varrho_j(a_0))/(\varrho_i(a_0)\varrho_j(z)) \) (recall \eqref{varrho}); \smallskip
\item[(c)] $\boldsymbol P_0(z)$ satisfies \hyperref[rhx]{\rhx}(c) within \( U_0 \); \smallskip
\item[(d)] it holds that \( \boldsymbol P_0(s) = \big(\boldsymbol I+\boldsymbol{\mathcal O}\big( n^{-1/3}\big)\big) \boldsymbol M^\mathsf{reg}(s) \Phi^{n\sigma_3}(s) \) uniformly on \( \partial U_0 \).
\end{itemize}

We solve \hyperref[rhpo]{\rhpo} in the next subsection.  The presence of the functions \( \kappa_{i,j}(z) \) is precisely what makes this parametrix approximate. {Of course, it would be preferable to solve \hyperref[rhx]{\rhx} locally around \( a_0 \) exactly, that is, with the jumps appearing in \hyperref[rhx]{\rhx}(b). However, currently, we do not know how to achieve this. Notice also that if \( \varrho_i(s)=c_i\varrho(s) \) for some function \( \varrho(z) \) analytic around \( a_0 \) and possibly different constants \( c_i \), which is exactly the case when the approximated function has the form \eqref{log-power}, see Appendix~\ref{ap:functions}, then \( \kappa_{i,j}(z)\equiv1 \) and the parametrix is exact.}

\subsection{Approximate Local Parametrix around \( a_0 \)}
\label{ssec:local}

Solution of \hyperref[rhpo]{\rhpo} is based on the solution of \hyperref[rhphi]{\rhphi} from Appendix~\ref{ap:34}. The knowledge of the statement of that Riemann-Hilbert problem is sufficient for understanding the construction further below.

Recall the definition of \( \Sigma_2,\Sigma_3 \) in Appendix~\ref{ap:34}. Define
\begin{equation}
\label{conf-map}
\zeta^3(z) := \left(\frac32\int_{a_0}^z\frac{(s-a_0)}{w(s)}\dd s\right)^2, \quad z\in \big(S_1\cup S_2\cup S_3\big)\cap U_0.
\end{equation}
Since \( w_+(s)=-w_-(s) \) for \( s\in\Delta \), \( \zeta^3(z) \) extends to a holomorphic function in \( U_0 \), which, upon easy verification, has a triple zero at \( a_0 \) (we already used this observation in \eqref{prop5-5}). Furthermore, it readily follows from \eqref{differential} that \( \zeta^3(s) \) is negative for \( s\in\Delta\cap U_0 \) (except for \( a_0 \), of course). Hence, we can select a branch \( \zeta(z) \) which is holomorphic in \(U_0 \) and such that \( \zeta(s)<0 \) when \( s\in\Delta_1\cap U_0 \). Since \( \zeta(z) \) has a simple zero at \( a_0 \), we can decrease the radius of \( U_0 \) if necessary so that \( \zeta(z) \) is univalent in \( U_0 \). Moreover, as we had some freedom in choosing the arcs \( \Delta_i^\star \), we now fix them so that 
\begin{equation}
\label{conf-map-rays}
\zeta(\Delta_1^\star)\subset(0,\infty) \qandq \zeta(\Delta_i^\star)\subset \Sigma_i, \quad i\in\{2,3\}.
\end{equation}
In what follows, all the roots of \( \zeta(z) \) are principal, that is, they remain positive on \( \Delta_1^\star \) and have a branch cut along \( \Delta_1 \). In particular,
\begin{equation}
\label{conf-quarter}
\zeta_+^{1/4}(s) = \ic\zeta_-^{1/4}(s), \quad s\in\Delta_1\cap U_0.
\end{equation}
Recall the definition of \( \Phi(\z) \) in \eqref{Phi} and that \( |\Phi(z)|>1 \) for \( z\in U_0\setminus\Delta \), see \eqref{Phi-Green}.
Thus, it also holds that 
\begin{equation}
\label{conf-map-Phi}
e^{-(2/3)\zeta^{3/2}(z)} = \left\{\begin{array}{rl}
\Phi^*(z), & z\in S_1\cap U_0, \smallskip \\
\Phi(z), & z\in (S_2\cup S_3 )\cap U_0.
\end{array}
\right.
\end{equation}

Let \( \varrho_i(s) \) be given by \eqref{varrho}, where for definiteness we fix a branch of \( (z-a_0)^{\alpha_0} \) that is analytic in \( U_0\setminus\overline{\Delta_1^\star} \). Define
\[
E_{\alpha_0}(z) :=  \left\{
\begin{array}{ll}
e^{\pm\pi\ic\alpha_0/2}, & z\in S_1\cap U_0^\pm, \smallskip \\
e^{\mp\pi\ic\alpha_0/2}, & z\in (S_2\cup S_3 )\cap U_0^\pm.
\end{array}
\right.
\]
Further, let \( r_i(z) \), \( i\in\{1,2,3\} \),  be a function holomorphic in \( S_i \) given by
\begin{equation}
\label{ri}
r_i(z) := \varsigma_i\ic (z-a_0)^{\alpha_0/2}E_{\alpha_0}(z)\left(\frac{\varrho_{i-1}(z)\varrho_{i+1}(z)}{\varrho_i(z)}\right)^{1/2},
\end{equation}
where \( \varsigma_1=\varsigma_2=-\varsigma_3=1 \) and we use the same square root branches of \( \varrho_j(z) \) in all \( r_i(z) \). Since \( E_{\alpha_0+}(s)E_{\alpha_0-}(s)\equiv1 \) for \( s\in\Delta_i^\circ\cap U_0 \), it holds that
\begin{equation}
\label{ri-Di}
r_{i-1-}(s)r_{i+1+}(s) = \varsigma_i\rho_i(s), \quad s\in\Delta_i^\circ\cap U_0.
\end{equation}
We further set \( r(z):=r_i(z) \) for \( z\in S_i \).

Recall that \( \rho\in W^\mathsf{reg} \). Let \( \boldsymbol\Phi_{\alpha_0}(\zeta;0) \) be the solution of \hyperref[rhphi]{\rhphi} with \( \alpha=\alpha_0 \) and the Stokes parameters given by \eqref{stokes}, see Appendix~\ref{ap:34} (it exists by the very definition of the class \( W^\mathsf{reg} \)). Then a solution of \hyperref[rhpo]{\rhpo} is given by
\begin{equation}
\label{P0}
\boldsymbol P_0(z) := \boldsymbol E_0^\mathsf{reg}(z)\boldsymbol\Phi_{\alpha_0}\big(n^{2/3}\zeta(z);0\big)\boldsymbol J(z)r^{-\sigma_3}(z),
\end{equation}
where \( \boldsymbol E_0^\mathsf{reg}(z) \) is a holomorphic matrix function in \( U_0 \) that we shall specify further below in \eqref{E0} and
\begin{equation}
\label{J}
\boldsymbol J(z) :=  \left\{
\begin{array}{rl}
\left(\begin{matrix} 0 & -1 \\ 1 & 0 \end{matrix}\right), & z\in S_1\cap U_0, \smallskip \\
\boldsymbol I, & z\in (S_2\cup S_3 )\cap U_0.
\end{array}
\right.
\end{equation}
Indeed, it follows from \eqref{conf-map-rays} that \( \boldsymbol P_0(z) \) is a holomorphic matrix function in \( U_0\setminus\Gamma \), i.e., \hyperref[rhpo]{\rhpo}(a) is satisfied. We further get from \hyperref[rhphi]{\rhphi}(b) and \eqref{ri-Di} that
\[
\boldsymbol P_{0-}^{-1}(s)\boldsymbol P_{0+}(s) = r_{3-}^{\sigma_3}(s) \left(\begin{matrix} 0 & 1 \\ -1 & 0 \end{matrix}\right) r_{2+}^{-\sigma_3}(s) = \left(\begin{matrix} 0 & \rho_1(s) \\ -1/\rho_1(s) & 0 \end{matrix}\right)
\]
for \( s\in\Delta_1^\circ\cap U_0 \),
\[
\boldsymbol P_{0-}^{-1}(s)\boldsymbol P_{0+}(s) = r_{1-}^{\sigma_3}(s)\boldsymbol J_-^{-1}(s)r_{3+}^{-\sigma_3}(s) = \left(\begin{matrix} 0 & \rho_2(s) \\ -1/\rho_2(s) & 0 \end{matrix}\right)
\]
for \( s\in\Delta_2^\circ\cap U_0 \), and
\[
\boldsymbol P_{0-}^{-1}(s)\boldsymbol P_{0+}(s) = r_{2-}^{\sigma_3}(s)\boldsymbol J_+(s)r_{1+}^{-\sigma_3}(s) = \left(\begin{matrix} 0 & \rho_3(s) \\ -1/\rho_3(s) & 0 \end{matrix}\right)
\]
for \( s\in\Delta_3^\circ\cap U_0 \). Moreover, using \eqref{stokes} we get that
\[
\boldsymbol P_{0-}^{-1}(s)\boldsymbol P_{0+}(s)= r_1^{\sigma_3}(s)\boldsymbol J^{-1}(s) \left(\begin{matrix} 1 & -b_1(\rho) \\ 0 & 1 \end{matrix}\right) \boldsymbol J(s)r_1^{-\sigma_3}(s) = \left(\begin{matrix} 1 & 0 \\ \frac{\kappa_{1,2}(s)}{\rho_{2+}(s)}+\frac{\kappa_{1,3}(s)}{\rho_{3-}(s)} & 1 \end{matrix}\right)
\]
for \( s\in\Delta_1^\star\cap U_0 \) (notice that the orientations of \( \zeta(\Delta_1^\star) \) and \( \Sigma_1 \) are opposite to each other), and
\[
\boldsymbol P_{0-}^{-1}(s)\boldsymbol P_{0+}(s) = r_l^{\sigma_3}(s)\left(\begin{matrix} 1 & 0 \\ b_l(\rho) & 1 \end{matrix}\right) r_l^{-\sigma_3}(s) = \left(\begin{matrix} 1 & 0 \\ \frac{\kappa_{l,l-1}(s)}{\rho_{l+1}(s)}+\frac{\kappa_{l,l+1}(s)}{\rho_{l-1}(s)} & 1 \end{matrix}\right)
\]
for \( s\in\Delta_l^\star\cap U_0 \), \( l\in\{2,3\} \). Hence, \( \boldsymbol P_0(z) \), defined in \eqref{P0}, fulfills \hyperref[rhpo]{\rhpo}(b). Next, it follows from \hyperref[rhphi]{\rhphi}(c), \eqref{P0}, and \eqref{ri} that
\[
\boldsymbol P_0(z) = \boldsymbol E_0^\mathsf{reg}(z)\left\{
\begin{array}{rl}
\boldsymbol{\mathcal O}\left(|\zeta(z)|^{-|\alpha_0|/2}\right), & \alpha_0\neq 0 \medskip \\
\boldsymbol{\mathcal O}\left(\log|\zeta(z)|\right), & \alpha_0= 0
\end{array}
\right\}\mathcal O\left(|z-a_0|^{-\alpha_0\sigma_3/2} \right).
\]
Since multiplication by \( \boldsymbol E_0^\mathsf{reg}(z) \) on the left does not mix the entries from different columns of the product of the other two factors above, \( \boldsymbol P_0(z) \) does indeed satisfy \hyperref[rhpo]{\rhpo}(c).

Finally, let
\begin{equation}
\label{E0}
\boldsymbol E_0^\mathsf{reg}(z) := \boldsymbol M^\mathsf{reg}(z)r^{\sigma_3}(z) \boldsymbol J^{-1}(z) \left(\begin{matrix} 1 & -\ic \\ -\ic & 1 \end{matrix}\right) \frac{\big(n^{2/3}\zeta(z)\big)^{\sigma_3/4}}{\sqrt2}.
\end{equation}
We get from \hyperref[rhphi]{\rhphi}(d), \eqref{Mr-bound}, \eqref{conf-map-Phi}, and \eqref{J} that
\begin{eqnarray}
\boldsymbol P_0(s) &=& \boldsymbol M^\mathsf{reg}(s) r^{\sigma_3}(s)  \left( \boldsymbol I + \boldsymbol{\mathcal O}\big(n^{-1/3}\big)\right) r^{-\sigma_3}(s) \Phi^{n\sigma_3}(s) \nonumber \\
\label{P0d}
& = &  \left( \boldsymbol I + \boldsymbol{\mathcal O}\big(n^{-1/3}\big) \right) \boldsymbol M^\mathsf{reg}(s) \Phi^{n\sigma_3}(s),
\end{eqnarray}
which is exactly what was required in \hyperref[rhpo]{\rhpo}(d). Thus, it only remains to prove that \( \boldsymbol E_0^\mathsf{reg}(z) \) is analytic in \( U_0 \). To this end, this matrix is clearly analytic in \( U_0\setminus\Delta \). It also holds that
\begin{eqnarray*}
(\boldsymbol E_{0-}^\mathsf{reg})^{-1}(s)\boldsymbol E_{0+}^\mathsf{reg}(s) &=& n^{-\sigma_3/6}\frac{\zeta_-^{-\sigma_3/4}(s)}{\sqrt2} \left(\begin{matrix} 1 & \ic \\ \ic & 1 \end{matrix}\right) \left(\begin{matrix} 0 & 1 \\ -1 & 0 \end{matrix}\right) \left(\begin{matrix} 1 & -\ic \\ -\ic & 1 \end{matrix}\right) \frac{\zeta_+^{\sigma_3/4}(s)}{\sqrt2}n^{\sigma_3/6} \\ 
&=& n^{-\sigma_3/6}\zeta_-^{-\sigma_3/4}(s) \left(\begin{matrix} -\ic & 0 \\ 0 & \ic \end{matrix}\right)\zeta_+^{\sigma_3/4}(s)n^{\sigma_3/6} = \boldsymbol I
\end{eqnarray*}
for \( s\in\Delta_1^\circ\cap U_0 \), where we used \hyperref[rhn]{\rhn}(b), \eqref{ri-Di}, \eqref{J}, and \eqref{conf-quarter}. Similarly, we have for \( s\in\Delta_l^\circ\cap U_0 \), \( l\in\{2,3\} \), that \( (\boldsymbol E_{0-}^\mathsf{reg})^{-1}(s)\boldsymbol E_{0+}^\mathsf{reg}(s)=\boldsymbol I \), because it holds there that
\[
\boldsymbol J_-(s)\left(\begin{matrix} 0 & (-1)^l \\ (-1)^{l-1} & 0 \end{matrix}\right)\boldsymbol J_+^{-1}(s)  = \boldsymbol I.
\]
Hence, \( \boldsymbol E_0^\mathsf{reg}(z) \) is analytic in \( U_0\setminus\{a_0\} \). Furthermore, we get from \hyperref[rhn]{\rhn}(c) and \eqref{ri} that
\[
 \boldsymbol E_0^\mathsf{reg}(z) = \boldsymbol{\mathcal O}\left(|z-a_0|^{-1/4}\right)\boldsymbol J^{-1}(s) \left(\begin{matrix} 1 & -\ic \\ -\ic & 1 \end{matrix}\right) \frac{\zeta^{\sigma_3/4}(z)}{\sqrt2}n^{\sigma_3/6} = \boldsymbol{\mathcal O}\left(|z-a_0|^{-1/2}\right).
\]
Therefore, \( a_0 \) cannot be a polar singularity of \( \boldsymbol E_0^\mathsf{reg}(z) \) and must be a point of analyticity. 

\subsection{Small-Norm Riemann-Hilbert Problem}

Let \( \boldsymbol C \) and \( \boldsymbol M^\mathsf{reg}(z) \) be matrices defined in \eqref{Mr} and \( \boldsymbol P_l(z) \), \( l\in\{0,1,2,3\} \), be the matrix functions solving \hyperref[rhpo]{\rhpo} and \hyperref[rhpi]{\rhpi}, \( i\in\{1,2,3\} \). Further, let \(
U := U_0 \cup U_1\cup U_2 \cup U_3 \), whose boundary we orient clockwise. 
\begin{figure}[!ht]
\includegraphics[scale=.9]{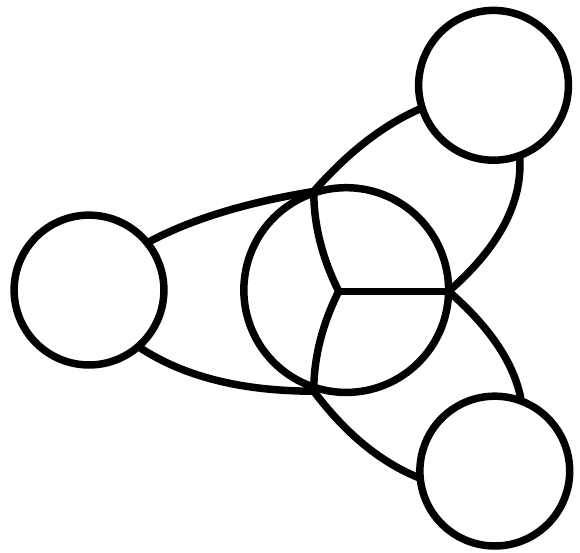}
\begin{picture}(0,0)
\put(-148,65){\(\partial U_1\)}
\put(-30,8){\(\partial U_2\)}
\put(-30,131){\(\partial U_3\)}
\put(-88,65){\(\partial U_0\)}
\end{picture}
\caption{Contour \( \Sigma \) and the circles \( \partial U_l \).}
\label{fig:3}
\end{figure}
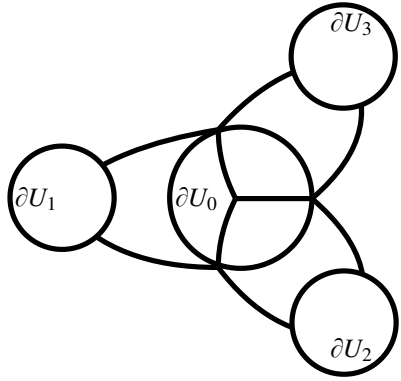
We are looking for the solution of \hyperref[rhx]{\rhx} in the form
\begin{equation}
\label{Z}
\boldsymbol X(z) = \boldsymbol C \boldsymbol Z(z) \left\{ 
\begin{array}{rl}
\boldsymbol M^\mathsf{reg}(z)\Phi^{n\sigma_3}, & z\not\in \left( \Gamma \cup \overline U \right), \medskip \\
\boldsymbol P_l(z), & z\in U_l\setminus\Gamma, \quad l\in\{0,1,2,3\}.
\end{array}
\right.
\end{equation}

Let \( \Sigma := \partial U \cup \big( \Gamma\setminus(\Delta^\circ\cup U_1\cup U_2\cup U_3)\big) \), see Figure~\ref{fig:3}. Combining the above equation with \hyperref[rhx]{\rhx}, we see that \( \boldsymbol Z(z) \) should solve the following Riemann-Hilbert problem (\rhz):

\begin{itemize}
\label{rhz}
\item[(a)] $\boldsymbol Z(z)$ is analytic in $\C \setminus \Sigma$ and $\lim_{z \to \infty} \boldsymbol Z(z) = \boldsymbol I$;
\item[(b)] $\boldsymbol Z(z)$ has continuous traces on $\Sigma^\circ$ that satisfy 
\[
\boldsymbol Z_+(s) = \boldsymbol Z_- (s)
\left\{
\begin{array}{rl}	
\boldsymbol P_l(s)\boldsymbol N_*^{-1}(s), & s\in \partial U, \medskip \\
\boldsymbol N_*(s)\left( \begin{matrix} 1 & 0 \\ \rho_i^{-1}(s) & 1 \end{matrix} \right)\boldsymbol N_*^{-1}(s), & s\in\Gamma_i^\pm\setminus\overline U,
\end{array}
\right.
\]
where \( \boldsymbol N_*(z) := \boldsymbol C^{-1}\boldsymbol N^\mathsf{reg}(z) = \boldsymbol M^\mathsf{reg}(z)\Phi^{n\sigma_3}(z) \) and we understand that \( \boldsymbol P_l(s) \) is used when \( s\in\partial U_l \), \( l\in\{0,1,2,3\} \), in the first jump relation, and
\[
\boldsymbol Z_+(s) = \boldsymbol Z_- (s) \boldsymbol P_{0-}(s)\left( \begin{matrix} 1 & 0 \\ \rho_{i+1+}^{-1}(s)+\rho_{i-1-}^{-1}(s) & 1 \end{matrix} \right)\boldsymbol P_{0+}^{-1}(s), \quad  s\in\Delta_i^\star.
\]
\end{itemize}

Note that if the jump relations \hyperref[rhpo]{\rhpo}(b) on \( \Delta_i^\star \) are exactly the same as the ones in \hyperref[rhx]{\rhx}(b), then the jump relations on \( \Delta_i^\star \cap U_0 \) in \hyperref[rhz]{\rhz}(b) above are absent.

As usual in this line of arguments, we shall show that all the jump matrices in \hyperref[rhz]{\rhz}(b) are uniformly close to \( \boldsymbol I \) along any sequence \( \N_\varepsilon^\mathsf{reg} \). Indeed, it follows from \hyperref[rhpo]{\rhpo}(d) and \hyperref[rhpi]{\rhpi}(d), \( i\in\{1,2,3\} \), that the jump of \( \boldsymbol Z(z) \) on \( \partial U \) can be estimated as
\begin{equation}
\label{R-jump-1}
 \boldsymbol I + \boldsymbol{\mathcal O}_\varepsilon\big(n^{-1/3}\big)  \quad \text{on} \quad \partial U,
\end{equation}
where  the constants in \(\boldsymbol{\mathcal O}_\varepsilon(\cdot) \) are independent of \( n \) but do depend on \( \varepsilon \). Furthermore, we get from \eqref{Mr-bound} as well as \eqref{Phi-Green} that the jump of \( \boldsymbol Z(z) \) on \( \Gamma_i^\pm\setminus \overline U \) satisfies
\begin{equation}
\label{R-jump-2}
\boldsymbol I + \rho_i^{-1}(s)\Phi^{-2n}(s)\boldsymbol M^\mathsf{reg}(s) \left( \begin{matrix} 0 & 0 \\ 1 & 0 \end{matrix} \right) (\boldsymbol M^\mathsf{reg})^{-1}(s) = \boldsymbol I + \boldsymbol{\mathcal O}_\varepsilon\big(e^{-cn}\big)
\end{equation}
for some constant \( c>0 \), where we also used  the fact that the arcs \( \Gamma_i^\pm \) lie fixed distance away from \( \Delta \) (again, the constant in \(\boldsymbol{\mathcal O}_\varepsilon(\cdot) \) is independent of \( n \) but depends on \( \varepsilon \)). 

{Now, if it holds that \( \kappa_{i,j}(z) \equiv 1 \), then \( \boldsymbol P_0(z) \) is an exact parametrix and \( \boldsymbol Z(z) \) does not have jumps on \( \Delta_i^\star \). Below, we assume that all of these functions are not identically \( 1 \), the case where only some of them are can be considered analogously. Let \( p \) be the largest integer such that \( |1-\kappa_{i,j}(s)| = \mathcal O( |s-a_0|^p )\) as \( s\to a_0 \) for all the pairs \( i,j \). It necessarily holds that \( p\geq1 \).}  It follows from \hyperref[rhpo]{\rhpo}(b) that the jump of \( \boldsymbol Z(z) \) on \( \Delta_i^\star \) can be written as
\begin{multline}
\label{R-jump-5}
\boldsymbol I + \left( \frac{1-\kappa_{i,i-1}(s)}{\rho_{i+1+}(s)}+ \frac{1-\kappa_{i,i+1}(s)}{\rho_{i-1-}(s)}\right)\boldsymbol P_{0+}(s)  \begin{pmatrix} 0 & 0 \\ 1 & 0 \end{pmatrix} \boldsymbol P_{0+}^{-1}(s) \\ = \boldsymbol I + \mathcal O\big(|s-a_0|^p\big) \boldsymbol E_0^\mathsf{reg}(s) \boldsymbol \Phi_{\alpha_0}\left(n^{2/3}\zeta(s);0\right) \boldsymbol J \begin{pmatrix} 0 & 0 \\ 1 & 0 \end{pmatrix} \boldsymbol J^{-1} \boldsymbol \Phi_{\alpha_0}^{-1}\left(n^{2/3}\zeta(s);0\right) \boldsymbol E_0^\mathsf{reg}(s)^{-1},
\end{multline}
where we used estimates \( |\rho_i(s)| \sim |s-a_0|^{\alpha_0}\sim |r(s)|^2 \), see \eqref{ri} and \hyperref[rhpo]{\rhpo}(b), as well as \eqref{P0} (notice also that we do not need to take boundary values as the relevant entries of \( \boldsymbol \Phi_{\alpha_0} \) are analytic across the considered arcs). {Let \( Z=n^{2/3}\zeta(s) \). We shall estimate the expression in \eqref{R-jump-5} separately in two regimes: when \( |Z|>1 \) and when \( |Z|\leq1 \).} Uniformity of the asymptotics in \hyperref[rhphi]{\rhphi}(d) and \eqref{J} yield that
\begin{equation}
\label{Phi-a-est1}
\boldsymbol \Phi_{\alpha_0}\left(Z;0\right) \boldsymbol J\left( \begin{matrix} 0 & 0 \\ 1 & 0 \end{matrix} \right) \boldsymbol J^{-1} \boldsymbol \Phi_{\alpha_0}^{-1}\left(Z;0\right) = \mathcal O\left( |Z|^{1/2} e^{ -(4/3)|Z|^{3/2} } \right) {= \mathcal O\big( |Z|^{-p} \big)}
\end{equation}
uniformly for \( |Z|>1 \), where one needs to observe that \( Z^{3/2}>0 \) for \( s\in\Delta_1^\star \) and \( Z^{3/2}<0 \) for \( s \in \Delta_2^\star\cup\Delta_3^\star \), see \eqref{conf-map-rays}  (recall also that \( \det(\boldsymbol\Phi_{\alpha_0}(Z;0)) \equiv1\)), {and notice that \( |Z|^{p+1/2}e^{-(4/3)|Z|^{3/2}} \) is a bounded above positive function of \( |Z| \) on \( [1,\infty) \)}. On the other hand, we get from \hyperref[rhphi]{\rhphi}(c) that
\begin{multline}
\label{Phi-a-est2}
\boldsymbol \Phi_{\alpha_0}\left(Z;0\right) \boldsymbol J\left( \begin{matrix} 0 & 0 \\ 1 & 0 \end{matrix} \right) \boldsymbol J^{-1} \boldsymbol \Phi_{\alpha_0}^{-1}\left(Z;0\right) = \\ 
\begin{cases} 
{\mathcal O}(1)|Z|^{\alpha_0\sigma_3/2}{\mathcal O}(1)|Z|^{-\alpha_0\sigma_3/2}\mathcal O(1) , & \alpha_0\neq 0, \\
{\mathcal O}(1)\mathcal O(\log|Z|){\mathcal O}(1)\mathcal O(\log|Z|)\mathcal O(1) , & \alpha_0 = 0,
\end{cases}
=
\begin{cases}
\mathcal O\left(|Z|^{-|\alpha_0|}\right), & \alpha_0\neq 0, \medskip \\
\mathcal O\left(\log^2|Z|\right), & \alpha_0=0,
\end{cases}
{ = \mathcal O\big( |Z|^{-p} \big)}
\end{multline}
uniformly for \( |Z| \leq 1 \) {because \( |Z|^{p-|\alpha_0|} \) and \( |Z|^p\log^2|Z| \) are bounded positive functions of \( |Z| \) on \( [0,1] \) since \( p\geq 1 \) and \( |\alpha_0|\leq 1 \).} Recall that \( |z-a_0|/|\zeta(z)| \) is bounded in \( U_0 \) and \( \boldsymbol E_0^\mathsf{reg}(z) \) is an analytic matrix function with the same determinant as \( \boldsymbol M^\mathsf{reg}(z) \). {Definition~\eqref{E0}, estimates \eqref{Mr-bound}, and the maximum modulus principle yield that \( \boldsymbol E_0^\mathsf{reg}(z) = \boldsymbol{\mathcal O}_\varepsilon(1)n^{\sigma_3/6}\) and \( \boldsymbol E_0^\mathsf{reg}(z)^{-1} = \boldsymbol{\mathcal O}_\varepsilon(1)n^{-\sigma_3/6}\).} It now follows from \eqref{R-jump-5}--\eqref{Phi-a-est2} that the jump of \( \boldsymbol Z(z) \) on \( \cup\Delta_i^\star \) can be estimated as
\begin{equation}
\label{R-jump-4}
{\boldsymbol I + \boldsymbol{\mathcal O}_\varepsilon\bigg(n^{1/3}\max_{s\in \cup\Delta_i^\star} |s-a_0|^p|Z|^{-p}\bigg) = \boldsymbol I + \boldsymbol{\mathcal O}_\varepsilon\big(n^{-(2p-1)/3}\big),}
\end{equation}
where \( Z=n^{2/3}\zeta(s) \). {As \( p\geq 1 \), the error rate in \eqref{R-jump-4} is not worse than \( \boldsymbol{\mathcal O}_\varepsilon(n^{-1/3})\)}.

Finally, by arguing as in \cite[Theorem~7.103 and Corollary~7.108]{Deift}, see also \cite[Theorem~8.1]{FokasItsKapaevNovokshenov}, we obtain from \eqref{R-jump-1}, \eqref{R-jump-2}, and \eqref{R-jump-4} that the matrix \( \boldsymbol Z(z) \) exists for all \( n\in\N_\varepsilon^\mathsf{reg} \) large enough and that
\[
\|\boldsymbol Z_\pm-\boldsymbol I\|_{2,\Sigma} = \mathcal O_\varepsilon \big( n^{-1/3}\big).
\]
Since the jumps of \( \boldsymbol Z(z) \) on \( \Sigma \) are restrictions of holomorphic matrix functions, the standard deformation of the contour technique and the above estimate yield that
\begin{equation}
\label{z-rate}
\boldsymbol Z = \boldsymbol I + \boldsymbol{\mathcal{O}}_{\delta,\varepsilon} \big( n^{-1/3}\big) \quad \text{uniformly in} \quad \overline\C\setminus\{a_0\}. 
\end{equation}

\subsection{Proof of Theorem~\ref{thm:asymp1}}

Let \( \boldsymbol Z(z) \) be a solution of \hyperref[rhz]{\rhz}, in which case \eqref{Z} holds. Given a closed set \( K\subset\overline\C\setminus\Delta  \), the contour \( \Sigma \) can always be adjusted so that \( K \) lies in the exterior domain of \( \Sigma \). Then it follows from \eqref{X} that \( \boldsymbol Y(z)=\boldsymbol X(z) \) on \( K \). Formulae \eqref{Asymp1} and \eqref{upsilons} now follow immediately from \eqref{Y}, \eqref{Mr} and \eqref{z-rate} since
\[
w^{i-1}(z)\left[\left(\boldsymbol{ZN}_*\right)(z)\right]_{1i} = (1 + \upsilon_{n1}(z))\Psi_n\big(z^{(i-1)}\big) + \upsilon_{n2}(z)\Psi_{n-1}\big(z^{(i-1)}\big)
\]
for \( z\in K \), where \( 1+\upsilon_{n1}(z),  \upsilon_{n2}(z) \) are the first row entries of \( \boldsymbol Z(z) \).  Similarly, if \( K \) is a compact subset of \( \Delta^\circ \), the lens \( \Sigma \) can be arranged so that \( K \) does not intersect \( \overline U \). As before, we get from \eqref{Y}, \eqref{X}, and \eqref{Z} that
\begin{multline}
\label{Qn-Delta}
\gamma_n^{-1}Q_n(z) =  \Big(\big(1 + \upsilon_{n1}(z) \big)\Psi_n(z) + \upsilon_{n2}(z)\Psi_{n-1}(z)\Big) \pm \\ (\rho_i w)^{-1}(z) \Big(\big(1 + \upsilon_{n1}(z)\big)\Psi_n^*(z) + \upsilon_{n2}(z)\Psi_{n-1}^*(z)\Big)
\end{multline}
for \( z\in \Omega_{i\pm} \setminus \overline U \). The top formula in \eqref{Asymp2} now follows by taking the trace of the right-hand side of the above equality on \( K \) and using the top relation in \eqref{BVP-Psi} (which yields that \( \Psi_{n\pm}^*(s)=\Psi_{n\mp}(s)(\rho w_+)(s) \) for \( s\in\Delta^\circ \)). Since \( R_n(z) = [\boldsymbol X]_{12}(z)\) for \( z\in \Omega_{i\pm} \setminus \overline U \), the bottom relation in \eqref{Asymp2} is even simpler to derive.

\section{Proof of Theorem~\ref{thm:asymp2}}
\label{s:RH2}

The main difference between proofs of Theorem~\ref{thm:asymp1} and \ref{thm:asymp2} is that we no longer can use the matrix \( \boldsymbol\Phi_{\alpha_0}(\zeta;0) \) and shall substitute it by \( \boldsymbol\Phi_{\alpha_0}^\mathsf{sing}(\zeta;0) \), which has slightly different behavior at infinity. This change necessitates modifications in matrices \( \boldsymbol M^\mathsf{reg}(z) \) and \( \boldsymbol E_0^\mathsf{reg}(z) \).

\subsection{On Matrices \( \boldsymbol E_0^\mathsf{reg}(z) \)}

To be able to introduce necessary modifications of \( \boldsymbol M^\mathsf{reg}(z) \) and \( \boldsymbol E_0^\mathsf{reg}(z) \), we need to discuss in more detail the behavior of \( \boldsymbol E_0^\mathsf{reg}(z) \) at \( a_0 \). To this end, let us write
\[
\boldsymbol E_0^\mathsf{reg}(z) = \left( \begin{matrix} E_{n,1}(z) & E_{n,2}(z) \smallskip \\ E_{n,3}(z) & E_{n,4}(z) \end{matrix} \right)n^{\sigma_3/6}.
\]
It follows readily from \eqref{E0} and \eqref{Mr} that
\begin{equation}
\label{fEn}
\left\{
\begin{array}{rcl}
\displaystyle E_{n,i}(z) &:=& \big(\ic\zeta^{1/2}(z)\big)^{1-i} G(z) \big( T_n(z) - (-1)^i (HT_n^*)(z)\big), \medskip \\
E_{n,i+2}(z) &:=& \big(\ic\zeta^{1/2}(z)\big)^{1-i} G(z) \left( \widetilde T_{n-1}(z) - (-1)^i \big(H\widetilde T_{n-1}^*\big)(z)\right),
\end{array}
\right. \quad z\in S_2,
\end{equation}
for \( i\in\{1,2\} \), where we set \( \widetilde T_m(\z) := (T_m/\Phi)(\z) \), and
\begin{equation}
\label{GandH}
G(z) := \frac{(\zeta^{1/4} S_\rho r)(z)}{\sqrt2}, \quad H(z) := \frac1{\ic\left((S_\rho r)^2w\right)(z)}, \quad z\in S_2.
\end{equation}
Notice that the functions \(  G(z) \) and \( H(z) \) remain bounded as \( S_2\ni z\to a_0 \) by \eqref{Srho-ai} and the very definitions of \( r(z) \) in \eqref{ri} and of \( \zeta(z) \) in \eqref{conf-map}.  We are interested in the quantity
\begin{equation}
\label{Dn}
D_n := E_{n,3}(a_0)E_{n,1}^\prime(a_0) - E_{n,1}(a_0)E_{n,3}^\prime(a_0).
\end{equation}
Since \( E_{n,i}(z) \) are analytic at \( a_0 \) and the functions \( T_n(\z) \) are analytic at \( \boldsymbol a_0 \), the lifts of \( G(z) \) and \( H(z) \) to \( \pi^{-1}(S_2)\cap \RS^{(0)} \) admit analytic continuations to some neighborhood of \( \boldsymbol a_0 \). Thus, the quantities \( G(a_0) \) and \( H(a_0) \) are well-defined. In fact, since \( E_{n,2}(a_0) \) and \( E_{n,4}(a_0) \) are finite while  \( T_n(a_0) = T_n^*(a_0) \) and \( \widetilde T_n(a_0) = \widetilde T_n^*(a_0) \), it must follow from \eqref{fEn} that
\begin{equation}
\label{H-exp}
H(z) = 1 + h_1(z-a_0)^{1/2} + h_2(z-a_0) + \cdots
\end{equation}
(from now on all the stated expansions are assumed to hold for \( z\in S_2\cap U_0 \)), where we use the same determination of \( (z-a_0)^{1/2} \) as in \eqref{Tn-exp}. It also clearly follows from \eqref{Tn-exp} and \eqref{prop5-5} that 
\[
\widetilde T_n(z) = t_0^{(n)} + t_1^{(n)}(z-a_0)^{1/2} + t_2^{(n)}(z-a_0) + \cdots,
\]
(that is, the first three coefficients in the Puiseux series of \( \widetilde T_n(z) \) and of \( T_n(z) \) are the same). Observe also that 
\[
T_n^*(z) = t_0^{(n)} - t_1^{(n)}(z-a_0)^{1/2} + t_2^{(n)}(z-a_0) + \cdots.
\]
Let us write \( G(z) = g_0 + g_1(z-a_0)^{1/2} + g_2(z-a_0) + \cdots \). Notice that \( g_0\neq0 \) by \eqref{Srho-ai}. Then a straightforward computation shows that
\begin{multline}
\label{Enkat0}
E_{n,2k+1}(z) = 2g_0 t_0^{(n-k)} + \left( \big( 2g_2 + g_1h_1 + g_0h_2 \big)t_0^{(n-k)} + \right. \\ + \left. g_0\left( 2t_2^{(n-k)} - h_1t_1^{(n-k)} \right) \right)(z-a_0) + \cdots
\end{multline}
for \( k\in\{0,1\} \) (in particular, it must hold that \( 2g_1+g_0h_1=0\)). Hence, we get from the definitions of the numbers \( D_n \) in \eqref{Dn} and the sequence \( \N_\varepsilon^\mathsf{sing} \) in \eqref{Nes} that
\begin{equation}
\label{Dn-lowb}
|D_n| = 2|g_0|^2|\det(\boldsymbol T_n)| \geq 2|g_0|^2\varepsilon>0, \quad n\in\N_\varepsilon^\mathsf{sing}.
\end{equation}

Let us also justify the remark made after Proposition~\ref{prop:lo}.  Analysis similar to the one above yields that
\[
E_{n,2k+2}(a_0) = -\ic g_0\zeta_1\left(2t_1^{(n-k)} - h_1t_0^{(n-k)} \right), \quad \zeta_1:= \left(\lim_{z\to a_0}\frac{z-a_0}{\zeta(z)}\right)^{1/2},
\]
for \( k\in\{0,1\} \). It further follows from \eqref{E0} and the display after \eqref{Mr} that
\[
\left(\gamma_n\gamma_{n-1}^*\right)^{-1} = \det\left(\boldsymbol E_0^\mathsf{reg}(z) \right) = \big(E_{n,1}E_{n,4}-E_{n,2}E_{n,3}\big)(z),
\]
which is true for any \( z\in U_0 \). Therefore, taking \( z=a_0 \) in the above equality gives that
\begin{equation}
\label{curious2}
\left(\gamma_n\gamma_{n-1}^*\right)^{-1} = -4\ic\zeta_1 g_0^2 \left(t_0^{(n)}t_1^{(n-1)} - t_0^{(n-1)}t_1^{(n)} \right).
\end{equation}
The desired claim now follows from \eqref{curious0}--\eqref{curious1}.

\subsection{Global and Local Parametrices}

Recall the definition of \( \boldsymbol M^\mathsf{reg}(z) \) in \eqref{Mr}. The new global parametrix is now defined as \( \boldsymbol N(z) = \boldsymbol {CM}^\mathsf{sing}(z)\Phi^{n\sigma_3}(z) \), where
\begin{equation}
\label{Ms}
 \boldsymbol M^\mathsf{sing}(z) := \left(\boldsymbol I + (z-a_0)^{-1}\boldsymbol L\right) \boldsymbol M^\mathsf{reg}(z)
\end{equation}
and the matrix \( \boldsymbol L \) is given by
\begin{equation}
\label{L}
\boldsymbol L := \frac1{D_n}\left( \begin{matrix} -(E_{n,1}E_{n,3})(a_0) & E_{n,1}^2(a_0) \medskip \\ -E_{n,3}^2(a_0) & (E_{n,1}E_{n,3})(a_0) \end{matrix}\right).
\end{equation}
Notice that \( \boldsymbol M^\mathsf{reg}(\infty) = \boldsymbol M^\mathsf{sing}(\infty) = \boldsymbol C^{-1} \) and, since \( \boldsymbol L \) has zero trace and determinant, \( \det(\boldsymbol M^\mathsf{reg}(z)) = \det(\boldsymbol M^\mathsf{sing}(z)) \). Observe also that the just defined matrix \( \boldsymbol N(z) \) satisfies \hyperref[rhn]{\rhn}(a,b) and \hyperref[rhn]{\rhn}(c) around \( a_1,a_2,a_3 \). However, its behavior at \( a_0 \) is obviously different. Most importantly for us, it is still true that
\begin{equation}
\label{Ms-bound}
\big|(\boldsymbol M^\mathsf{sing})^{\pm1}(z)\big| = \boldsymbol{\mathcal O}_{K,\varepsilon}(1)
\end{equation}
for \( z\in K \) as \( \N_\varepsilon^\mathsf{sing}\ni n\to\infty \), where \( K \) is any compact set that avoids \( a_0,a_1,a_2,a_3 \). Indeed, since \( \N_\varepsilon^\mathsf{sing} \subseteq \N_\varepsilon^\mathsf{reg} \), estimates \eqref{Mr-bound} still hold. It was computed in the previous subsection that
\begin{equation}
\label{En0-bound}
\big| E_{n,2k+1}(a_0) \big| = 2\big| g_0t_0^{(n-k)} \big| \leq 2|g_0|M_R
\end{equation}
for all \( n\in \N \), where the last conclusion is a consequence of \eqref{T-bounds1}. Thus, \eqref{Ms-bound} follows from \eqref{Dn-lowb}.

Local parametrix \( \boldsymbol P_0(z) \) is now constructed similarly to \eqref{P0} as
\begin{equation}
\label{P0-s}
\boldsymbol P_0(z) = \boldsymbol E_0^\mathsf{sing}(z)\boldsymbol\Phi_{\alpha_0}^\mathsf{sing}\big(n^{2/3}\zeta(z);0\big)\boldsymbol J(z)r^{-\sigma_3}(z),
\end{equation}
where \( \boldsymbol\Phi_{\alpha_0}^\mathsf{sing}(\zeta;x) \) is the solution of \hyperref[rhphis]{\rhphis} from Appendix~\ref{ap:34} and \( \boldsymbol E_0^\mathsf{sing}(z) \) is given by
\begin{equation}
\label{E0s}
\boldsymbol E_0^\mathsf{sing}(z) := \Big( \boldsymbol I + (z-a_0)^{-1}\boldsymbol L\Big) \boldsymbol E_0^\mathsf{reg}(z)\big(n^{2/3}\zeta(z)\big)^{-\sigma_3}.
\end{equation}
It follows from \eqref{E0}, \eqref{Phis-carried}, and \eqref{Ms} that \eqref{P0d} holds with \( \boldsymbol M^\mathsf{reg}(z) \) replaced by  \( \boldsymbol M^\mathsf{sing}(z) \). Hence, provided \( \boldsymbol E_0^\mathsf{sing}(z) \) is analytic in \( U_0 \), the just constructed matrix \( \boldsymbol P_0(z) \) solves  \hyperref[rhpo]{\rhpo} with \( \boldsymbol M^\mathsf{reg}(z) \) replaced by  \( \boldsymbol M^\mathsf{sing}(z) \) in \hyperref[rhpo]{\rhpo}(d) (again, entries of \( \boldsymbol L \) are bounded by \eqref{Dn-lowb} and \eqref{En0-bound} along \( \N_\varepsilon^\mathsf{sing}\)). To show analyticity of \( \boldsymbol E_0^\mathsf{sing}(z) \) in \( U_0 \), observe that
\[
\boldsymbol E_0^\mathsf{sing}(z) = \left( \begin{matrix} 1 - \frac{(E_{n,1}E_{n,3})(a_0)}{(z-a_0)D_n} & \frac{E_{n,1}^2(a_0)}{(z-a_0)D_n} \medskip \\ -\frac{E_{n,3}^2(a_0)}{(z-a_0)D_n} & 1 + \frac{(E_{n,1}E_{n,3})(a_0)}{(z-a_0)D_n} \end{matrix} \right) \left( \begin{matrix} \frac{E_{n,1}(z)}{z-a_0} & (z-a_0)* \medskip \\ \frac{E_{n,3}(z)}{z-a_0} & (z-a_0)* \end{matrix} \right) \left(\frac{n^{1/2}\zeta(z)}{z-a_0}\right)^{-\sigma_3}
\]
by \eqref{L}. Recall that \( \zeta(z) \) has a simple zero at \( a_0 \) and therefore the last matrix above is analytic at \( a_0 \). Thus, we only need to investigate analyticity of the first column of the product of the first two matrices above. The \((1,1)\)-entry of this product is equal to
\[
\frac1{z-a_0}\left(E_{n,1}(z) - E_{n,1}(a_0)\frac{E_{n,3}(a_0)E_{n,1}(z)-E_{n,1}(a_0)E_{n,3}(z)}{(z-a_0)D_n}\right),
\]
which is indeed analytic at \( a_0 \) by the very choice of \( D_n \) in \eqref{Dn}. The fact that \((2,1)\)-entry of the product is analytic at \( a_0 \) can be checked analogously. 

\subsection{A Discussion} In this subsection we briefly digress from the proof of Theorem~\ref{thm:asymp2} and provide a broader view of the construction of global and local paramatrices presented in this and the previous sections.

In this section we are forced to use matrix function \( \boldsymbol\Phi_{\alpha_0}^\mathsf{sing}(\zeta;0) \) instead of \( \boldsymbol\Phi_{\alpha_0}(\zeta;0) \) as the latter does not exist in the considered case due to polar singularities at \( x=0 \) of the functions \( U_k(x) \) from \eqref{UV}. However, one can clearly see from \eqref{Phi-special} that \( \boldsymbol\Phi_{\alpha_0}^\mathsf{sing}(\zeta;x) \) is actually well-defined around zero as long as \( U_1(0)\neq 0 \). A natural question arises whether one could use \( \boldsymbol\Phi_{\alpha_0}^\mathsf{sing}(\zeta;0) \) instead of \( \boldsymbol\Phi_{\alpha_0}(\zeta;0) \) in all cases leading to \( U_1(0)\neq 0 \). Let us explain that this is indeed possible.

If we want to use \( \boldsymbol\Phi_{\alpha_0}^\mathsf{sing}(\zeta;0) \) as model local parametrix, we do need to introduce modified global parametrix \( \boldsymbol M^\mathsf{sing}(z) \) via \eqref{Ms} with  \( \boldsymbol L \) given by \eqref{L}, but differently defined constants \( D_n \). Indeed, if \( U_1(0)\neq 0 \), then \eqref{Phis-carried} no longer holds and needs to be replaced by
\[
\left( \boldsymbol I - \frac1\zeta \left(\begin{matrix} 0 & 0 \\ 2^{-1/3}\ic U_1^{-1}(0) & 0 \end{matrix}\right) \right) \boldsymbol \Phi_\alpha^\mathsf{sing}(\zeta;0) = \frac{\zeta^{3\sigma_3/4}}{\sqrt2} \left(\begin{matrix} 1 & \ic \\ \ic & 1 \end{matrix}\right) \left(\boldsymbol I + \boldsymbol{\mathcal{O}}\Big(\zeta^{-1/2}\Big)\right) \exp\left\{ -\frac23 \zeta^{3/2}\sigma_3\right\}
\]
(in fact, the above formula extends \eqref{Phis-carried} since \( U_1^{-1}(0)=0 \) for weights in \( W^\mathsf{sing} \)). This change necessitates the following modification in the definition of \( \boldsymbol E_0^\mathsf{sing}(z) \):
\begin{align*}
\boldsymbol E_0^\mathsf{sing}(z) & := \left( \boldsymbol I + \frac{\boldsymbol L}{z-a_0}\right) \boldsymbol E_0^\mathsf{reg}(z)\big(n^{2/3}\zeta(z)\big)^{-\sigma_3} \left( \boldsymbol I - \frac{2^{-1/3}\ic U_1^{-1}(0)}{n^{2/3}\zeta(s)} \left(\begin{matrix} 0 & 0 \\ 1 & 0 \end{matrix}\right) \right) \\
&  = \left( \boldsymbol I + \frac{\boldsymbol L}{z-a_0}\right) \left[ \left( \begin{matrix} \frac{E_{n,1}(z)}{z-a_0} & (z-a_0)E_{n,2}(z) \medskip \\ \frac{E_{n,3}(z)}{z-a_0} & (z-a_0)E_{n,4}(z) \end{matrix} \right) + n^{1/3}\Lambda(z) \left( \begin{matrix} E_{n,2}(z) & 0 \medskip \\ E_{n,4}(z) & 0 \end{matrix} \right) \right] \left(\frac{n^{1/2}\zeta(z)}{z-a_0}\right)^{-\sigma_3},
\end{align*}
where \( \Lambda(z) := -2^{-1/3}\ic U_1^{-1}(0)\big(\zeta(z)/(z-a_0)\big)^2 \). Therefore, \( \boldsymbol E_0^\mathsf{sing}(z) \) is analytic at \( a_0 \) if
\[
\frac{\boldsymbol L}{(z-a_0)^2} \left( \begin{matrix} E_{n,1}(z) \medskip \\ E_{n,3}(z) \end{matrix} \right) + \frac{1}{z-a_0} \left( \begin{matrix} E_{n,1}(z) \medskip \\ E_{n,3}(z) \end{matrix} \right) + \frac{n^{1/3}\Lambda(z) \boldsymbol L}{z-a_0} \left( \begin{matrix} E_{n,2}(z) \medskip \\ E_{n,4}(z) \end{matrix} \right)
\]
is analytic at \( a_0 \). To achieve the latter, one needs to take \( \boldsymbol L \) as in \eqref{L} with
\[
D_n = n^{1/3}\Lambda(0)\big(E_{n,3}E_{n,2} - E_{n,1}E_{n,4}\big)(a_0) + \big(E_{n,3}E_{n,1}^\prime - E_{n,1}E_{n,3}^\prime\big)(a_0).
\]
Observe that when \( U_1^{-1}(0)=0 \), \( \Lambda(0)=0 \) and we recover the original definition of \( D_n \) in \eqref{Dn}. Hence, when \( U_1^{-1}(0)\neq0 \), it follows from \eqref{curious2} and the display before that \( D_n \) is non-zero for all \( n\in\N_\varepsilon^\mathsf{reg} \) large enough and the above construction goes through. It also holds in this case that \( \boldsymbol L=\boldsymbol{\mathcal O}\big(n^{-1/3}\big) \) and therefore matrices \( \boldsymbol M^\mathsf{reg}(z) \) and \( \boldsymbol M^\mathsf{sing}(z) \) are asymptotically the same.

\subsection{Proof of Theorem~\ref{thm:asymp2}} Consider now \hyperref[rhz]{\rhz} with \( \boldsymbol N_*(z) = \boldsymbol M^\mathsf{sing}(z)\Phi^{n\sigma_3}(z) \) and \( \boldsymbol P_0(z) \) given by \eqref{P0-s}. Again, let us estimate the size of the jump matrices in \hyperref[rhz]{\rhz}(b). It follows from \eqref{Ms-bound} that estimates \eqref{R-jump-1}--\eqref{R-jump-2} remain valid in the present case.  {Using \eqref{Phis-carried} instead of \hyperref[rhphi]{\rhphi}(d) gives the same estimate as in \eqref{Phi-a-est1}, one only needs to use boundedness of \( |Z|^{p+3/2}e^{-(4/3)|Z|^{3/2}} \). Since \hyperref[rhphis]{\rhphis}(c) is the same as \hyperref[rhphi]{\rhphi}(c), the estimate \eqref{Phi-a-est2} remain unchanged. It follows from the definition of \( \boldsymbol E_0^\mathsf{sing}(z) \) in \eqref{E0s}, its holomorphy in \( U_0 \), the maximum modulus principle, and \eqref{Ms-bound} that \( \boldsymbol E_0^\mathsf{sing}(z) = \boldsymbol{\mathcal O}_\varepsilon(1)n^{-\sigma_3/2}\) and \( \boldsymbol E_0^\mathsf{sing}(z)^{-1} = \boldsymbol{\mathcal O}_\varepsilon(1)n^{\sigma_3/2}\). Therefore, \eqref{R-jump-4} now is replaced by
\[
\boldsymbol I + \boldsymbol{\mathcal O}_\varepsilon\big(n\max |s-a_0|^p|Z|^{-p}\big) = \boldsymbol I + \boldsymbol{\mathcal O}_\varepsilon\big(n^{-(2p-3)/3}\big).
\] 
Recall that \( p \) is the largest integer such that \( |1-\kappa_{i,j}(s)| = \mathcal O( |s-a_0|^p )\) as \( s\to a_0 \). It can be readily checked that condition \eqref{bad-cond} implies that \( p\geq2\) and therefore \( (2p-3)/3 \geq 1/3 \).} Thus, the jump matrices in  \hyperref[rhz]{\rhz}(b) still satisfy the estimate \( \boldsymbol I + \boldsymbol{\mathcal O}_\varepsilon(n^{-1/3}) \). Hence, we again can conclude that \( \boldsymbol Z(z) \) exists for all \( n\in\N_\varepsilon^\mathsf{sing} \) large enough and satisfies \eqref{z-rate}.

As in the proof of Theorem~\ref{thm:asymp1}, given a closed set \( K\subset\overline\C\setminus\Delta  \), the contour \( \Sigma \) can always be adjusted so that \( K \) lies in the exterior domain of \( \Sigma \). Then we get from \eqref{X} and \eqref{Z} that
\[
\boldsymbol Y(z) = \boldsymbol X(z) = \boldsymbol C\boldsymbol M^\mathsf{sing}(z)\Phi^{n\sigma_3}(z), \quad z\in K.
\]
Therefore, we get from \eqref{Y} that
\begin{multline*}
\gamma_n^{-1}Q_n(z) = \big(1+\upsilon_{n1}^\mathsf{sing}(z)\big) \left( \left(1 + \frac{[\boldsymbol L]_{11}}{z-a_0} \right)\Psi_n(z) + \frac{[\boldsymbol L]_{12}}{z-a_0}\Psi_{n-1}(z) \right)  \\ + \upsilon_{n2}^\mathsf{sing}(z) \left( \left(1+\frac{[\boldsymbol L]_{22}}{z-a_0}\right)\Psi_{n-1}(z) + \frac{[\boldsymbol L]_{21}}{z-a_0}\Psi_n(z) \right), \quad z\in K,
\end{multline*}
where \( 1+\upsilon_{n1}^\mathsf{sing}(z),  \upsilon_{n2}^\mathsf{sing}(z) \) are the first row entries of \( \boldsymbol Z(z) \). For \( i\in\{1,2\} \), let
\begin{equation}
\label{new-upsilons}
\upsilon_{ni}(z) := \upsilon_{ni}^\mathsf{sing}(z) \left(1+\frac{[\boldsymbol L]_{ii}}{z-a_0}\right) + \upsilon_{n3-i}^\mathsf{sing}(z) \frac{[\boldsymbol L]_{3-ii}}{z-a_0}.
\end{equation}
It follows from \eqref{L}, \eqref{En0-bound}, and \eqref{Dn-lowb} that \( \upsilon_{ni}(z) \) satisfy bounds as in \eqref{upsilons}. Moreover, we have that
\[
\gamma_n^{-1}Q_n(z) = \left(1 + \frac{[\boldsymbol L]_{11}}{z-a_0} + \upsilon_{n1}(z) \right) \Psi_n(z) + \left( \frac{[\boldsymbol L]_{12}}{z-a_0} + \upsilon_{n2}(z) \right) \Psi_{n-1}(z), \quad z\in K.
\]
We further get from \eqref{L}, \eqref{Enkat0}, and \eqref{Dn-lowb} that
\[
[\boldsymbol L]_{11} = -\frac{(E_{n,1}E_{n,3})(a_0)}{D_n} = -\frac{2t_0^{(n-1)}t_0^{(n)}}{\det(\boldsymbol T_n)} \qandq [\boldsymbol L]_{12} = \frac{E_{n,1}^2(a_0)}{D_n} = \frac{2(t_0^{(n)})^2}{\det(\boldsymbol T_n)}.
\]
This finishes the proof of the asymptotic formula for \( Q_n(z) \) while the formula for \( (wR_n)(z) \) can be shown absolutely analogously. Similarly, we get from \eqref{Y}, \eqref{X}, and \eqref{Z} with \( \boldsymbol M^\mathsf{reg}(z) \) replaced by \( \boldsymbol M^\mathsf{sing}(z) \) that
\begin{multline*}
\gamma_n^{-1}Q_n(z) = \left( \left(1 + \frac{[\boldsymbol L]_{11}}{z-a_0} + \upsilon_{n1}(z) \right) \Psi_n(z) + \upsilon_{n2}(z) \Psi_{n-1}(z) \right) \\ \pm (\rho_iw)^{-1}(z)\left( \left(1 + \frac{[\boldsymbol L]_{11}}{z-a_0} + \upsilon_{n1}(z) \right) \Psi_n^*(z) + \upsilon_{n2}(z) \Psi_{n-1}^*(z) \right)
\end{multline*}
for \( z\in\Omega_i^\pm\setminus \overline U \) with \( \upsilon_{ni}(z) \) as in \eqref{new-upsilons}. After that the proof of the analog of the first relation in \eqref{Asymp2} proceeds exactly as after \eqref{Qn-Delta}. The proof of the second relation in \eqref{Asymp2} can be obtained analogously.

\appendix

\section{Riemann-Hilbert Problem for Painlev\'e II}
\label{sec:painleve}

\subsection{Lax Pair and the Corresponding Riemann-Hilbert Problem}

The material of this section originates in \cite{FlNew80}. However, we essentially follow the presentation in \cite[Section~5.0]{FokasItsKapaevNovokshenov}, see also \cite[Section~1.0]{FokasItsKapaevNovokshenov} for the relevant facts of the general theory of differential equations. Let \( q(s) \) be a solution of Painlev\'e II equation 
\begin{equation}
\label{p2}
q^{\prime\prime}(s) = sq(s) + 2q^3(s) - \nu,  \quad \nu>-1/2.
\end{equation}
It is known that \( q(s) \) is a meromorphic function in the entire complex plane. For each \( s \), which is not a pole of \( q \), consider the following system of differential equations:
\begin{equation}
\label{pIIa}
\partial_\tau \boldsymbol\Psi(\tau;s) = \boldsymbol A^{FN}(\tau;s)\boldsymbol\Psi(\tau;s),
\end{equation}
where
\begin{equation}
\label{pIIb}
\boldsymbol A^{FN}(\tau;s) := \left(\begin{matrix} -\ic(4\tau^2+s+2q^2(s)) & 4\tau q(s) + 2\ic q^\prime(s)+\nu/\tau \smallskip \\ 4\tau q(s) - 2\ic q^\prime(s)+\nu/\tau & \ic(4\tau^2+s+2q^2(s)) \end{matrix}\right)
\end{equation}
and
\begin{equation}
\label{pIIc}
\partial_s \boldsymbol\Psi(\tau;s) = \boldsymbol U^{FN}(\tau;s)\boldsymbol\Psi(\tau;s), \quad \boldsymbol U^{FN}(\tau;s) := \left(\begin{matrix} -\ic\tau & q(s) \smallskip \\ q(s) & \ic\tau \end{matrix}\right).
\end{equation}
In general, such a system would be overdetermined, but the compatibility condition
\begin{equation}
\label{compat-cond}
\partial_s \boldsymbol A^{FN}(\tau;s) - \partial _\tau \boldsymbol U^{FN}(\tau;s) + \big[\boldsymbol A^{FN}(\tau;s),\boldsymbol U^{FN}(\tau;s)\big] = 0  
\end{equation}
of \eqref{pIIa}--\eqref{pIIc} exactly reduces to the fact that \( q(s) \) solves \eqref{p2}. The general theory of differential equations implies that equation \eqref{pIIa}--\eqref{pIIb} has a set of canonical solutions \( \boldsymbol\Psi_k(\tau;s) \), \( k\in\{1,\ldots,6\} \), that are uniquely determined by the conditions \( \det(\boldsymbol\Psi_k(\tau;s))\equiv1 \) and
\begin{equation}
\label{Psi-sector}
\boldsymbol\Psi_k(\tau;s) = \left(\boldsymbol I + \boldsymbol{\mathcal O}\left(\tau^{-1}\right) \right) e^{-\ic(\frac43\tau^3+s\tau)\sigma_3} \qasq \tau\to\infty,~~\arg(\tau)\in\left(\frac{(k-2)\pi}3,\frac{k\pi}3\right).
\end{equation}
Given any solution \( \boldsymbol\Psi(\tau;s) \) of \eqref{pIIa}--\eqref{pIIb} one can readily check by using \eqref{compat-cond} that
\[
\partial_s \boldsymbol\Psi(\tau;s) - \boldsymbol U^{FN}(\tau;s)\boldsymbol\Psi(\tau;s)
\]
must also satisfy \eqref{pIIa}--\eqref{pIIb}. Hence, this difference is equal to \( \boldsymbol\Psi(\tau;s)\boldsymbol S(s) \), where \( \boldsymbol S(s) \) is a matrix that does not depend on \( \tau \) and can be expressed as
\begin{equation}
\label{right-multiplier}
\boldsymbol S(s) = \boldsymbol\Psi^{-1}(\tau;s)\partial_s \boldsymbol\Psi(\tau;s) - \boldsymbol\Psi^{-1}(\tau;s)\boldsymbol U^{FN}(\tau;s)\boldsymbol\Psi(\tau;s).
\end{equation}
A further straightforward computation leads to the fact that the whole system \eqref{pIIa}--\eqref{pIIc} is solved by
\begin{equation}
\label{general-solution}
\boldsymbol\Psi(\tau;s)\boldsymbol Q(s), \quad \boldsymbol \partial_s \boldsymbol Q(s) = - (\boldsymbol{SQ})(s).
\end{equation}
By explicitly computing the term next to \( \tau^{-1} \) in \eqref{Psi-sector} (see \eqref{pII-exp} and \eqref{AB12} further below and note that the asymptotic expansion of \( \boldsymbol\Psi_k(\tau;s) \) holds in the whole Stokes sector \( ((k-2)\pi/3,k\pi/3) \)), plugging this asymptotic formula into \eqref{right-multiplier}, and using the fact that \eqref{right-multiplier} must be independent of \( \tau \), one can deduce that \( \boldsymbol S(s) \) corresponding to any \( \boldsymbol\Psi_k(\tau;s) \) is identically zero, i.e., we can take \( \boldsymbol Q(s)=\boldsymbol I \) in \eqref{general-solution}, and therefore \( \boldsymbol\Psi_k(\tau;s) \) satisfies \eqref{pIIc} as well (see the paragraph containing \cite[Equation~(5.0.7)]{FokasItsKapaevNovokshenov}). 

 Define the rays \( \Sigma_k^{FN} \) and sectors \( S_k^{FN} \) by
 \[
\Sigma_k^{FN} := \bigg\{\tau:\arg(\tau)=(2k-1)\pi/6\bigg\}  \qandq S_k^{FN}:= \left\{ z:~\arg(z)\in\left(\frac{(2k-3)\pi}6,\frac{(2k-1)\pi}6\right)\right\},
\]
where the rays are oriented away from the origin. Further, put
\begin{equation}
\label{FN}
\boldsymbol\Psi_\nu^{FN}(\tau;s) := \boldsymbol\Psi_k(\tau;s), \quad \tau\in S_k^{FN}.
\end{equation}
Since the matrices \( \boldsymbol\Psi_k(\tau;s) \) solve the same differential equation, they are related to each other through the right multiplication by a constant matrix. By studying their behavior at infinity, one can conclude that
\begin{equation}
\label{FN-jump}
\boldsymbol\Psi_{\nu+}^{FN}(\tau;s) = \boldsymbol\Psi_{\nu-}^{FN}(\tau;s)
\left\{
\begin{array}{rl}
\begin{pmatrix} 1 & 0 \\ s_{2j-1} & 1 \end{pmatrix}, & \tau\in\Sigma_{2j-1}^{FN}, \bigskip \\
\begin{pmatrix} 1 & s_{2j} \\ 0 & 1 \end{pmatrix}, & \tau\in\Sigma_{2j}^{FN},
\end{array}
\right.
\end{equation}
for some constants \( s_1,\ldots,s_6 \) (Stokes parameters). Since \( \sigma_1 \boldsymbol A^{FN}(-\tau;s)\sigma_1 = -\boldsymbol A^{FN}(\tau;s)\) and due to the uniqueness of the solution of \eqref{pIIa}--\eqref{pIIb} satisfying \eqref{Psi-sector} in a given sector, it holds that
\begin{equation}
\label{sigma1-sym}
\boldsymbol\Psi_\nu^{FN}(\tau;s) = \sigma_1\boldsymbol\Psi_\nu^{FN}(-\tau;s)\sigma_1 , \quad \sigma_1 := \left( \begin{matrix} 0 & 1 \\ 1 & 0 \end{matrix}\right), \quad \sigma_2 := \left( \begin{matrix} 0 & -\ic \\ \ic & 0 \end{matrix}\right).
\end{equation}
This implies that there are only three independent Stokes parameters and \( s_1=s_4 \), \( s_2=s_5 \), \(s_3=s_6 \).

On the other hand, since the residue of \( \boldsymbol A^{FN}(\tau;s) \) at \( \tau=0 \) is equal to
\[
\begin{pmatrix} 0 & \nu \\ \nu & 0 \end{pmatrix} = \frac1{\sqrt 2} \begin{pmatrix} 1 & -1 \\  1 & 1 \end{pmatrix} \begin{pmatrix} \nu & 0 \\ 0 & -\nu \end{pmatrix} \frac1{\sqrt 2} \begin{pmatrix} 1 & 1 \\  -1 & 1 \end{pmatrix},
\]
where all matrices have unit determinants, it holds that the canonical solution of \eqref{pIIa}--\eqref{pIIb} at the origin that has unit determinant\footnote{Of course, \( -\boldsymbol\Psi^{(0)}(\tau;s) \) also solves \eqref{pIIa}--\eqref{pIIb} and has unit determinant.} is of the form
\[
\boldsymbol\Psi^{(0)}(\tau;s) = \frac1{\sqrt 2} \begin{pmatrix} 1 & -1 \\  1 & 1 \end{pmatrix}\left( \boldsymbol I + \sum_{k=1}^\infty \boldsymbol\Psi_{0,k}(s) \tau^k \right) \tau^{\nu\sigma_3} \begin{pmatrix} 1 & \chi_\nu( \Pi(s) + \kappa(s)\log\tau) \\ 0 & 1 \end{pmatrix},
\]
where \( \chi_\nu = 0 \) in the non-resonant cases \( \nu+1/2\notin\N \) and \( \chi_\nu = 1 \) in the resonant cases \( \nu+1/2\in\N \). The matrices \(  \boldsymbol\Psi_{0,k}(s)  \) and the number \( \kappa(s) \) are uniquely determined by \eqref{pIIa}--\eqref{pIIb} under the condition that we set \( (1,2) \)-entry of \( \boldsymbol\Psi_{0,2\nu}(s) \) to be zero in the resonant cases, while \( \Pi(s) \) is a free parameter for each fixed \( s \) (equivalently, we could set \( \Pi(s)\equiv0 \) and say that \( (1,2) \)-entry of \( \boldsymbol\Psi_{0,2\nu}(s) \) is a free parameter, whose choice necessarily affects the values of the second column of every \( \boldsymbol\Psi_{0,k}(s) \), \( k> 2\nu\)). Let us put \( \eta(\tau;s) := \chi_\nu(\Pi(s) + \kappa(s)\log\tau) \). Then we get that
\[
\big(\boldsymbol\Psi^{(0)}\big)^{-1}\partial_s \boldsymbol\Psi^{(0)} = \begin{pmatrix} 1 & -\eta \\ 0 & 1 \end{pmatrix} \tau^{-\nu\sigma_3} \boldsymbol{\mathcal O}(\tau) \tau^{\nu\sigma_3} \begin{pmatrix} 1 & \eta \\ 0 & 1 \end{pmatrix} + \begin{pmatrix} 0 & \partial_s\eta \\ 0 & 0 \end{pmatrix}.
\]
Similarly, we obtain from \eqref{pIIc} that
\[
\big(\boldsymbol\Psi^{(0)}\big)^{-1}\boldsymbol U^{FN}\boldsymbol\Psi^{(0)} = q\sigma_3 + \begin{pmatrix} 0 & 2q\eta \\ 0 & 0 \end{pmatrix} + \begin{pmatrix} 1 & -\eta \\ 0 & 1 \end{pmatrix} \tau^{-\nu\sigma_3} \boldsymbol{\mathcal O}(\tau) \tau^{\nu\sigma_3} \begin{pmatrix} 1 & \eta \\ 0 & 1 \end{pmatrix}.
\]
Plugging the above expressions into \eqref{right-multiplier} gives
\[
\boldsymbol S(s) = -q\sigma_3 + \begin{pmatrix} 0 & \partial_s\eta-2q\eta \\ 0 & 0 \end{pmatrix} + \tau \begin{pmatrix} O_1-\eta \tau^{2\nu}O_3 & -\eta^2 \tau^{2\nu} O_3 + \eta O_1 -\eta O_4 + \tau^{-2\nu}O_2 \\ \tau^{2\nu} O_3 & O_4 + \eta\tau^{2\nu} O_3 \end{pmatrix},
\]
where \( O_l(\tau;s) \) are some analytic functions of \( \tau \) coming from the terms \( \tau^{-\nu\sigma_3} \boldsymbol{\mathcal O}(\tau) \tau^{\nu\sigma_3} \) in the preceeding two relations. Since the above expression must be independent of \( \tau \) we see immediately that \( O_1 = O_3 = O_4 \equiv 0 \) and \( O_2(\tau;s) = O(s)\tau^{2\nu-1} \) for some suitable function \( O(s) \), in which case
\[
\boldsymbol S(s) = -q(s)\sigma_3 + \begin{pmatrix} 0 & (\kappa^\prime(s)-2(q\kappa)(s))\log\tau + \Pi^\prime(s) - 2(q\Pi)(s) + O(s) \\ 0 & 0 \end{pmatrix}.
\] 
Let \( Q(s) \) be the exponential of some fixed antiderivative of \( q(s) \). The independence from \( \tau \) yields that \( \kappa(s) = \kappa_\nu Q^2(s) \) for some constant \( \kappa_\nu \). Further, let \( Q_*(s) \) be a function (particular solution) such that \( Q_*^\prime(s) = (Q_*q)(s) - (OQ^{-1})(s) \). Then the differential equation in \eqref{general-solution} for \( \boldsymbol S(s) \) as above is solved by
\[
\boldsymbol Q(s) = \begin{pmatrix} Q(s) & hQ(s) - \Pi(s)Q^{-1}(s) + Q_*(s) \\ 0 & Q^{-1}(s) \end{pmatrix},
\]
where \( h \) is an arbitrary constant (we look for solutions of determinant one and therefore use exponential of the same antiderivative of \( q(s) \) in the diagonal entries). Therefore, the solution of the whole system \eqref{pIIa}--\eqref{pIIc} around the origin assumes the form
\begin{equation}
\label{Psi-hat0}
\boldsymbol\Psi_0(\tau;s) = \frac{\boldsymbol I -\ic\sigma_2}{\sqrt 2} \big( \boldsymbol I + \boldsymbol{\mathcal O}(\tau) \big) \tau^{\nu\sigma_3} \begin{pmatrix} 1 & \chi_\nu( (\kappa_\nu\log\tau + h)Q^2(s) + (Q_*Q)(s)) \\ 0 & 1 \end{pmatrix} Q^{\sigma_3}(s).
\end{equation}

Since the matrices \( \boldsymbol\Psi_0(\tau;s) \) and \( \boldsymbol\Psi_1(\tau;s) \) solve the same system of differential equations, namely, \eqref{pIIa}--\eqref{pIIc}, they are connected via the right multiplication by a constant matrix. Thus, we can write
\begin{equation}
\label{Psi1-hat0}
\boldsymbol\Psi_1(\tau;s) = \boldsymbol\Psi_0(\tau;s) \boldsymbol E
\end{equation}
in some disk around the origin. Notice that \( \boldsymbol\Psi_0(\tau;s) \) must also obey the symmetry in \eqref{sigma1-sym}. Hence, \( \sigma_1\boldsymbol\Psi_0(-\tau;s)\sigma_1 \) is connected to \( \boldsymbol\Psi_0(\tau;s) \) through the right multiplication by a constant matrix as well. In fact, it is not hard to compute that
\[
\sigma_1\boldsymbol\Psi_0(-\tau;s)\sigma_1 = \boldsymbol\Psi_0(\tau;s) e^{\pi\ic\nu\sigma_3}  \begin{pmatrix} \pi\ic\kappa_\nu & 1 \\ -1 & 0 \end{pmatrix}. 
\]
This relation together with \eqref{FN-jump} and \eqref{Psi1-hat0} yields that
\[
\boldsymbol E = e^{\pi\ic\nu\sigma_3}  \begin{pmatrix} \pi\ic\kappa_\nu & 1 \\ -1 & 0 \end{pmatrix} \sigma_1 \boldsymbol E  \begin{pmatrix} 1 & 0 \\ s_1 & 1 \end{pmatrix} \begin{pmatrix} 1 & s_2 \\ 0 & 1 \end{pmatrix} \begin{pmatrix} 1 & 0 \\ s_3 & 1 \end{pmatrix} \sigma_1,
\]
see \cite[Equation~(5.0.14)]{FokasItsKapaevNovokshenov}. Solving the above system of relations gives that
\begin{equation}
\label{connection}
\boldsymbol E = \begin{pmatrix} 1 & \chi_\nu p \\ 0 & 1 \end{pmatrix} d^{\sigma_3} \boldsymbol A_1^{FN}
\end{equation}
where \( d\neq0,p \) are some constants and the matrix \( \boldsymbol A_1^{FN} \) is explicitly known and involves expressions depending on \( s_1,s_2,s_3 \) and \( \nu \), see \cite[Equations (5.0.17)--(5.0.18)]{FokasItsKapaevNovokshenov} excluding the exceptional resonant cases \( \nu+1/2\in\N \) and \( s_1=s_2=s_3 = \ic(-1)^{\nu+1/2} \) in which it is a lower-triangular matrix with one's on the main diagonal (notice that we renamed \( \sqrt{-J_+} \) from \cite[Equations (5.0.18)]{FokasItsKapaevNovokshenov} as \( d \) here). Moreover, this computation also yields that
\begin{equation}
\label{pII-stokes}
s_1 + s_2 + s_3 + s_1s_2s_3=-2\ic\sin(\nu\pi),
\end{equation}
and that \( \kappa_\nu = (-1)^{\nu+1/2}/\pi \) in the non-exceptional resonant cases and \( \kappa_\nu=0 \) otherwise. Combining \eqref{Psi-hat0} with \eqref{Psi1-hat0} then gives that  in the neighborhood of \( \tau =0 \) the function \( \boldsymbol \Psi_{\nu}^{FN}(\tau;s) \) admits a representation
\begin{equation}
\label{FN-origin}
\boldsymbol\Psi_\nu^{FN}(\tau;s) = \boldsymbol H^{FN}(\tau;s) \tau^{\nu\sigma_3} \left(\begin{matrix} 1 & \kappa_\nu\log\tau \\ 0 & 1 \end{matrix}\right) \boldsymbol A_k^{FN}, \quad \tau\in S_k^{FN},
\end{equation}
where $\boldsymbol H^{FN}(\tau;s)$ is a holomorphic matrix function of \( \tau \) around the origin while the matrices \( \boldsymbol A_k^{FN} \) are connected to each other via the jump relations \eqref{FN-jump} with \( \boldsymbol A_1^{FN} \) as in \eqref{connection}.

Altogether, we see that \( \boldsymbol\Psi_\nu^{FN}(\tau;s) \) solves the following Riemann-Hilbert problem (\rhpsi):

\begin{itemize}
\label{rhpsi}
\item[(a)] $\boldsymbol\Psi_\nu(\tau;s)$ is holomorphic in $\C\setminus\big(\{0\}\cup\Sigma^{FN}\big)$, \( \Sigma^{FN} := \cup_{k=1}^6\Sigma_k^{FN} \);
\item[(b)] $\boldsymbol\Psi_\nu(\tau;s)$ has continuous traces on $\Sigma^{FN}$ that satisfy \eqref{FN-jump};
\item[(c)] $\boldsymbol\Psi_\nu(\tau;s)$ satisfies \eqref{FN-origin} as $\tau\to0$;
\item[(d)] $\boldsymbol\Psi_\nu(\tau;s)$ satisfies \( \boldsymbol\Psi_\nu(\tau;s) =  \left(\boldsymbol I+\boldsymbol{\mathcal{O}}\left(\tau^{-1}\right)\right) e^{-\ic(\frac43\tau^3+s\tau)\sigma_3} \) uniformly in $\{|\tau|\geq1\}\setminus\Sigma^{FN}$.
\end{itemize}

Conversely, it follows from \cite[Theorem~5.1]{FokasItsKapaevNovokshenov} that given \( \nu>-1/2 \), parameters \( s_1,s_2,s_3 \) satisfying \eqref{pII-stokes}, and, in the exceptional resonant cases, the value of the free parameter of the matrix \( \boldsymbol A_1^{FN} \) from \eqref{connection}, a solution of \hyperref[rhpsi]{\rhpsi} uniquely exists as a meromorphic function of \( s \) and
\[
q(s) = 2\ic \lim_{\R\ni\tau\to\infty} \tau \left[\boldsymbol\Psi_\nu(\tau;s)e^{\ic(\frac43\tau^3+s\tau)\sigma_3}\right]_{12}
\]
solves \eqref{p2}, where \( [\cdot]_{12}\) is the \( (1,2) \)-entry of the corresponding matrix. In particular, Stokes parameters \( s_1,s_2,s_3 \) uniquely parametrize solutions of \eqref{p2} excluding the exceptional resonant cases each of which corresponds to a family of solutions parametrized by the free parameter of \( \boldsymbol A_1^{FN} \).

\subsection{Asymptotic Expansion}

Later we shall need the first four terms of the full asymptotic expansion of \( \boldsymbol\Psi_\nu^{FN}(\tau;s) \) at infinity (the first is the identity matrix, so, we need to find the other three). Following classical ideas, see \cite[Proposition~1.1]{FokasItsKapaevNovokshenov}, one can write the expansion of \( \boldsymbol\Psi_\nu^{FN}(\tau;s) \) as a series in \( 1/\tau \) with off-diagonal coefficients times the exponential of a series with diagonal coefficients. Symmetry \eqref{sigma1-sym} yields that
\[
\boldsymbol\Psi_\nu^{FN}(\tau;s) =  \left( \boldsymbol I + \sum_{n=1}^\infty\tau^{-n}\psi_n(s)\boldsymbol \Sigma_n \right) \exp\left\{ \sum_{n=-3,n\neq0}^\infty \frac{\tau^{-n}}{n}\lambda_n(s)\boldsymbol\Lambda_n \right\},
\]
where \( \boldsymbol \Sigma_{2m} = \sigma_1 \) and \( \boldsymbol \Sigma_{2m+1}=\ic\sigma_2 \), while \( \boldsymbol\Lambda_{2m}=\boldsymbol I \) and \( \boldsymbol\Lambda_{2m+1} = \sigma_3 \). Plugging the above expansion into \eqref{pIIa}--\eqref{pIIb} gives
\begin{multline*}
- \sum_{n=2}^\infty (n-1)\tau^{-n}\psi_{n-1}\boldsymbol \Sigma_{n-1} - \left( \boldsymbol I + \sum_{n=1}^\infty\tau^{-n}\psi_n\boldsymbol \Sigma_n \right) \sum_{n=-2,n\neq1}^\infty \tau^{-n}\lambda_{n-1}\boldsymbol\Lambda_{n-1} = \\ \left( -4\ic\tau^2\sigma_3 + 4q\tau\sigma_1 - \ic\big(s+2q^2\big)\sigma_3 -2q^\prime\sigma_2 + \nu\tau^{-1}\sigma_1  \right) \left( \boldsymbol I + \sum_{n=1}^\infty\tau^{-n}\psi_n\boldsymbol \Sigma_n \right).
\end{multline*}
From the coefficients next to \( \tau^2 \) and \( \tau \) we get that
\[
\left\{
\begin{array}{rcl}
 \lambda_{-3}(s) & = & 4\ic, \smallskip \\
\psi_1(s) & = & -\ic q(s)/2, \smallskip \\
\lambda_{-2}(s) & \equiv & 0.
\end{array}
\right.
\]
By comparing the constant terms on both sides as well as the terms next to \( \tau^{-1} \), we get
\[
\left\{
\begin{array}{rcl}
\lambda_{-1}(s) & = & \ic s, \smallskip \\
\psi_2(s) & = & q^\prime(s)/4, \smallskip \\
\psi_3(s) & = & \ic (q^3(s)+ sq(s) - \nu)/8.
\end{array}
\right.
\]
By equating the terms next to \( \tau^{-2} \), we obtain
\[
\left\{
\begin{array}{rcl}
\lambda_1(s) & = & \ic \left(q^4(s)+sq^2(s)-2\nu q(s)-(q^\prime(s))^2\right)/2, \smallskip \\
\psi_4(s) & = & -\left(sq^\prime(s)+q^2(s)q^\prime(s)+q(s)\right)/16.
\end{array}
\right.
\]
Lastly, the coefficients next \( \tau^{-3} \) and \( \tau^{-4} \) give us
\[
\left\{
\begin{array}{rcl}
\lambda_2(s) & = & q^2(s)/4, \smallskip \\
\psi_5(s) & = & -\ic\left( 2q^\prime(s) + q(s)H(s) + \left(q^2(s)+s\right)\left(q^3(s)+sq(s)-\nu \right) \right)/32, \smallskip \\
\lambda_3(s) & = & -\ic\left( q^\prime(s)q(s) - sH(s) + \nu^2 \right)/8,
\end{array}
\right.
\]
where we denote by \( H(s) := \big( q^\prime(s) \big)^2 - sq^2(s) - q^4(s) + 2\nu q(s) \) the Hamiltonian of Painlev\'e II equation. Now, we would like to rewrite the expansion of  \( \boldsymbol\Psi_\nu^{FN}(\tau;s) \) at infinity as
\begin{equation}
\label{pII-exp}
\boldsymbol\Psi_\nu^{FN}(\tau;s) =  \left(\boldsymbol I + \sum_{n=1}^\infty \frac1{(-\ic \tau)^n} \left( \begin{matrix} A_n(s) & B_n(s) \smallskip \\ (-1)^nB_n(s) & (-1)^nA_n(s) \end{matrix} \right) \right) \exp\left\{-\ic\left(\frac43\tau^3+s\tau\right)\sigma_3\right\}
\end{equation}
for some functions \( A_n(s),B_n(s) \), \( n\in\N \), where the form of the coefficients follows from \eqref{sigma1-sym}. Since
\begin{multline*}
\exp\left\{ \sum_{n=-3,n\neq0}^\infty \frac{\tau^{-n}}{n}\lambda_n(s)\boldsymbol\Lambda_n \right\} = \left( \boldsymbol I + \frac{\lambda_1(s)}\tau \sigma_3 + \frac{\lambda_2(s)+\lambda_1^2(s)}{2\tau^2}\boldsymbol I \right. + \\ \left. \frac{2\lambda_3(s)+ 3(\lambda_1\lambda_2)(s)+\lambda_1^3(s)}{6\tau^3} \sigma_3 + \boldsymbol{\mathcal O} \left(\frac1{\tau^4}\right) \right) \times\exp\left\{-\ic\left(\frac43\tau^3+s\tau\right)\sigma_3\right\},
\end{multline*}
we get that
\[
\left\{
\begin{array}{ll}
\ic A_1(s)\sigma_3 - B_1(s)\sigma_2 & = \lambda_1(s)\sigma_3 + \ic\psi_1(s)\sigma_2, \smallskip \\
- A_2(s)\boldsymbol I - B_2(s)\sigma_1 & = \left(\lambda_2+\lambda_1^2\right)(s) \boldsymbol I/2 + (\psi_2-\psi_1\lambda_1)(s)\sigma_1, \smallskip \\
-\ic A_3(s)\sigma_3 + B_3(s)\sigma_2 & = \left(2\lambda_3+3\lambda_2\lambda_1+\lambda_1^3\right)(s)\sigma_3/6 + \smallskip \\
& \hspace{.5in} \ic\left( 2\psi_3-2\psi_2\lambda_1+\psi_1\left(\lambda_2+\lambda_1^2\right) \right)(s)\sigma_2/2.
\end{array}
\right.
\]
Then it holds that
\begin{equation}
\label{AB12}
\left\{
\begin{array}{rcl}
-2B_1(s) & = & q(s), \smallskip \\ 
-2A_1(s) & =  & H(s),
\end{array}
\right.
\quad
\left\{
\begin{array}{rcl}
-4B_2(s) & = & q^\prime(s)+q(s)H(s), \smallskip \\ 
-8A_2(s) & = & q^2(s) - H^2(s),
\end{array}
\right.
\end{equation}
and, by recalling \eqref{p2}, that
\[
\left\{
\begin{array}{rcl}
-16B_3(s) & = & 2q^{\prime\prime}(s) - 3q^3(s) + 2q^\prime(s)H(s) + q(s)H^2(s), \smallskip  \\ 
48A_3(s) & = & 2q^\prime(s)q(s) + \big(3 q^2(s)-2s\big)H(s) - H^3(s) + 2\nu^2.
\end{array}
\right.
\]
One also can verify either directly or by using \eqref{pIIc} that \( A_n^\prime(s) = (-1)^nq(s)B_n(s) \).

Our primary interest is the behavior of \( A_n(s), B_n(s) \) around a pole of \( q(s) \). It can be readily checked that if \( q(s) \) has a pole at  \( s_0 \), then
\begin{equation}
\label{pIIpole}
q(s) = \displaystyle \frac{q_{-1}}{s-s_0}+ \mathcal O\big(s-s_0), \quad q_{-1}\in\{-1,1\}, \qandq H(s) = \displaystyle \frac1{s-s_0} + \mathcal O(1)
\end{equation}
as \( s\to s_0 \). We shall need more precise behavior of \( q(s) \) and \( H(s) \) around \( s_0=0 \) when \( q_{-1}=-1 \).  It is known, see \cite[Equation (17.1)]{GromakLaineShimomura}, or can be obtained directly from \eqref{p2}, that
\begin{equation}
\label{pIIpole0}
\left\{
\begin{array}{rcl}
q(s) & = & \displaystyle -\frac1s + \frac{\nu+1}4 s^2 + \frac{\mathfrak q}5s^3 + \mathcal O\left(s^5\right), \medskip  \\ 
H(s) & = & \displaystyle \frac1s + 2\mathfrak q + \frac{\nu+1}4 s^2 + \frac{2\mathfrak q}{15}s^3 + \mathcal O\left(s^5\right), 
\end{array}
\right.
\end{equation}
where \( \mathfrak q \) is a free parameter and we used the identity \( H^\prime(s)=-q^2(s) \) to find the expansion of \( H(s) \) with the constant term found through its original definition. Hence, in this case it is also true that
\[
\left\{
\begin{array}{rcl}
B_2(s) & = & \displaystyle \frac{\mathfrak q}{2s} - \frac{\nu+1}8s + \mathcal O\left(s^2\right), \medskip \\ 
A_2(s) & = & \displaystyle \frac{\mathfrak q}{2s}  + \frac{\mathfrak q^2}2 + \frac{\nu+1}8s + \mathcal O\left(s^2\right),
\end{array}
\right.
\]
and
\begin{equation}
\label{pIIpole2}
B_3(s) = \frac{\mathfrak q^2}{4s} + \mathcal O(s) \qandq A_3(s)  =  -\frac{\mathfrak q^2}{4s} + \mathcal O(1).
\end{equation}

\section{Riemann-Hilbert Problem for Painlev\'e XXXIV} 
\label{ap:34}

The primary source for the material of this appendix is \cite{IKOs08}. 

\subsection{Riemann-Hilbert Problem} We retain the notation of the previous appendix. In what follows, we always assume that
\[
\tau = \ic 2^{-1/3}\zeta^{1/2} \qandq s= -2^{1/3}x,
\]
where \( \arg\big(\zeta\big)\in(-\pi,\pi) \) and the square root is principal (that is, \( \zeta_{\pm}^{1/2} = \pm\ic|\zeta|^{1/2} \) for \( \zeta\in(-\infty,0) \)). Set
\begin{equation}
\label{UV}
\left\{
\begin{array}{l}
U_n(x) := A_n(s) + (-1)^n B_n(s), \medskip \\ 
V_n(x) := A_n(s) - (-1)^n B_n(s).
\end{array}
\right.
\end{equation}

Given \( \alpha>-1 \), let \( q \) be a solution of \eqref{p2} with \( \nu=\alpha+1/2 \) and \( x \) such that \( q(s) \) is finite. Let \( \boldsymbol\Psi_{\alpha+1/2}^{FN}(\tau;s) \) be given by \eqref{FN} and \( s_1,s_2,s_3 \) be the corresponding Stokes parameters satisfying \eqref{pII-stokes}. Set
\begin{equation}
\label{p2-34}
\boldsymbol\Phi_\alpha(\zeta;x) := \left( \begin{matrix} 1 & 0 \\  -2^{1/3}\ic V_1(x) & 1 \end{matrix} \right) \frac{\zeta^{-\sigma_3/4}}{\sqrt2} \begin{pmatrix} 1 & \ic \\  \ic & 1 \end{pmatrix} e^{\pi\ic\sigma_3/4} \boldsymbol\Psi_{\alpha+1/2}^{FN}(\tau;s) e^{-\pi\ic\sigma_3/4}.
\end{equation}
Then \( \boldsymbol\Phi_\alpha(\zeta;x) \) solves the following Riemann-Hilbert problem (\rhphi):

\begin{itemize}
\label{rhphi}
\item[(a)] $\boldsymbol\Phi_\alpha(\zeta;x)$ is holomorphic in $\C\setminus\big(\Sigma_2\cup \Sigma_3\cup(-\infty,\infty)\big)$, where
\[
\Sigma_k := \left\{e^{2(k-1)\pi\ic/3}x:x\in(0,\infty)\right\}, \quad k\in\{2,3\},
\]
are oriented towards the origin (the real line and its subsets are oriented from left to right as standard);
\item[(b)] $\boldsymbol\Phi_\alpha(\zeta;x)$ has continuous traces on $(-\infty,0)\cup(0,\infty)\cup\Sigma_2\cup\Sigma_3$ that satisfy
\[
\boldsymbol\Phi_{\alpha+}(\tau;x) = \boldsymbol\Phi_{\alpha-}(\tau;x)
\left\{
\begin{array}{rl}
\left(\begin{matrix} 0 & 1 \\ -1 & 0 \end{matrix}\right), &  \tau\in(-\infty,0), \medskip \\
\left(\begin{matrix} 1 & b_1 \\ 0 & 1 \end{matrix}\right), & \tau\in(0,\infty), \medskip \\
\left(\begin{matrix} 1 & 0 \\ b_k & 1 \end{matrix}\right), & \tau\in\Sigma_k, ~~ k\in\{2,3\},
\end{array}
\right.
\]
where \( b_1 := \ic s_2 \), \( b_2 = \ic s_3 \), and \( b_3 = \ic s_1 \) (\( b_1+b_2+b_3 - b_1b_2b_3=2\cos(\alpha\pi) \));
\item[(c)] as $\zeta\to0$  it holds that
\[
\boldsymbol\Phi_\alpha(\zeta;x) = \boldsymbol H(\zeta;x) \zeta^{\alpha\sigma_3/2}\left(\begin{matrix} 1 & \ic\kappa_{\alpha+1/2} \log\zeta \\ 0 & 1 \end{matrix}\right)\boldsymbol A_j, \quad \zeta\in Q_j,
\]
where \( Q_j \) is the connected component of \( \C\setminus\big(\Sigma_2\cup \Sigma_3\cup(-\infty,\infty)\big) \) contained in or containing the \( j \)-th quadrant, $\boldsymbol H(\zeta;x)$ is holomorphic around \( \zeta=0 \), and \( \boldsymbol A_j= e^{\pi\ic\sigma_3/4} \boldsymbol A_j^{FN} e^{-\pi\ic\sigma_3/4} \);
\item[(d)] $\boldsymbol\Phi_\alpha(\zeta;x)$ has the following asymptotic expansion near $\infty$:
\begin{multline*}
\boldsymbol\Phi_\alpha(\zeta;x) =  \left(\boldsymbol I + \sum_{n=1}^\infty \frac{2^{2n/3}}{\zeta^n} \left(\begin{matrix} U_{2n}(x) & -2^{-1/3}\ic U_{2n-1}(x) \\ 2^{1/3}\ic(V_{2n+1}-U_{2n}V_1)(x) & (V_{2n}-V_1U_{2n-1})(x) \end{matrix}\right) \right) \\ \times \frac{\zeta^{-\sigma_3/4}}{\sqrt2} \left(\begin{matrix} 1 & \ic \\ \ic & 1 \end{matrix}\right) \exp\left\{-\frac23 \left(\zeta^{3/2}+x\zeta^{1/2}\right)\sigma_3\right\},
\end{multline*}
which holds uniformly in $\{|\zeta|\geq1\}\setminus\big(\Sigma_2\cup \Sigma_3\cup(-\infty,\infty)\big)$.
\end{itemize}

The only claim that is not contained in \cite{IKOs08} is the explicit form of the series in \hyperref[rhphi]{\rhphi}(d). The latter follows easily from \eqref{pII-exp} since
\begin{multline*}
\frac1{\sqrt2}\left(\begin{matrix} 1 & \ic \\ \ic & 1 \end{matrix}\right) e^{\pi\ic\sigma_3/4} \left( \begin{matrix} A_n(s) & B_n(s) \smallskip \\ (-1)^nB_n(s) & (-1)^nA_n(s) \end{matrix} \right) e^{-\pi\ic\sigma_3/4} \frac1{\sqrt2}\left(\begin{matrix} 1 & -\ic \\ -\ic & 1 \end{matrix}\right) \\ = \frac12\left( \begin{matrix} \big(1+(-1)^n\big) U_n(x) & -\ic \big(1-(-1)^n\big) U_n(x) \smallskip \\ \ic \big(1-(-1)^n\big) V_n(x) & \big(1+(-1)^n\big) V_n(x) \end{matrix} \right)
\end{multline*}
and for \( n=2m \) it holds that
\[
\zeta^{-\sigma_3/4}\frac{2^{2m/3}}{\zeta^m}\left( \begin{matrix} U_{2m}(x) & 0 \smallskip \\ 0 &  V_{2m}(x) \end{matrix} \right)\zeta^{\sigma_3/4} = \frac{2^{2m/3}}{\zeta^m}\left( \begin{matrix} U_{2m}(x) & 0 \smallskip \\ 0 &  V_{2m}(x) \end{matrix} \right)
\]
while for \( n=2m-1 \) it holds that
\begin{multline*}
\zeta^{-\sigma_3/4}\frac{2^{(2m-1)/3}}{\zeta^{m-1/2}}\left( \begin{matrix} 0 & -\ic U_{2m-1}(x) \smallskip \\ \ic V_{2m-1}(x) &  0 \end{matrix} \right)\zeta^{\sigma_3/4} = \\  \frac{2^{2(m-1)/3}}{\zeta^{m-1}}\left( \begin{matrix} 0 & 0 \smallskip \\ \ic 2^{1/3} V_{2(m-1)+1}(x) & 0 \end{matrix} \right) + \frac{2^{2m/3}}{\zeta^m}\left( \begin{matrix} 0 & -\ic 2^{-1/3} U_{2m-1}(x) \smallskip \\ 0 & 0 \end{matrix} \right)
\end{multline*}
(notice that the identity matrix plus the first matrix on the right-hand side of the last equality for \( m=1 \) is exactly the inverse of the first matrix in \eqref{p2-34}).

\subsection{Lax Pair} Define
\begin{equation}
\label{Uu}
\left\{
\begin{array}{l}
\mathcal U(x) := 2^{-2/3}(H-q)(s), \medskip \\
u(x) := \mathcal U^\prime(x) - x/2 = 2^{-1/3}\left(q^\prime(s)+q^2(s)+s/2\right).
\end{array}
\right.
\end{equation}
Then it can be verified that \( u(x) \) is a solution of Painlev\'e XXXIV equation with parameter \( \alpha^2 \):
\begin{equation}
\label{p34}
u^{\prime\prime}(x) = 4u^2(x) + 2xu(x) + \frac{(u^\prime(x))^2-\alpha^2}{2u(x)}.
\end{equation}
It has been shown in \cite[Lemma~3.3]{IKOs08} that  \( \boldsymbol\Phi_\alpha(\zeta;x) \) satisfies
\begin{equation}
\label{p34a}
\left\{
\begin{array}{l}
\displaystyle \partial_\zeta \boldsymbol\Phi_\alpha(\zeta;x) = \boldsymbol A(\zeta;x)\boldsymbol\Phi_\alpha(\zeta;x), \bigskip \\
\displaystyle \partial_x \boldsymbol\Phi_\alpha(\zeta;x) = \boldsymbol U(\zeta;x)\boldsymbol\Phi_\alpha(\zeta;x),
\end{array}
\right.
\end{equation}
in each sector \( Q_j \), where 
\begin{multline}
\label{p34b}
\boldsymbol A(\zeta;x) =  \left(\begin{matrix} \mathcal U(x) & \ic \smallskip \medskip\\ -\ic\big(\zeta + u(x) + x - \mathcal U^2(x) \big) &  -\mathcal U(x) \end{matrix}\right) +  \\
\frac1\zeta\left(\begin{matrix} u^\prime(x)/2-(u\mathcal U)(x) & -\ic u(x) \medskip \\ \displaystyle -\ic \left((\mathcal U^2u)(x) - (\mathcal Uu^\prime)(x) + \frac{(u^\prime(x))^2-\alpha^2}{4u(x)}\right) & (u\mathcal U)(x) - u^\prime(x)/2 \end{matrix}\right)
\end{multline}
and
\begin{equation}
\label{p34c}
\boldsymbol U(\zeta;x) := \left(\begin{matrix} \mathcal U(x) & \ic \smallskip \medskip\\ \ic\big(-\zeta + \mathcal U^2(x) - \mathcal U^\prime(x) \big) &  -\mathcal U(x) \end{matrix}\right).
\end{equation}
Expressions \eqref{p34b} and \eqref{p34c} are not explicitly stated in \cite[Lemma~3.3]{IKOs08}. What follows from the lemma is that \( \boldsymbol A(\zeta;x) \) is equal to
\[
\left( \begin{matrix} 1 & 0 \medskip\\  -2^{1/3}\ic V_1(x) & 1 \end{matrix} \right) \left( \begin{matrix} -2^{1/3}q(s) + \alpha/(2\zeta) &  \ic(1-u(x)/\zeta) \medskip\\ -\ic\big(\zeta-u(x)+2^{2/3}q^\prime(s)\big)  & 2^{1/3}q(s)-\alpha/(2\zeta) \end{matrix} \right) \left( \begin{matrix} 1 & 0 \medskip\\  2^{1/3}\ic V_1(x) & 1 \end{matrix} \right)
\]
and that \( \boldsymbol U(\zeta;x) \) is equal to
\[
\left[ \left( \begin{matrix} 0 & 0 \medskip\\  -2^{1/3}\ic V_1^\prime(x) & 0 \end{matrix} \right) + \left( \begin{matrix} 1 & 0 \medskip \\  -2^{1/3}\ic V_1(x) & 1 \end{matrix} \right) \left( \begin{matrix} -2^{1/3}q(s) & \ic \medskip\\  -\ic\zeta & 2^{1/3}q(s) \end{matrix} \right) \right] \left( \begin{matrix} 1 & 0 \medskip\\  2^{1/3}\ic V_1(x) & 1 \end{matrix} \right),
\]
which yield \eqref{p34b} and \eqref{p34c} (identity \( \alpha - u^\prime(x) = 2^{4/3}u(x)q(s) \) is useful in deriving \eqref{p34b}; recall also that \( -2V_1(x)=q(s)+H(s) \) by \eqref{AB12} and \eqref{UV}).

As mentioned in the previous subsection, it follows from the general theory of Riemann-Hilbert problems that the solution of  \hyperref[rhphi]{\rhphi} is a meromorphic function of \( x \). In fact, if \( s=-2^{1/3}x \) is not a pole of \( q \), then the solution of \hyperref[rhphi]{\rhphi} exists and is given by \eqref{p2-34}. Moreover, if \( s_0 \) is a pole of \( q \) with the residue \( 1 \), that is, \( q_{-1}=1 \) in \eqref{pIIpole}, then \eqref{pIIpole} yields that the functions \( \mathcal U(x) \) and \( u(x) \) are, in fact, analytic at \( x_0=-2^{-1/3}s_0 \). In this case, the matrices \( \boldsymbol A(\zeta;x) \) and \( \boldsymbol U(\zeta;x) \) are analytic at \( x_0 \) as well (notice that the fraction in \( (2,1) \)-entry of \(  \boldsymbol A(\zeta;x) \) is equal to \( u^{\prime\prime}(x)/2 - 2u^2(x) - xu(x) \) by \eqref{p34}).Then we get from the second relation in \eqref{p34a} that \( \boldsymbol \Phi_\alpha(\zeta;x) \) also must be analytic at \( x_0 \). On the other hand, if \( q_{-1} = -1 \), then \( x_0 \) is a simple pole of \( \mathcal U(x) \) and a double pole of \( u(x) \). Moreover, the Riemann-Hilbert problem \hyperref[rhphi]{\rhphi} cannot be solvable since it follows from \hyperref[rhphi]{\rhphi}(d) that
\begin{equation}
\label{ant-u}
\lim_{\zeta\to\infty} \zeta\left[ \boldsymbol \Phi_\alpha(\zeta;x) \exp\left\{\frac23 \left(\zeta^{3/2}+xz^{1/2}\right)\sigma_3\right\} \left(\begin{matrix} 1 & -\ic \\ -\ic & 1 \end{matrix}\right) \frac{\zeta^{\sigma_3/4}}{\sqrt2}   \right]_{12} = -2^{1/3}\ic U_1(x) = \ic \mathcal U(x)
\end{equation}
(analyticity of \( \boldsymbol \Phi_\alpha(\zeta;x) \) at \( x_0 \) implies analyticity of \( \mathcal U(x) \) at \( x_0 \)).

\subsection{Modified Riemann-Hilbert Problem}

Assume now that \( \alpha \) and Stokes parameters \( b_1,b_2,b_3 \) from \hyperref[rhphi]{\rhphi}(b) are such that \( x_0=0 \) is a pole of \( \boldsymbol \Phi_\alpha(\zeta;x) \) (hence, \(s_0=0 \) is a pole of \( q(s) \) with residue \( -1 \)). Note that in this case
\[
\left\{
\begin{array}{rcl}
U_1(x) & = & \displaystyle 2^{-1/3}x^{-1} - \mathfrak q + \mathcal O\big(x^3\big), \medskip \\
U_2(x) & = & \displaystyle -2^{-1/3}\mathfrak qx^{-1} + \mathfrak q^2/2 + \mathcal O\big(x^2\big), \medskip \\
U_3(x) & = & \displaystyle 2^{-4/3}\mathfrak q^2x^{-1}  + \mathcal O(1),
\end{array}
\right.
\]
as evident by \eqref{pIIpole0}--\eqref{pIIpole2} and \eqref{UV}. Recall that \( \mathcal U(x) = -2^{1/3}U_1(x) \) and set
\[
\left\{
\begin{array}{rcl}
\displaystyle \mathcal S(x) &:=& 2^{2/3}\big(U_2+U_3U_1^{-1}\big)(x), \medskip \\
\displaystyle \mathcal W(x) &:=& (\mathcal U^2-\mathcal U^\prime - \mathcal S)(x).
\end{array}
\right.
\]
 One can readily verify that
\begin{equation}
\label{USW}
\left\{
\begin{array}{rcl}
\displaystyle \mathcal U(x) &=& -x^{-1} + 2^{1/3}\mathfrak q + \mathcal O\left(x^3\right), \medskip \\
\displaystyle \mathcal S(x) &=& -2^{1/3}\mathfrak q x^{-1} + 2^{2/3}\mathfrak q^2 + S_*x + \mathcal O\left(x^2\right), \medskip \\
\displaystyle \mathcal W(x) &=& -2^{1/3}\mathfrak qx^{-1} - S_*x + \mathcal O\left(x^2\right),
\end{array}
\right.
\end{equation}
as \( x\to0 \) for some constant \( S_* \) that we can avoid computing. To get rid of the pole of \( \boldsymbol \Phi_\alpha(\zeta;x) \) at \( x_0=0 \), let us define
\begin{equation}
\label{Phi-special}
\boldsymbol \Phi_\alpha^\mathsf{sing}(\zeta;x) := \left( \begin{matrix} \displaystyle \zeta - \mathcal S(x) & -\ic\mathcal U(x) \medskip \\ -\ic\mathcal U^{-1}(x) & 0 \end{matrix} \right)\boldsymbol \Phi_\alpha(\zeta;x).
\end{equation}
Denote the prefactor above by \( \boldsymbol S(\zeta;x) \) (the way it can be found will become clear at the end of this subsection). Trivially, 
\begin{equation}
\label{PhisUs}
\partial_x \boldsymbol\Phi_\alpha^\mathsf{sing}(\zeta;x) = \boldsymbol U^\mathsf{sing}(\zeta;x)\boldsymbol\Phi_\alpha^\mathsf{sing}(\zeta;x), \quad \boldsymbol U^\mathsf{sing}:= \partial_x\boldsymbol S\boldsymbol S^{-1}+\boldsymbol S\boldsymbol U\boldsymbol S^{-1}.
\end{equation}
Since the matrix \( \boldsymbol S(\zeta;x) \) has determinant identically equal to \( 1 \), it is a lengthy but straightforward computation to find that
\[
\boldsymbol U^\mathsf{sing}(\zeta;x) = \ic\left( \begin{matrix} \ic(\zeta+\mathcal W(x))\mathcal U^{-1}(x) & \zeta^2 + \zeta\big(\mathcal W-\mathcal S\big)(x) + \big( \mathcal W^2 +\mathcal U^\prime\mathcal W-\mathcal U\mathcal S^\prime\big)(x) \medskip \\ \mathcal U^{-2}(x) &  -\ic(\zeta + \mathcal W(x))\mathcal U^{-1}(x) \end{matrix}\right).
\]
It readily follows from \eqref{USW} that \( \boldsymbol U^\mathsf{sing}(\zeta;x) \) is in fact analytic at \( x_0=0 \). Hence, it follows from \eqref{PhisUs} that \( \boldsymbol\Phi_\alpha^\mathsf{sing}(\zeta;x) \) is analytic at \( x_0=0 \) as well. It further follows from \hyperref[rhphi]{\rhphi}, the definition of \( \boldsymbol \Phi_\alpha^\mathsf{sing}(\zeta;x) \) in \eqref{Phi-special}, the above proven analyticity in the parameter \( x \) in some neighborhood of \( 0 \), and the explanation further below that \( \boldsymbol \Phi_\alpha^\mathsf{sing}(\zeta;x) \) solves the following Riemann-Hilbert problem (\rhphis):

\begin{itemize}
\label{rhphis}
\item[(a,b,c)] \( \boldsymbol \Phi_\alpha^\mathsf{sing}(\zeta;x) \) satisfies \hyperref[rhphi]{\rhphi}(a,b,c);
\item[(d)] $\boldsymbol\Phi_\alpha^\mathsf{sing}(\zeta;x)$ has the following behavior near $\infty$:
\[
\boldsymbol\Phi_\alpha^\mathsf{sing}(\zeta;x) =  \left(\boldsymbol I+\boldsymbol{\mathcal{O}}\Big(\zeta^{-1}\Big)\right) \frac{\zeta^{3\sigma_3/4}}{\sqrt2} \left(\begin{matrix} 1 & \ic \\ \ic & 1 \end{matrix}\right) \exp\left\{ -\frac23 \left( \zeta^{3/2}+x\zeta^{1/2}\right)\sigma_3\right\}
\]
uniformly in $\C\setminus\big(\Sigma_2\cup \Sigma_3\cup(-\infty,\infty)\big)$.
\end{itemize}

Indeed, since \( \boldsymbol S(\zeta;x) \) is an analytic (in fact, linear) in \( \zeta \), \hyperref[rhphis]{\rhphis}(a,b,c) are obvious. It also holds  that
\begin{multline*}
\left(\begin{matrix} \zeta & 0 \\ 0 & 0 \end{matrix}\right) \left(\boldsymbol I + \sum_{n=1}^\infty \frac{2^{2n/3}}{\zeta^n} \left(\begin{matrix} * & -2^{-1/3}\ic U_{2n-1}(x) \\ * & * \end{matrix}\right) \right) = 
\\  \left(\begin{matrix} \zeta & 0 \\ 0 & 0 \end{matrix}\right) + \left(\begin{matrix} * & -2^{1/3}\ic U_1(x) \\ 0 & 0 \end{matrix}\right) + \frac1\zeta\left(\begin{matrix} * & -2\ic U_3(x) \\ 0 & 0 \end{matrix}\right) + \left(\begin{matrix} \mathcal O\big(\zeta^{-2}\big) & \mathcal O\big(\zeta^{-2}\big) \\ 0 & 0 \end{matrix}\right) = 
\\ \left( \zeta\left(\begin{matrix} 0 & -2^{1/3}\ic U_1(x) \\ 0 & 0 \end{matrix}\right) + \left(\begin{matrix} 1 & -2\ic U_3(x) \\ 0 & 0 \end{matrix}\right) + \left(\begin{matrix} \mathcal O\big(\zeta^{-1}\big) & \mathcal O\big(\zeta^{-1}\big) \\ 0 & 0 \end{matrix}\right) \right) \zeta^{\sigma_3}.
\end{multline*}
Furthermore, since \( -U_2(x)=(V_2-V_1U_1)(x) \), which can be verified using \eqref{AB12} and \eqref{UV}, it holds  that
\begin{multline*}
\left(\begin{matrix} - 2^{2/3}\big(U_2+U_3U_1^{-1}\big)(x) & 2^{1/3}\ic U_1(x) \\ 2^{-1/3}\ic U_1^{-1}(x) & 0 \end{matrix}\right) \left(\boldsymbol I + \sum_{n=1}^\infty \frac{2^{2n/3}}{\zeta^n} \left(\begin{matrix} * & -2^{-1/3}\ic U_{2n-1}(x) \\ * & (V_{2n}-V_1U_{2n-1})(x) \end{matrix}\right) \right) = 
\\ \left(\begin{matrix} * & 2^{1/3}\ic U_1(x) \\ 2^{-1/3}\ic U_1^{-1}(x) & 0 \end{matrix}\right) + \frac1\zeta\left(\begin{matrix} * & 2\ic U_1(x)\big(U_2+U_3U_1^{-1}+V_2-V_1U_1\big)(x) \\ * & 1 \end{matrix}\right) + \boldsymbol{\mathcal O}\left(\frac1{\zeta^2}\right) = 
\\ \left( \zeta\left(\begin{matrix} 0 & 2^{1/3}\ic U_1(x) \\ 0 & 0 \end{matrix}\right) + \left(\begin{matrix} 0 & 2\ic U_3(x) \\ 0 & 1 \end{matrix}\right) + \frac1\zeta\left(\begin{matrix} * & * \\ 2^{-1/3}\ic U_1^{-1}(x) & * \end{matrix}\right) + \boldsymbol{\mathcal O}\left(\frac1{\zeta^2}\right) \right)\zeta^{\sigma_3}
\end{multline*}
(these are exactly the relations that define \( \boldsymbol S(\zeta;x) \)). \hyperref[rhphis]{\rhphis}(d) now follows from \hyperref[rhphi]{\rhphi}(d) and the definition of \( \boldsymbol S(\zeta;x) \) in \eqref{Phi-special}. Since \( U_1^{-1}(0)=0 \), the last two computations also show that
\begin{multline}
\label{Phis-carried}
\boldsymbol \Phi_\alpha^\mathsf{sing}(\zeta;0) = \left(\boldsymbol I + \frac1\zeta \left(\begin{matrix} * & * \\ 0 & * \end{matrix}\right) + \boldsymbol{\mathcal{O}}\Big(\zeta^{-2}\Big)\right) \frac{\zeta^{3\sigma_3/4}}{\sqrt2} \left(\begin{matrix} 1 & \ic \\ \ic & 1 \end{matrix}\right) \exp\left\{ -\frac23 \zeta^{3/2}\sigma_3\right\}
\\ = \frac{\zeta^{3\sigma_3/4}}{\sqrt2} \left(\begin{matrix} 1 & \ic \\ \ic & 1 \end{matrix}\right) \left(\boldsymbol I + \boldsymbol{\mathcal{O}}\Big(\zeta^{-1/2}\Big)\right) \exp\left\{ -\frac23 \zeta^{3/2}\sigma_3\right\}
\end{multline}
uniformly in $\{|\zeta|\geq1\}\setminus\big(\Sigma_2\cup \Sigma_3\cup(-\infty,\infty)\big)$.

\section{An Example of Weights in \( W^\mathsf{sing} \)}
\label{ap:example}

In this appendix we discuss a group of examples of weights in \( W^\mathsf{sing} \) for which \( \N_\varepsilon^\mathsf{sing} = \varnothing \) and therefore our results do not apply. As follows from Proposition~\ref{prop:lo}, this can happen only on \( \Delta_\mathsf{sym} \).

\subsection{Orthogonal Polynomials on a Segment} Let \( \widehat\rho(z) \) be an analytic and non-vanishing function in some neighborhood of \( [-1,0] \) and \( \widehat Q_n(z) \) be the non-identically zero monic polynomial of smallest degree, which is necessarily at most \( n \), such that
\begin{equation}
\label{ex-1}
\int_{-1}^0 x^k\widehat Q_n(x)\widehat\rho(x)\dd x =0, \quad k\in\{0,\ldots,n-1\}
\end{equation}
(we could also introduce Jacobi-type singularities at \( -1 \) and \( 0 \) into the weight \( \widehat\rho(x) \), but choose not to for simplicity of the exposition). Such polynomials were studied in \cite{KMcLVAV04} using Riemann-Hilbert approach under an additional assumption of positivity on \( [-1,0] \) (this assumption is really not needed for the approach to work). To describe the asymptotics of these polynomials, let
\[
\widehat\Phi(z) := 2z+1+2\widehat w(z), \quad \widehat w(z):=\sqrt{z(z+1)},
\]
where the branch of the square root is chosen so that \( \widehat w(z) \) is analytic in \( \C\setminus[-1,0] \) and \( \widehat w(z)=z+\mathcal O(1) \) as \( z\to\infty \), while \( \widehat\Phi(z) \) is nothing but the conformal map of \( \overline\C\setminus[-1,0] \) onto \( \big\{|z|>1\big\} \) such that \( \widehat\Phi(\infty)=\infty \) and \( \widehat\Phi^\prime(\infty)>0 \). Further, let
\begin{equation}
\label{ex-2a}
D_{\widehat\rho}(z) := \exp\left\{ \frac{\widehat w(z)}{2\pi\ic} \int_{-1}^0 \frac{\log\widehat\rho(x)}{z-x}\frac{\dd x}{\widehat w_+(x)} \right\}, \quad z\not\in[-1,0],
\end{equation}
which is simply the Szeg\H{o} function of \( \widehat\rho(x) \) (one must choose a continuous determination of \( \log\widehat\rho(x) \) on \( [-1,0] \); Szeg\H{o} functions for different continuous determinations will either coincide or differ by a sign), i.e., it is a non-vanishing and holomorphic function in the domain of its definition with continuous traces on \( [-1,0] \) that satisfy \( D_{\widehat\rho+}(x)D_{\widehat\rho-}(x) = 1/\widehat\rho(x) \). Lastly, let
\[
D(z) := \left( \widehat\Phi(z)/\widehat w(z) \right)^{1/2}
\] 
be the branch holomorphic in \( \overline\C\setminus[-1,0] \) that is positive for \( z>0 \) (this is the Szeg\H{o} function of \( \widehat w_+(x) \), that is, \( D(z) \) is analytic and non-vanishing in \( \overline\C\setminus[-1,0] \) and \( D_+(x)D_-(x) = 1/\widehat w_+(x) \) for \( x\in(-1,0) \)). Then it holds that
\begin{equation}
\label{ex-2}
\widehat Q_n(z) = \big(1+\mathcal O\big(1/n\big)\big) \widehat\gamma_n \widehat\Phi^n(z)\big( DD_{\widehat\rho}\big)(z)
\end{equation}
locally uniformly in \( \overline\C\setminus[-1,0] \), where \( \widehat\gamma_n^{-1} = 4^n(DD_{\widehat\rho})(\infty) \) is the normalizing constant that makes the right-hand side of \eqref{ex-2} behave like \( z^n + \mathcal O(z^{n-1}) \) around infinity. The reader should see clear parallels between \eqref{ex-2} and \eqref{Psin}, \eqref{Asymp1} (the product \( (DD_{\widehat\rho})(z) \) also can be defined via an integral representation \eqref{ex-2a} with \( \widehat\rho(x) \) replaced by \( (\widehat\rho \widehat w_+)(x) \)).

\subsection{Rotationally Symmetric Weights} Let \( \Delta_\mathsf{sym} \) be given by \eqref{Delta_sym}. In this case we let \( a_1:=-1 \), \( a_2=-\eta \), and \( a_3=-\eta^2 \), where \( \eta:=e^{2\pi\ic/3} \). Define a weight function \( \rho(s) \) on \( \Delta_\mathsf{sym} \) by setting \(  \rho(s) := \widehat\rho\big(s^3\big) \), \( s\in\Delta_\mathsf{sym} \). Let \( Q_n(s) \) be the minimal degree non-identically zero polynomial satisfying \eqref{ortho} with the just defined weight \( \rho(s) \). Then
\[
Q_{3n}(z) = Q_{3n+1}(z) = Q_{3n+2}(z) = \widehat Q_n\big(z^3\big).
\]
Indeed, as all the legs of \( \Delta_\mathsf{sym} \) are oriented towards the origin, it holds for \( k\leq 3n+1 \) that
\[
\int_{\Delta_\mathsf{sym}} s^k\widehat Q_n\big(s^3\big)\rho(s)\dd s = \int_{-1}^0\left(1+\eta^{k+1} + \eta^{2(k+1)}\right)x^k\widehat Q_n\big(x^3\big) \widehat\rho\big(x^3\big)\dd x.
\]
The above expression is equal to \( 0 \) when \( k\neq 3m+2 \) due to the sum \( 1+\eta^{k+1} + \eta^{2(k+1)} \). When \( k = 3m+2 \) we get that  
\[
\int_{\Delta_\mathsf{sym}} s^{3m+2}\widehat Q_n\big(s^3\big)\rho(s)\dd s = 3\int_{-1}^0x^{3m+2}\widehat Q_n\big(x^3\big) \widehat\rho\big(x^3\big)\dd x = \int_{-1}^0x^m\widehat Q_n(x)\widehat\rho(x)\dd x =0,
\]
where the last conclusion holds by \eqref{ex-1} and one needs to observe that \( 3m+2=k\leq 3n+1 \) implies that \( m\leq n-1 \). Thus, it follows from \eqref{ex-2} that
\begin{equation}
\label{ex-4}
Q_{3n}(z) =  \widehat\gamma_n(1+o(1))\widehat\Phi^n\big(z^3\big) \big( DD_{\widehat\rho}\big)\big(z^3\big)
\end{equation}
locally uniformly in \( \overline\C\setminus\Delta_\mathsf{sym} \). Let \( \RS \) be the Riemann surface of \( w(z) = \sqrt{z(z^3+1)} \) as defined in Section~\ref{ssec:RS} and \( \bd \) be the lift of \( \Delta_\mathsf{sym} \) to this surface, where \( w(z) \) is the branch from \eqref{w}. As before, \( \bd \) is the boundary of the sheets \( \RS^{(0)} \) and \( \RS^{(1)} \). Set
\[
\widetilde\Psi(\z) := \left( \widehat\Phi^n\big(z^3\big) \big( DD_{\widehat\rho}\big)\big(z^3\big) \right)^{(-1)^k}, \quad \z\in\RS^{(k)}\setminus\bd.
\]
Then \( \widetilde\Psi(\z) \) is a sectionally meromorphic function on \( \RS\setminus\bd \) that has a pole of order \( 3n \) at \( \infty^{(0)} \), a zero of order \( 3n \) at \( \infty^{(1)} \), and otherwise is non-vanishing and finite. These are exactly the properties exhibited by \( \widetilde\Psi_n(\z) \) from Theorem~\ref{thm:BANS}. However, it holds that
\[
\widetilde\Psi_+(\s) = \widetilde\Psi_-(\s)/\big(\widehat\rho\widehat w_+\big)\big(s^3\big) = \widetilde\Psi_-(\s)/\big(s(\rho w_+)(s)\big)
\]
for \( \s\in\bd \). That is the jump relations satisfied by \( \widetilde\Psi(\s) \) are different from the ones in \eqref{BVP-Psi}. Hence, formula \eqref{ex-4} is distinct from \eqref{Asymp1}. In what follows we find \( \N^\mathsf{reg} \) (as we explained in Proposition~\ref{prop:jip} in this case the sequence is \( 3 \)-periodic and therefore is independent of \( \varepsilon\)), show that the above weight \( \rho(s) \) belongs to \( W^\mathsf{sing} \), and that \( \N^\mathsf{sing}=\varnothing \). Thus, formulae \eqref{Asymp1} are indeed not applicable.

\subsection{Sequence \( \N^\mathsf{reg} \)} 

As before, we shall write \( \Delta_\mathsf{sym} = \Delta_1\cup\Delta_2\cup\Delta_3 \), where \( \Delta_i \) has endpoints \( a_i=-\eta^{i-1} \) and \( 0 \). It can be readily checked that
\begin{equation}
\label{ex-5}
w_+(x) = |w(x)| \left\{ \begin{array}{rl} e^{3\pi\ic/2}, & x\in\Delta_1, \\ e^{-7\pi\ic/6}, & x\in\Delta_2, \\ e^{\pi\ic/6}, & x\in\Delta_3, \end{array}\right. =:  |w(x)|e^{a_w(x)}.
\end{equation}
One can also check that \( \boldsymbol\Delta_1 \), \( \boldsymbol\Delta_2 \), and \( \boldsymbol\Delta_3 \) are homologous to \( \boldsymbol\alpha+\boldsymbol\beta \), \( - \boldsymbol\beta \), and \( -\boldsymbol\alpha \), respectively, see Figure~\ref{fig:torus}. Therefore,
\[
\oint_{\boldsymbol\alpha}\frac{\dd s}{w(\s)} = -2\int_{\Delta_3}\frac{\dd s}{w_+(s)} = -2e^{-\pi\ic/6}\int_{\Delta_3}\frac{\dd s}{|w(s)|} = 2e^{\pi\ic/6}\int_{-1}^0\frac{\dd x}{|w(x)|},
\] 
where we used \eqref{ex-5} for the second equality. Similarly, it holds that
\[
\oint_{\boldsymbol\beta}\frac{\dd s}{w(\s)} = -2\int_{\Delta_2}\frac{\dd s}{w_+(s)} = 2e^{\pi\ic/6}\int_{\Delta_2}\frac{\dd s}{|w(s)|} = 2e^{5\pi\ic/6}\int_{-1}^0\frac{\dd x}{|w(x)|}.
\]
That is, we have that \( \mathsf B=e^{2\pi\ic/3} \), see \eqref{B}. Moreover, we get from \eqref{prop3a} and the fact that \( I_i=1/3 \) for \( \Delta_\mathsf{sym} \) that \( \omega=-\tau= 2/3 \). In particular, \( \Phi^3(\z) \) must be a rational function on the whole surface \( \RS \), see \eqref{Phi-jump}. In fact, it is not hard to see that
\[
\Phi^3\big(z^{(k)}\big) = \widehat\Phi^{(-1)^k}\big(z^3 \big) = \left(2z^3+1+2zw(z)\right)^{(-1)^k}.
\]

Define \( \sqrt w(z) \) to be the branch holomorphic in \( \C\setminus\Delta_\mathsf{sym} \) satisfying \( \sqrt w(z) = z + \mathcal O(1) \) as \( z\to\infty \). It holds that
\[
\sqrt w_\pm(s) = \sqrt{|w(s)|} \left\{ \begin{array}{rl} e^{\pm3\pi\ic/4}, & s\in\Delta_1, \\ e^{-\pi\ic/3\mp\pi\ic/4}, & s\in\Delta_2, \\ e^{\pi\ic/3\mp\pi\ic/4}, & s\in\Delta_3. \end{array}\right.
\]
Therefore, we have that \( -\ic \sqrt w_+(s)\sqrt w_-(s) =w_+(s) \), \( s\in\Delta_\mathsf{sym} \), see \eqref{ex-5}. Recall \eqref{a}. Let
\begin{equation}
\label{ex-9}
S\big(z^{(k)}\big) := \exp \left\{-2\pi\ic(\tau+2) a\big(z^{(k)}\big) \right\} \Phi\big(z^{(k)}\big) \left(e^{-\pi\ic/4}\sqrt w(z)\right)^{-(-1)^k},
\end{equation}
for \( z\in\C\setminus(\Delta_\mathsf{sym}\cup\pi(\boldsymbol \alpha)\cup \pi(\boldsymbol \beta)) \), \( k\in\{0,1\} \). This function is holomorphic and non-vanishing in the domain of its definition. Notice that it extends holomorphically to \( \infty^{(0)} \) and \( \infty^{(1)} \) and has non-zero values there. We also get from \eqref{Phi-jump} and \eqref{a-jump} that \( S(\z) \) is holomorphic and non-vanishing \( \RS\setminus(\boldsymbol\Delta\cup\boldsymbol \alpha) \) and satisfies
\[
S_+(\s) = S_-(\s) \left\{ \begin{array}{rl} \exp\{2\pi\ic(\omega+\mathsf B\tau+2\mathsf B)\}, & \s\in\boldsymbol\alpha\setminus \boldsymbol A, \smallskip \\ 1/w_+(s), & \s\in\boldsymbol\Delta\setminus\boldsymbol A. \end{array}\right.
\]
One should observe as well that \( S(\z)S(\z^*) \equiv 1 \) for \( \z\in\RS \). Thus, if we can show that \( c_\rho \), defined in \eqref{veccrho}, is equal to \( \omega+\mathsf B\tau+2\mathsf B = 2\ic/\sqrt3 \), then we will get that 
\begin{equation}
\label{ex-7}
S_\rho\big(z^{(k)}\big) = S\big(z^{(k)}\big)D_{\widehat\rho}^{(-1)^k}\big(z^3\big),  \quad z\in\overline\C\setminus\Delta.
\end{equation}
Indeed, it is straightforward to see that this function is holomorphic and non-vanishing in \( \RS\setminus(\boldsymbol\Delta\cup\boldsymbol\alpha) \) and \( S_\rho(\z)S_\rho(\z^*) \equiv 1 \) there. Moreover, when \( \z\to\s\in\boldsymbol\Delta_+\setminus\boldsymbol A \) and \( \tr\to\s\in\boldsymbol\Delta_-\setminus\boldsymbol A \), it holds that \( z\to s\in\Delta^\circ_{\mathsf{sym}\pm} \) and \( t\to s\in\Delta^\circ_{\mathsf{sym}\mp} \). Therefore, with a slight abuse of notation, we can write
\[
S_{\rho+}(\s) = S_+(\s)D_{\widehat\rho\pm}\big(s^3\big)=S_-(\s)/\left(w_+(s)D_{\widehat\rho\mp}\big(s^3\big)\widehat\rho\big(s^3\big)\right) = S_{\rho-}(\s)/(\rho w_+)(s)
\]
for \( s\in\boldsymbol\Delta\setminus\boldsymbol A \). Since the jumps of \( S_\rho(\s) \) and of \( S(\z) \) across \( \boldsymbol\alpha\setminus\boldsymbol A \) coincide, the desired claim follows from the uniqueness part of Proposition~\ref{prop:szego}. Thus, to see that functions defined via \eqref{szego} and \eqref{ex-7} coincide, we only need to show that \( c_\rho=2\ic/\sqrt 3 \). Of course, this claim relies on the choice of the branch of \( \log(\rho w_+)(s) \). We shall assume that \eqref{szego} uses the following determination:
\[
\log(\rho w_+)(s) = \log\widehat\rho\big(s^3\big) + \log|w_+(s)| + a_w(s), \quad s\in\Delta^\circ,
\]
for some continuous determination of \( \log\widehat\rho(s) \), where \( a_w(s) \) was defined in \eqref{ex-5}. Then
\begin{multline*}
\frac1{2\pi\ic}\oint_{\boldsymbol\Delta} \log(\rho w_+)(s)\frac{\dd s}{w(\s)} = \frac1{\pi\ic}\int_{\Delta_\mathsf{sym}} a_w(s)\frac{\dd s}{w_+(s)} + \\ \frac{e^{-3\pi\ic/2}+ e^{7\pi\ic/6}+e^{-\pi\ic/6} }{\pi\ic}\int_{-1}^0\log\left(\widehat\rho\big(x^3\big)|w_+(x)|\right)\frac{\dd x}{|w(x)|},
\end{multline*}
where one can readily check that the constant in front of the last integral is equal to zero. Hence,
\begin{multline*}
\frac1{2\pi\ic}\oint_{\boldsymbol\Delta} \log(\rho w_+)(s)\frac{\dd s}{w(\s)} = \frac32\int_{\Delta_1}\frac{\dd s}{w_+(s)} - \frac76\int_{\Delta_2}\frac{\dd s}{w_+(s)} + \frac16\int_{\Delta_3}\frac{\dd s}{w_+(s)} \\ = \int_{-1}^0\left(\frac32e^{-3\pi\ic/2} - \frac76 e^{11\pi\ic/6} + \frac16e^{7\pi\ic/6} \right)\frac{\dd x}{|w(x)|} = \frac{4\eta}{\sqrt3}\int_{-1}^0\frac{\dd x}{|w(x)|}.
\end{multline*}
Thus, we get from \eqref{veccrho} and the first computation after \eqref{ex-5} that indeed \( c_\rho = 2\ic/\sqrt 3\).

The Jacobi inversion problem \eqref{jip} now reads
\[
\mathfrak a(\z_n) = \big[(n+1/2)(\omega + \mathsf B\tau)\big] = \big[(2n+1)(1-\mathsf B)/3\big].
\]
From this we can conclude that \( \z_0=\infty^{(0)} \), \( \z_1=\boldsymbol a_0=\boldsymbol 0 \), and \( \z_2=\infty^{(1)} \). Indeed, the middle claim follows immediately from \eqref{Amap} and the unique solvability of \eqref{jip}. The other two equalities follows from the fact that
\[
\mathfrak a\left(\infty^{(0)}\right) = \left[\frac{\omega+\mathsf B\tau}2\right] = \left[\frac13(1-\mathsf B)\right] \qandq \mathfrak a\left(\infty^{(1)}\right)= \left[-\frac{\omega+\mathsf B\tau}2\right] = \left[\frac23(1-\mathsf B)\right],
\]
which was shown in \eqref{0infty}. Thus, \( \N^\mathsf{reg} = \N\setminus(3\N) \).

\subsection{Sequence \( \N^\mathsf{sing} \)} 
\label{sap:sing}

It readily follows from \eqref{stokes} that \( b_i(\rho)=-2 \). It was shown in \cite{Kit95} that in this case \hyperref[rhphi]{\rhphi} with \( \alpha=0 \) is not solvable for \( x=0 \), that is,  \( \rho\in W^\mathsf{sing} \). More precisely, in \cite[Equations (18) and (19)]{Kit95} it was shown that the Stokes parameters\footnote{A slightly different Lax pair than the one presented in \eqref{pIIa}--\eqref{pIIc} is used in \cite{Kit95} so that the Stokes parameters appearing there are equal to \( -\ic s_1,\ic s_2,-\ic s_3 \).} \( s_1=s_2=s_3 =2\ic \) in \hyperref[rhpsi]{\rhpsi} with \( \nu=1/2 \) correspond to the solution \( q(s) \) of \eqref{p2} with a pole at the origin of residue \( -1 \). The desired conclusion then follows from the explanation given before \eqref{ant-u}. 
 
Due to the periodicity of \( \z_n \), we only need to compute \( \det(\boldsymbol T_1) \) and \( \det(\boldsymbol T_2) \). Since \( \z_1=\boldsymbol a_0 \), it holds that \( t_0^{(1)} =0 \) and therefore
\[
\det(\boldsymbol T_1) = t_0^{(0)}\left( 2t_2^{(1)}-h_1t_1^{(1)} \right) \qandq \det(\boldsymbol T_2) = t_0^{(2)}\left( 2t_2^{(1)}-h_1t_1^{(1)} \right).
\]
Since \( t_0^{(0)}t_0^{(2)} \neq 0 \), to show that \( \N^\mathsf{sing}=\varnothing \), we need to prove that \(  2t_2^{(1)} =h_1t_1^{(1)} \). To this end, set
\[
A(\z) := \exp \big\{-2\pi\ic(\tau+2) a(\z) \big\}
\]
and, as usual, we shall write \( A(z) \) for the pull-back of \( A(\z) \) from \( \RS^{(0)} \). It now follows from \eqref{GandH}, \eqref{ri}, \eqref{ex-9}, and \eqref{ex-7} that
\[
H^{-1}(z) = \ic (S_\rho r)^2(z)w(z) = (A\Phi)^2(z) \big(\widehat\rho D_{\widehat\rho}^2\big)\big(z^3\big).
\]
Let \( \boldsymbol U_r \) be two copies of \( \{|z|<r\}\setminus[-r,0] \), \( r<1 \), glued crosswise across the cuts (double cyclic neighborhood of \( 0 \)). It is quite easy to see that the function \( \big(\widehat\rho D_{\widehat\rho}^2\big)(z) \) lifted to one of the copies of \( \{|z|<r\}\setminus[-r,0] \) extends analytically to \( \boldsymbol U_r \) by lifting its reciprocal to the other copy. Hence, \( \big(\widehat\rho D_{\widehat\rho}^2\big)\big(z^3\big) \) possesses a Puiseux series at zero and satisfies there
\[
\big(\widehat\rho D_{\widehat\rho}^2\big)\big(z^3\big) = 1 + \mathcal O\big(z^{3/2}\big) \qasq z\to 0.
\]
Recall that the same is true about \( \Phi(z) \), see \eqref{prop5-5}. Thus, we can deduce from \eqref{H-exp} that if we develop \( A(z) \) into a Puiseux series in the sector \( S_2 \), then
\[
A(z) = 1 - (h_1/2) z^{1/2} + \mathcal O(z) \qasq z\to 0.
\]
On the other hand, it follows from \eqref{Phi-jump}, \eqref{T-jump}, \eqref{a-jump}, and the fact that \( c_\rho=\omega+\mathsf B\tau +2\mathsf B \) that \( K(\z) := \big(AT_1\Phi^2\big)(\z) \) is a rational function on \( \RS \) with the zero/pole divisor \( \boldsymbol a_0 + \infty^{(1)} - 2\infty^{(0)}\). The Puiseux series of \( K(z) \) in the sector \( S_2 \) is equal to
\begin{equation}
\label{ex-10}
K(z) = t_1^{(1)} z^{1/2} + \left( t_2^{(1)} - h_1t_1^{(1)}/2 \right) z + \mathcal O\left(z^{3/2}\right) \qasq z\to 0.
\end{equation}
It follows from the symmetries of \( \Delta_\mathsf{sym} \) that \( K(\eta\z) \) and \( K(\overline\eta\z) \), \( \eta=e^{2\pi\ic/3} \), are rational functions on \( \RS \) as well. Moreover, they have the same zero/pole divisors as \( K(\z) \). Thus, they must be constant multiples of \( K(\z) \). By looking at their Puiseux series in \( S_2 \) we then can conclude that the coefficient next \( z \) in \eqref{ex-10} must be equal to \( 0 \), which is exactly the desired claim.

\section{Examples of Approximated Functions}
\label{ap:functions}

{
In this appendix we show how to write functions \eqref{log-power} in the form \eqref{CI}. We also compute their Stokes parameters \eqref{stokes} and show that Theorem~\ref{thm:asymp2} is always applicable to these functions.

Logarithmic functions in \eqref{log-power} are the easiest example. Let \( c_1,c_2,c_3 \) be arbitrary non-zero complex numbers. Set \( c_0 := -c_1-c_2-c_3 \) and define \( \rho_i(s) = \rho_{\Delta_i^\circ}(s)=2\pi\ic c_i \), \( i\in\{1,2,3\} \) (unless \( c_1=c_2=c_3 \), \( \rho(s) \) is not well defined at \( a_0 \)). Then it holds that
\[
f(z) = \frac1{2\pi \ic}\int_\Delta \frac{\rho(s)ds}{s-z}  = \sum_{i=1}^3 c_i \log (z-s) \bigg|_{a_i}^{a_0} = \sum_{k=0}^3 c_k\log (z-a_k),
\]
where the branches of the logarithms are chosen so that the right-hand side above is analytic in \( \overline\C\setminus\Delta \). Since each function \( \rho_i(s) \) is constant on \( \Delta_i \), condition \eqref{bad-cond} is satisfied automatically. In fact, it holds that \( \kappa_{i,j}(z) \equiv 1\) and therefore \hyperref[rhpo]{\rhpo} is not an approximate but an exact parametrix for these functions. In particular, estimate \eqref{R-jump-5}  and an analogous estimate in the proof of Theorem~\ref{thm:asymp2} are not needed. This class of functions does include weights in both \( W^\mathsf{reg} \) and \(W^\mathsf{sing} \). Indeed, it holds that
\[
b_1(\rho) = -\frac{c_2+c_3}{c_1}, \quad b_2(\rho) = -\frac{c_1+c_3}{c_2}, \qandq b_3(\rho) = -\frac{c_1+c_2}{c_3}.
\]
Therefore, if \( -c_0=c_1+c_2+c_3 =0 \), then \( b_i(\rho)=1 \) and \( \rho\in W^\mathsf{reg} \) (\( f(z) \) has no logarithmic singularity at \( a_0 \) and such functions are also covered by the results of \cite{ApYa15}). On the other hand, if \( c_1=c_2=c_3 \), then \( b_i(\rho)=-2 \) and, as mentioned at the beginning of Section~\ref{sap:sing},  \( \rho\in W^\mathsf{sing} \). 

Let \( \alpha_k>-1 \), \( k\in\{0,1,2,3\} \), be such that their sum is an integer necessarily bigger or equal to \( -2 \). Recall our notation from Section~\ref{ss:ol}, where we denoted by \( S_1,S_2,S_3 \) the connected components of \( U_0\setminus \Delta \) labeled in such a way that \( S_k \) does  not border \( \Delta_k^\circ \). Let \( L \) be any curve extending to infinity from \( a_0 \) that partially lies within \( S_1 \). We take branches \( (z-a_i)^{\alpha_i} \) that are analytic off \( \Delta_i\cup L \), \( i\in\{1,2,3 \} \), and \( (z-a_0)^{\alpha_0} \) that is analytic off \( L \). Define
\[
f(z) = \prod_{k=0}^3(z-a_k)^{\alpha_k} - p(z),
\]
where \( p(z) \) is the polynomial part of the product at infinity. That is, \( f(z) \) is analytic in \( \overline\C\setminus\Delta \) and vanishes at infinity. Let \( \rho(s) = f_+(s) - f_-(s) \), which is a smooth function on \( \Delta\setminus\{a_0,a_1,a_2,a_3\} \). Let \( C_\rho(z) \) be the right-hand side of \eqref{CI}. Then it follows from Plemelj-Sokhotski formulae that
\[
C_{\rho+}(s) - f_+(s) = C_{\rho-}(s) + \rho(s) - f_+(s) = C_{\rho-}(s) - f_-(s), \quad s\in\Delta\setminus\{a_0,a_1,a_2,a_3\}.
\]
Since \( C_\rho(z)-f(z) \) has smooth traces on each \( \Delta_i^\circ \), this function is analytic across each \( \Delta_i^\circ \) by Morera's theorem. As this difference vanishes at infinity and cannot have polar or essential singularities at \( \{a_0,a_1,a_2,a_3\} \) (this follows from the known behavior of the Cauchy integrals around the endpoints of the contours of integration), it must holds that \( f(z) \equiv C_\rho(z) \). Let
\[
\varrho(z) = \prod_{i=1}^3(z-a_i)^{\alpha_i} \begin{cases}
1, & z\in S_2, \\
e^{-2\pi\ic\alpha_1}, & z\in S_3,
\end{cases}
\quad 
\varrho(z) = \prod_{i=1}^3(z-a_i)^{\alpha_i} \begin{cases}
e^{-2\pi\ic(\alpha_1+\alpha_2)}, & z\in S_{12}, \\
e^{2\pi\ic\alpha_3}, & z\in S_{13},
\end{cases}
\]
where \( S_{1j} \) is the connected component of \( S_1\setminus L \) that borders \( \Delta_j \). It can be readily seen that \( \varrho(z) \) is in fact analytic in \( U_0 \). Recall \eqref{varrho}. Then it follows from the first set of definitions above that
\[
\varrho_1(s) = \frac{f_+(s)-f_-(s)}{(s-a_0)^{\alpha_0}} = \rho(s) \big(1-e^{2\pi\ic\alpha_1}\big), \quad s\in\Delta_1^\circ.
\]
Similarly, we get that
\[
\varrho_2(s) = \varrho(s) e^{2\pi\ic\alpha_1}\big(1-e^{2\pi\ic\alpha_2} \big) \qandq \varrho_3(s) = \varrho(s) \big(e^{-2\pi\ic\alpha_3} -1\big).
\]
Hence, the weights \( \varrho_i(s) \) are constant multiples of the same function analytic in \( U_0 \) and the condition \eqref{bad-cond} is satisfied (in fact, all the functions \( \kappa_{i,j}(z)\equiv 1 \) and \hyperref[rhpo]{\rhpo} is an exact parametrix). Moreover, by plugging the above expressions into \eqref{stokes} we get that
\[
b_i(\rho) = \frac{\sin(\alpha_0+\alpha_i)\pi}{\sin \alpha_i\pi} = \cos \alpha_0\pi + \sin \alpha_0\pi \cot \alpha_i\pi
\]
for each \( i\in\{1,2,3\} \). Hence, if \( \alpha_0 = 0 \), then \( b_i(\rho)=1 \) and \( \rho\in W^\mathsf{reg} \) (\( f(z) \) has no  singularity at \( a_0 \) and such functions are also covered by the results of \cite{ApYa15}). On the other hand, for any other \( \alpha_0\neq 0 \) fixed, every real solution of \eqref{bi} can be expressed in the above form with an appropriate choice of \( \alpha_1,\alpha_2,\alpha_3\in(-1,0) \). It is an interesting questions whether the class of real Stokes parameters contains some that lead to the matrices \( \boldsymbol\Psi_{\alpha_0} \) that have a pole at \( x=0 \).
}

\end{document}